\documentclass[a4paper,reqno, 10pt]{amsart}
      \usepackage{etex}
      \usepackage[OT2, T1]{fontenc}
      \usepackage{url}
      \usepackage{amsmath}
      \usepackage{lmodern}
      \usepackage{array}
      \usepackage{graphicx}
      \usepackage{amsfonts}
      \usepackage{amssymb}
      \usepackage{extpfeil}
\usepackage{xargs}
      \usepackage{xcolor}
      \definecolor{imperialred}{RGB}{237, 41, 57}
      \definecolor{royalblue}{RGB}{64, 106, 212}
      \definecolor{link}{RGB}{11,0,128}
      \definecolor{olivegreen}{RGB}{128, 128, 0}
      \usepackage[colorlinks=true, citecolor=royalblue,linkcolor=olivegreen,urlcolor=link]{hyperref}

      \usepackage{amsthm}
      \usepackage{xkeyval}
      \usepackage{colonequals}		
      \usepackage{mathabx}		
      \usepackage[left]{lineno}
      \usepackage{amscd}
         \usepackage{parskip}
          \usepackage{ragged2e}
           \usepackage{graphicx}
      \usepackage{mathrsfs} 			
            \usepackage{cleveref}
      \usepackage{enumitem}
      \usepackage{tikz-cd}
      \usepackage{tikz}
        \usetikzlibrary{shapes.geometric}
      \usetikzlibrary{automata, bending,matrix,arrows,decorations.pathmorphing,backgrounds,positioning,fit,petri,arrows.meta}
      \usetikzlibrary{arrows}
      \usepackage{enumitem}
      \tikzset{commutative diagrams/.cd,arrow style=tikz,diagrams={>=latex'}}
      \usepackage{filecontents}
      \usepackage[alphabetic,lite]{amsrefs} 	
      \usepackage{stmaryrd}	
      \usetikzlibrary{arrows}

\usepackage{marginnote}
\newcommand{\mytodo}[2][]{{%
 \let\marginpar\marginnote
 \reversemarginpar
 \renewcommand{\baselinestretch}{0.8}%
 \todo[#1]{#2}}}

         \newcommand{\gL}{\lambda}

         \newcommand{\GG}{\Gamma}

         \newcommand{\bA}{\mathbb{A}}

         \newcommand{\bG}{\mathbb{G}}

         \newcommand{\cH}{\mathcal{H}}

         \newcommand{\cM}{\mathcal{M}}

         \newcommand{\cP}{\mathcal{P}}
         \newcommand{\cQ}{\mathcal{Q}}
         \newcommand{\cR}{\mathcal{R}}

         \newcommand{\fm}{\mathfrak{m}}
         
         \newcommand{\fp}{\mathfrak{p}}
         \newcommand{\fq}{\mathfrak{q}}

         \newcommand{\fS}{\mathfrak{S}}

         \newcommand{\sA}{\mathscr{A}}
         \newcommand{\sB}{\mathscr{B}}

         \newcommand{\sE}{\mathscr{E}}
         \newcommand{\sF}{\mathscr{F}}
         \newcommand{\sG}{\mathscr{G}}
         \newcommand{\sH}{\mathscr{H}}
         \newcommand{\sI}{\mathscr{I}}
         
         \newcommand{\sK}{\mathscr{K}}
         \newcommand{\sL}{\mathscr{L}}
         \newcommand{\sM}{\mathscr{M}}
         \newcommand{\sN}{\mathscr{N}}
         \newcommand{\sO}{\mathscr{O}}

         \newcommand{\sV}{\mathscr{V}}

         \newcommand{\ra}{\rightarrow}

         \newcommand{\hra}{\hookrightarrow}

         \newcommand{\wt}{\widetilde}
         \newcommand{\wh}{\widehat}

         \newcommand{\pr}{^{\prime}}
         \newcommand{\prpr}{^{\prime\prime}}

         \newcommand{\ce}{\colonequals}
         \newcommand{\ov}{\overline}
         \newcommand{\un}{\underline}

         \renewcommand{\b}{\textbf}
         \newcommand{\surjects}{\twoheadrightarrow}
         
         \newcommand{\isoto}{\overset{\sim}{\longrightarrow}}
         
         \newcommand{\fppf}{\mathrm{fppf}}		                                           
         \newcommand{\et}{\mathrm{\acute{e}t}}	                                           
         \newcommand{\Zar}{\mathrm{Zar}}		                                               
         \newcommand{\sh}{\mathrm{sh}}		                                              

         \newcommand{\specializes}{\rightsquigarrow}                                               
         
         \newcommand{\ddual}{\vee\!\vee}

         \providecommand{\abs}[1]{\left\lvert#1\right\rvert}

         \providecommand{\p}[1]{\left(#1\right)}

         \providecommand{\SP}[1]{\cite[\href{https://stacks.math.columbia.edu/tag/#1}{#1}]{SP}}

         \providecommand{\f}[2]{\frac{#1}{#2}}

\newextarrow{\xbigtoto}{{15}{15}{15}{12}}
   {\bigRelbar\bigRelbar{\bigtwoarrowsleft\rightarrow\rightarrow}}
         
         \DeclareMathOperator{\Ker}{Ker}			                       
         \DeclareMathOperator{\Coker}{Coker}		                       
         \DeclareMathOperator{\im}{Im}			                       
         \DeclareMathOperator{\Spec}{Spec}		                       
         \DeclareMathOperator{\wdim}{wdim}
         \DeclareMathOperator{\fldim}{fl.dim}
         \DeclareMathOperator{\pd}{proj.dim}
         \DeclareMathOperator{\rad}{rad}			                       
         \DeclareMathOperator{\Hom}{Hom}			                       
         \DeclareMathOperator{\Ann}{Ann}			                       
         \DeclareMathOperator{\Ass}{Ass}			                       
         \DeclareMathOperator{\Frac}{Frac}		                       
         \DeclareMathOperator{\depth}{depth}		                       
         \DeclareMathOperator{\Supp}{Supp}		                       
         \DeclareMathOperator{\Ext}{Ext}		                           	
         \DeclareMathOperator{\Tor}{Tor}		                           	
         \DeclareMathOperator{\Set}{\textbf{Set}}
         \DeclareMathOperator{\Sch}{\textbf{Sch}}		                                                  
         \DeclareMathOperator{\GL}{GL}		                                                  
         \DeclareMathOperator{\End}{End}		                                                  
         \DeclareMathOperator{\Pic}{Pic}		                                                  
         \DeclareMathOperator{\codim}{codim}		                                                  
\makeatletter
\newcommand{\csub[1]}{\centering\subsection{\text{#1} \nopunct} \hfill \justify}
\makeatother
         \newcommand{\ba}{\begin{aligned}}
         \newcommand{\ea}{\end{aligned}}
         \newcommand{\be}{\begin{equation}}
         \newcommand{\ee}{\end{equation}}
         \newcommand{\pf}{\begin{proof}}
         \newcommand{\bpf}{\begin{proof}}
         \newcommand{\epf}{\end{proof}}
         \newcommand{\bsol}{\begin{solution}}
         \newcommand{\esol}{\end{solution}}
         \newcommand{\bthm}{\begin{thm}}
         \newcommand{\ethm}{\end{thm}}
         \newcommand{\bthmt}{\begin{thm-tweak}}
         \newcommand{\ethmt}{\end{thm-tweak}}
         \newcommand{\bprop}{\begin{prop}}
         \newcommand{\eprop}{\end{prop}}
         \newcommand{\bcor}{\begin{cor}}
         \newcommand{\ecor}{\end{cor}}
         \newcommand{\bcort}{\begin{cor-tweak}}
         \newcommand{\ecort}{\end{cor-tweak}}
         \newcommand{\brem}{\begin{rem}}
         \newcommand{\erem}{\end{rem}}
         \newcommand{\bremt}{\begin{rem-tweak}}
         \newcommand{\eremt}{\end{rem-tweak}}
         \newcommand{\brems}{\begin{rems} \hfill \begin{enumerate}[label=\b{\thenumberingbase.},ref=\thenumberingbase]}
         
         \newcommand{\erems}{\end{enumerate} \end{rems}}
         \newcommand{\begs}{\begin{egs} \hfill \begin{enumerate}[label=\b{\thenumberingbase.},ref=\thenumberingbase]}
         \newcommand{\egi}{\addtocounter{numberingbase}{1} \item}
         \newcommand{\eegs}{\end{enumerate} \end{egs}}
         \newcommand{\eremstweak}{\end{enumerate} \end{rems-tweak}}
         \newcommand{\eremst}{\end{enumerate} \end{rems-tweak}}
         \newcommand{\blem}{\begin{lemma}}
         \newcommand{\elem}{\end{lemma}}
         \newcommand{\blemt}{\begin{lemma-tweak}}
         \newcommand{\elemt}{\end{lemma-tweak}}
         \newcommand{\bconj}{\begin{conj}}
         \newcommand{\econj}{\end{conj}}
         \newcommand{\bprob}{\begin{Problem}}
         \newcommand{\eprob}{\end{Problem}}

         \newcommand{\bq}{\begin{Q}}
         \newcommand{\eq}{\end{Q}}
         \newcommand{\benum}{\begin{enumerate}[label={{\upshape(\alph*)}}]}
         \newcommand{\benuma}{\begin{enumerate}[label={{\upshape(\arabic*)}}]}
         \newcommand{\benumr}{\begin{enumerate}[label={{\upshape(\roman*)}}]}
         \newcommand{\eenum}{\end{enumerate}}
         \newcommand{\bc}{}
         \newcommand{\bd}{\begin{defn}}
         \newcommand{\ed}{\end{defn}}
         \newcommand{\bque}{\begin{que}}
         \newcommand{\eque}{\end{que}}
         \newcommand{\bfct}{\begin{fact}}
         \newcommand{\efct}{\end{fact}}
         \newcommand{\bdt}{\begin{defn-tweak}}
         \newcommand{\edt}{\end{defn-tweak}}
         \newcommand{\beg}{\begin{eg}}
         \newcommand{\eeg}{\end{eg}}
         \newcommand{\begt}{\begin{eg-tweak}}
         \newcommand{\eegt}{\end{eg-tweak}}
         \newcommand{\bcl}{\begin{claim}}
         \newcommand{\ecl}{\end{claim}}

         \newcommand{\x}{\text}

         \newcommand{\q}{\quad}
         \newcommand{\qq}{\quad\quad}

         \newcommand{\tst}{\textstyle}

         \newcommand{\sHom}{\mathscr{H}\! om}
         
\providecommand{\SPD}[2]{\cite[\href{https://stacks.math.columbia.edu/tag/#1}{#1}, \href{https://stacks.math.columbia.edu/tag/#2}{#2}]{SP}}
     \providecommand{\SPN}[9]{\cite[\href{https://stacks.math.columbia.edu/tag/#1}{#1}, \href{https://stacks.math.columbia.edu/tag/#2}{#2},
     \href{https://stacks.math.columbia.edu/tag/#3}{#3},
     \href{https://stacks.math.columbia.edu/tag/#4}{#4},
     \href{https://stacks.math.columbia.edu/tag/#5}{#5},
     \href{https://stacks.math.columbia.edu/tag/#6}{#6},
     \href{https://stacks.math.columbia.edu/tag/#7}{#7},
     \href{https://stacks.math.columbia.edu/tag/#8}{#8},
     \href{https://stacks.math.columbia.edu/tag/#9}{#9}]{SP}}

\tikzset{
    labl/.style={anchor=south, rotate=90, inner sep=.5mm}
}

      \usepackage{aliascnt}
\newaliascnt{numberingbase}{subsubsection}
\numberwithin{equation}{numberingbase}

\newtheoremstyle{thms}{0pt}{0pt}{\itshape}{}{\bfseries}{.}{ }{}
\theoremstyle{thms}
\newtheorem{conj}[numberingbase]{Conjecture}

\newtheorem{cor}[numberingbase]{Corollary}
\newtheorem{lemma}[numberingbase]{Lemma}

\newtheorem{prop}[numberingbase]{Proposition}

\newtheorem{Q}[numberingbase]{Question}
\newtheorem{thm}[numberingbase]{Theorem}

\newtheoremstyle{claims}{0pt}{0pt}{}{}{\itshape}{.}{ }{}
\theoremstyle{claims}
\newtheorem{claim}[equation]{Claim}
\newtheorem{cl-tweak}[subsubsection]{Claim}
\Crefname{cl-tweak}{Claim}{Claims}

\newtheoremstyle{defs}{0pt}{0pt}{}{}{\bfseries}{.}{ }{}
\theoremstyle{defs}
\newtheorem{defn}[numberingbase]{Definition}

\newtheorem{eg}[numberingbase]{Example}

\newtheorem*{egs}{Examples}
\newtheorem{rem}[numberingbase]{Remark}

\newtheorem*{rems}{Remarks}

\Crefname{claim}{Claim}{Claims}
\Crefname{sublemma}{Lemma}{Lemmas}
\Crefname{conj}{Conjecture}{Conjectures}
\Crefname{cor}{Corollary}{Corollaries}
\Crefname{defn}{Definition}{Definitions}
\Crefname{eg}{Example}{Examples}
\Crefname{prop}{Proposition}{Propositions}
\Crefname{Q}{Question}{Questions}
\Crefname{rem}{Remark}{Remarks}
\Crefname{thm}{Theorem}{Theorems}
\Crefname{variant}{Variant}{Variants}

\theoremstyle{thms}
\newtheorem{thm-tweak}[subsection]{Theorem}
\Crefname{thm-tweak}{Theorem}{Theorems}
\newtheorem{lemma-tweak}[subsection]{Lemma}
\Crefname{lemma-tweak}{Lemma}{Lemmas}
\newtheorem{cor-tweak}[subsection]{Corollary}
\Crefname{cor-tweak}{Corollary}{Corollaries}
\newtheorem{prop-tweak}[subsection]{Proposition}
\Crefname{prop-tweak}{Proposition}{Propositions}
\newtheorem{conj-tweak}[subsection]{Conjecture}
\Crefname{conj-tweak}{Conjecture}{Conjectures}

\theoremstyle{defs}
\newtheorem{defn-tweak}[subsection]{Definition}
\Crefname{defn-tweak}{Definition}{Definitions}
\newtheorem{eg-tweak}[subsection]{Example}
\Crefname{eg-tweak}{Example}{Examples}
\newtheorem*{rems-tweak}{Remarks}
\newtheorem{rem-tweak}[subsection]{Remark}
\Crefname{rem-tweak}{Remark}{Remarks}

\newtheoremstyle{subsection-tweak}
   {0pt}
   {0pt}%
   {}
   {}%
   {\bfseries}
   {}%
   {.5em}
   {\thmnumber{\@{#1}{}\@{#2}.}%
    \thmnote{~{\bfseries#3.}}}

\theoremstyle{subsection-tweak}
\newtheorem{pp}[subsection]{}
\newcommand{\bpp}{\begin{pp}}
\newcommand{\epp}{\end{pp}}
\newtheorem{pp-t}[subsubsection]{}
\newcommand{\bppt}{\begin{pp-t}}
\newcommand{\eppt}{\end{pp-t}}

\theoremstyle{subsection-tweak}

\theoremstyle{subsection-tweak}
\newtheorem{pp-tweak}{}

      \makeatletter
      \def\@tocline#1#2#3#4#5#6#7{
          \begingroup
          \@ifempty{#4}{}{}

          \parindent\z@ \leftskip#3\relax \advance\leftskip\@tempdima\relax
          #5\hskip-\@tempdima
            \ifcase #1
             \or\or \hskip 2em \or \hskip 1em \else \hskip 3em \fi%
            #6\nobreak\relax
          \dotfill\hbox to\@pnumwidth{\@tocpagenum{#7}}\par
          \nobreak
          \endgroup
        }
       \def\l@section{\@tocline{1}{0pt}{1pc}{}{}}

      \renewcommand{\tocsection}[3]{%
        \indentlabel{\@ifnotempty{#2}{\makebox[1.3em][l]{%
          \ignorespaces#1 \bfseries{#2}.\hfill}}}\bfseries{#3}
          \vspace{-3.5pt}}

      \renewcommand{\tocsubsection}[3]{%
        \indentlabel{\@ifnotempty{#2}{\hspace*{-0.5em}\makebox[2.1em][l]{%
          \ignorespaces#1#2.\hfill}}}#3
          \vspace{-4.5pt}}

\makeatother

\makeatletter 
\newcommand\appendix@section[1]{%
  \refstepcounter{section}%
  \orig@section*{Appendix \@Alph\c@section. #1}%
}
\let\orig@section\section
\g@addto@macro\appendix{\let\section\appendix@section}
\makeatother

      \DeclareMathVersion{normal2}

\author{Ning Guo}
\address{St. Petersburg branch of V. A. Steklov Mathematical Institute, Fontanka 27, 191023 St. Petersburg, Russia}
\email{guo.ning@eimi.ru}
\author{Fei Liu}
\address{Department of Mathematics, Southern University of Science and Technology, Shenzhen, China}
\email{liufei54@pku.edu.cn}
\date{\today}

\def\UTFviii@defined#1{%
  \ifx#1\relax
      \PackageError{inputenc}{Unicode\space char\space\expandafter
                              \UTFviii@splitcsname\string#1\relax
                              \MessageBreak
                              not\space set\space up\space
                              for\space use\space with\space LaTeX}\@eha
  \else\expandafter
    #1%
  \fi
}

\makeatletter
\def\UTFviii@defined#1{%
  \ifx#1\relax
      ?%
  \else\expandafter
    #1%
  \fi
}

\makeatother


\begin{document}
\textheight = 21cm
\textwidth = 13cm

\title{Purity and quasi-split torsors over Pr\"ufer bases}

\begin{abstract}
\setlength{\parindent}{0pt}
We establish an analogue of the Zariski--Nagata purity theorem for finite \'etale covers on smooth schemes over Pr\"ufer rings by demonstrating Auslander's flatness criterion in this non-Noetherian context.
We derive an Auslander--Buchsbaum formula for general local rings, which provides a useful tool for studying the algebraic structures involved in our work.
Through the analysis of reflexive sheaves, we prove various purity theorems for torsors under certain group algebraic spaces, such as the reductive ones.
Specifically, using results from \cite{EGAIV4} on parafactoriality on smooth schemes over normal bases, we prove the purity for cohomology groups of multiplicative type groups at this level of  generality.
Subsequently, we leverage the aforementioned purity results to resolve the Grothendieck--Serre conjecture for torsors under a quasi-split reductive group scheme over schemes smooth over Pr\"ufer rings.
Along the way, we also prove a version of the Nisnevich purity conjecture for quasi-split reductive group schemes in our Pr\"uferian context, inspired by the recent work of $\check{\mathrm{C}}$esnavi$\check{\mathrm{c}}$ius \cite{Ces22b}.
\par
\medskip
\textbf{Résumé.} (Pureté et torseurs quasi-déployés sur les bases de Prüfer) Nous établissons un analogue du théorème de pureté de Zariski--Nagata pour les revêtements étales sur les schémas lisses sur les anneaux de Prüfer en démontrant le critère de platitude d'Auslander dans ce contexte non-noethérien. Nous dérivons une formule d'Auslander--Buchsbaum pour les anneaux locaux généraux, qui fournit un outil utile pour étudier les structures algébriques impliquées dans notre travail. Grâce à l'analyse des faisceaux réflexives, nous prouvons divers théorèmes de pureté pour les torseurs sous certains espaces algébriques en groupes, tels que les réductifs. En particulier, en utilisant des résultats de \cite{EGAIV4} sur la parafactorialité sur les schémas lisses sur des bases normales, nous prouvons la pureté pour les groupes de cohomologie des groupes de type multiplicatif à ce niveau de généralité. Par la suite, nous utilisons les résultats de pureté susmentionnés pour résoudre la conjecture de Grothendieck--Serre pour les torseurs sous un schéma en groupes réductifs quasi-déployés sur des schémas lisses sur des anneaux de Prüfer. Nous prouvons également une version de la conjecture de pureté de Nisnevich pour les schémas en groupes réductifs quasi-déployés dans notre contexte Prüferien, inspirée par les travaux récents de $\check{\mathrm{C}}$esnavi$\check{\mathrm{c}}$ius \cite{Ces22b}.
\end{abstract}

\maketitle 
\hypersetup{
    linktoc=page,     
}
\newpage
\renewcommand*\contentsname{}
\q\\

\tableofcontents

\section{{Introduction}}
\bpp[Purity and regularity]
In algebraic geometry, purity refers to a diverse range of phenomena in which certain invariants or categories associated to geometric objects are insensitive to the removal of closed subsets of large codimensions.
In the classical Noetherian world, purities, say, for vector bundles (and even torsors), or for finite \'etale covers, are intimately related to the regularities  measured by lengths of regular sequences of geometric objects.
For a concrete instance, the Auslander--Buchsbaum formula \cite{AB57}*{Theorem~3.7}
\[
   \depth_RM+\pd_RM=\depth_RR
\]
controls the projective dimension of the finite type module $M$ over the Noetherian local ring $R$ via depths, leading to the purity for vector bundles on regular local rings of dimension $2$ (\cite{Sam64}*{Proposition~2}).
Granted this, Colliot-Thélène and Sansuc \cite{CTS79}*{Théorème~6.13} established the purity for reductive torsors over arbitrary regular local ring $R$ of dimension $2$ by bootstrapping from the vector bundle case: for every reductive $R$-group scheme $G$, the restriction map
\[
 \text{$H^1_{\et}(\Spec R,G)\isoto H^1_{\et}(\Spec R\backslash \{\fm_R\},G)$\q is bijective.}
\]
Nevertheless, not only does the term `regularity' make sense for Noetherian rings, its non-Noetherian generalization can still enlighten us to contemplate purity problems.
\epp
\bpp[Regularity and its avatar over Pr\"ufer bases]\label{par-coh-reg}
The concept of regularity for non-Noetherian rings was first introduced by Bertin, as found in references \cite{Ber71} and \cite{Ber72}*{Définition~3.5}, specifically for coherent local rings. A ring
$R$ is termed \emph{regular} if every finitely generated ideal of $R$ possesses a finite projective dimension.
According to Serre's homological characterization \cite{Ser56}*{Théorème~3}, this definition aligns with the traditional understanding of regularity in the context of Noetherian rings.
A typical non-Noetherian example can be found in \emph{Prüfer domains}. By definition, these are domains whose all local rings as valuation rings.
Recall that an integral domain $V$ is a \emph{valuation ring} if every pair $a,b\in V\backslash\{0\}$ satisfies either $a\in (b)$ or $b\in (a)$.
As a case in point, Noetherian valuation rings are precisely {either} discrete valuation rings {or fields}.
The regularity of Prüfer domains is a direct consequence of the characteristic that all finitely generated ideals of valuation rings are principal\footnote{To elucidate, given a Prüfer domain
$R$ (which is coherent), any partial resolution
$0\to N \to R^n \to I \to 0$ of a finitely generated ideal $I\subset R$ results in $N$ being finitely presented. As $N$ is finite free over each local ring of $R$, it is consequently finite projective over
$R$.}.
{As a further instance, every smooth algebra over a Prüfer domain is coherent regular (\Cref{coh of V-flat ft alg}).}
For recent research on the regularity of schemes over  Pr\"ufer domains, see \cite{Kna03}, \cite{Kna08}.
In addition to their regularity and other properties (see, e.g., \Cref{geom}), Pr\"ufer domains are ubiquitous in the study of nonarchimedean geometry, Zariski--Riemann spaces, and other fields, which motivates further investigation of their algebro-geometric properties.
\epp
\bpp[Auslander--Buchsbaum for general local rings]
Let $A$ be a local ring with quasi-compact punctured spectrum and $M$ an $A$-module having a finite resolution by finite free $A$-modules. 
We have the following Auslander--Buchsbaum formula
\[
     \mathrm{proj.}\dim_A (M)+\depth_A(M)=\depth_A(A).
     \]
Here $\mathrm{proj.}\dim_A(0)=-\infty$ and $\depth_AM$ is the smallest $i$ such that the $i$-th local cohomology of $M$ is nonzero (\S\ref{setup of depth & proj dim}).
Our proof is significantly different from the classical case \cite{AB57}*{Theorem~3.7}.
Specifically, we bypass the interpretation of projective dimensions in terms of Tor functors, which is a crucial ingredient in Auslander--Buchsbaum's argument.
\epp
\bpp[Basic setup I]\label{intro-setup}
The purity part of the present article focuses on a semilocal Pr\"ufer domain $R$ with $\dim R>0$ (and with $\dim R<\infty$ if necessary), an $R$-flat finite type algebraic space $X$ with regular $R$-fibres, and a closed subspace $Z\subset X$ such that $j\colon X\backslash Z\hra X$ is quasi-compact.
For a point $x\in X$ lying in an open subscheme, the local ring of $X$ at $x$ makes sense and we denote $A\ce \sO_{X,x}$.
When involving torsors on $X$, we let $G$ be an $X$-group algebraic space that \'etale locally permits an embedding $G\hra \GL_n$ such that $\GL_n/G$ is $X$-affine.
This condition is fulfilled if $G$ is $X$-reductive\footnote{Namely, it is a smooth affine $X$-group algebraic space $G$ whose geometric $X$-fibres are (connected) reductive algebraic groups.}, or finite locally free.
\epp
\bpp[Purity for torsors on smooth relative curves over Pr\"ufer bases]\label{par-Purity for torsors on smooth relative curves over Pr\"ufer bases}
Once the projective dimensions of reflexive sheaves on $X$ are controlled, by imposing codimensional constraints on $Z$, we may extend vector bundles on $X\backslash Z$ to $X$, as in Noetherian scenarios.
Subsequently, this allows us to obtain the purity \Cref{purity for rel. dim 1} for $G$-torsors {(where $G$ is as in the Basic setup I)}: if $X$ is an $S$-curve and $Z$ satisfies
\[
\text{$Z_{\eta}= \emptyset$\q for each generic point $\eta\in S$ \q and\q  $\codim(Z_s,X_s)\ge 1$ for all $s\in S$,}
\]
then restriction induces the following equivalence of categories of $G$-torsors
\[
\mathbf{Tors}(X_{\fppf},G) \isoto \mathbf{Tors}((X\backslash Z)_{\fppf},G).
\]
In particular, passing to isomorphism classes of objects, we have the bijection 
\[
H^1_{\fppf}(X,G)\simeq H^1_{\fppf}(X\backslash Z,G)
\]
of nonabelian cohomology pointed sets.
Meanwhile,  a local version \Cref{extends across codim-2 points} allows us to loose  constraints on the relative dimension of $X$: if either
\begin{align*}
 x\in X_{\eta}  & \quad  \text{with $\dim \sO_{X_{\eta},x} =2$, or }\\ 
 x\in X_s  & \quad \text{with $s \neq \eta$ and $\dim \sO_{X_s,x} =1$,}
\end{align*}
then every $G$-torsor over $\Spec \sO_{X,x} \backslash \{x\}$ extends uniquely to a $G$-torsor over $ \sO_{X,x}$. This permits us to iteratively extend reductive torsors beyond a closed subset of higher fibrewise codimensions.
\epp
\bpp[Zariski--Nagata over Pr\"ufer bases]\label{par-Zariski--Nagata over Pr\"ufer bases}
The Zariski--Nagata purity, known as ``purity of branch locus'', states that every finite extension $A\subset B$ of rings with $A$ regular Noetherian and $B$ normal is unramified if and only if so it is in codimension one on $\Spec B$.
This purity was settled by Zariski \cite{Zar58} in a geometric context, and more algebraically by Nagata \cite{Nag59}  based on Chow's local Bertini theorem.
In contrast to them, Auslander gave an alternative proof \cite{Aus62}*{Theorem~1.4} by skillful homological methods leading to a criterion for flatness.
In \cite{SGA2new}*{Exposé~X, \S3}, Grothendieck reformulated their results into a purity concerning finite \'etale covers and proved this purity on {a Noetherian local ring} that is a complete intersection of dimension $\geq 3$ by reducing the assertion to hypersurfaces via several passages involving formal completions.
Nevertheless, a practical deficiency of the {latter} argument is that, even over a rank-one valuation ring $V$ with pseudo-uniformizer $\varpi$, the coherence of the $\varpi$-adic completion $\wh{A}$ of a certain local $V$-algebra $A$ is unknown to us, not to mention the crucial primary decomposition on it (that will guarantee a certain finiteness result).
To circumvent this technical obstacle, we revert to Auslander's argument by establishing a Pr\"uferian counterpart \Cref{Auslander's flatness criterion} of the criterion for flatness \cite{Aus62}*{Theorem~1.3}.
Granted this, we acquire the Pr\"uferian Zariski--Nagata \Cref{fet cover}: the pullback
\[
 \text{$\mathrm{F\'Et}_X\isoto \mathrm{F\'Et}_{X\backslash Z}$ \q is an equivalence }
\]
for every closed subset $Z\subset X$ in the Basic setup I \S\ref{intro-setup} that satisfies the conditions
\[
\begin{cases}
    \codim(Z_{\eta}, X_{\eta})\geq 2  & \text{for each generic point } \eta\in S, \text{ and}  \\
    \codim(Z_s,X_s)\geq 1             & \text{for all } s\in S.
\end{cases}
\]
In particular, if $X$ is connected and $\ov{x}\colon \Spec \Omega\ra X\backslash Z $ is a geometric point, then
\[
  \text{the map \q $\pi_1^{\et}(X\backslash Z, \ov{x})\ra \pi_1^{\et}(X,\ov{x})$\q is bijective.}
\]
\epp




\bpp[Basic setup II]\label{intro basic setup ii}
The rest of this section deals mainly with the following.
For a semilocal Pr\"ufer domain $R$ with fraction field $K$, an integral $R$-smooth scheme $X$, the semilocalization $A\ce \sO_{X,\b{x}}$ of $X$ at a finite subset $\b{x}\subset X$ contained in a single affine open of $X$, and a \emph{quasi-split} reductive $A$-group scheme $G$, we study the trivialization behaviour of $G$-torsors.
\epp

\bpp[Grothendieck--Serre for quasi-split groups]
The Grothendieck--Serre conjecture predicts that every torsor under a reductive group scheme $G$ over a regular local ring $A$ is trivial if it becomes trivial over $\Frac A$.
This conjecture was settled in the affirmative when $A$ contains a field (that is, $A$ is of equicharacteristic), but in the mixed characteristic case, except for several sporadic or low dimensional cases, the conjecture remains open beyond quasi-split groups \cite{Ces22a}. For a detailed review of the state of the art in this area, see \cite{Pan18}*{\S 5} , as well as \cite{GL23}*{\S 1.2} for a summary of recent developments.
In this paper, we prove \Cref{qs-torsors}~\ref{qs-GS-rk1}, thereby generalizing the main result of \cite{Ces22a} to the Pr\"uferian context: in the basic setup \S\ref{intro basic setup ii}, assume in addition that the Pr\"ufer ring $R$ is of Krull dimension 1, then every generically trivial $G$-torsor is trivial, that is, we have
    \[
     \Ker \left(H^1(A,G)\to H^1(\Frac A,G)\right) =\{*\}.
    \]
The proof follows a similar strategy of \cite{Ces22a} (with its earlier version given by Fedorov \cite{Fed22b}), and the key input is our toral version of purity \Cref{purity for gp of mult type} and Grothendieck--Serre type \Cref{G-S type results for mult type} in this context. More precisely, by the valuative criterion of properness, a generically trivial torsor on $X$, say, reduces to a generically trivial torsor under a Borel $B$ away from a closed subset $Z$ of $X$ that has codimension $\ge 2$ (resp., $\ge 1$) in the generic (resp., non-generic) $R$-fibre. Further, thanks to the aforementioned toral purity and Grothendieck--Serre type results, the above $B$-torsor even reduces to a $\text{rad}^u(B)$-torsor on $X \backslash Z$. Then, with the help of the geometric \Cref{geome. input} (unfortunately, whose validity imposes the dimension-1 constraint on $R$), we can reduce to studying torsors over the relative affine line via excision and patchings, and we then conclude by \cite{GL23}*{Theorem 5.1}.

\epp

\bpp[A version of Nisnevich's purity conjecture for quasi-split groups]
Now, we turn to Nisnevich's purity conjecture, where we require the total isotropicity of group schemes.
A reductive group scheme $G$ over a scheme $S$ is \emph{totally isotropic} at a point $s\in S$ if in the following canonical decomposition
\[
   \tst G^{\mathrm{ad}}_{\sO_{S,s}}\cong \prod_i\mathrm{Res}_{A_i/\sO_{S,s}}(G_i)
\]
(\emph{c.f.}~\cite{SGA3IIInew}*{Exposé~XXIV, Proposition~5.10~(i)}) every $G_i$ 
contains a $\bG_{m,A_i}$, where $\sO_{S,s}\to A_i$ is finite \'etale, and $G_i$ is an adjoint semisimple $A_i$-group scheme whose geometric $A_i$-fibres have connected Dynkin diagram of fixed type $i$.
If this holds for all $s\in S$, then $G$ is \emph{totally isotropic}. For instance, tori and quasi-split group schemes are totally isotropic.

Proposed by Nisnevich \cite{Nis89}*{Conjecture~1.3} and modified due to the anisotropic counterexamples of Fedorov \cite{Fed22b}*{Proposition~4.1}, the Nisnevich conjecture predicts that, for a regular semilocal ring $R$, a regular parameter $r$ (i.e., $r\in \fm\backslash \fm^2$ for every maximal ideal $\fm\subset R$), and a reductive $R$-group scheme $G$ such that $G_{R/rR}$ is totally isotropic, every generically trivial $G$-torsor on $R[\f{1}{r}]$ is trivial, namely,
\[
  \tst \text{ $\Ker(H^1(R[\f{1}{r}],G)\ra H^1(\Frac R,G))=\{\ast\}$.}
\]
The case when $R$ is a local ring of a regular affine variety over a field and $G=\GL_n$ was settled by Bhatwadekar--Rao \cite{BR83} and was subsequently extended to arbitrary regular local rings containing fields by Popescu \cite{Pop02}*{Theorem~1}.
 Nisnevich in \cite{Nis89} proved the conjecture in dimension two, assuming that $R$ is a local ring with infinite residue field and that $G$ is quasi-split.
For the state of the art, the conjecture was settled in equicharacteristic case and in several mixed characteristic case by {\v{C}}esnavi{\v{c}}ius in \cite{Ces22b}*{Theorem~1.3} (previously, Fedorov \cite{Fed21} proved the case when $R$ contains an infinite field).
Besides, the toral case and some low dimensional cases are known and surveyed in \cite{Ces22}*{Section~3.4.2~(1)} including Gabber's result \cite{Gab81}*{Chapter~I, Theorem~1} for the local case $\dim R\leq 3$ when $G$ is either $\GL_n$ or $\mathrm{PGL}_n$.
We prove \Cref{qs-torsors}~\ref{qs-Nis-conj}: in the Basic Setup II \S\ref{intro basic setup ii},
\[
    \Ker \left(H^1(A\otimes_RK,G)\to H^1(\Frac A,G)\right) =\{*\}.
\]
\epp

\bpp[Notations and conventions]\label{notation-convention}
All rings in this paper are commutative with units, unless stated otherwise.
{Also, we adopt the notion of \emph{normal schemes} as in \cite{EGAI}*{4.1.4}, that is, they are schemes whose all local rings are integrally closed domains.
A \emph{Pr\"ufer scheme} is a scheme that is covered by spectra of Pr\"ufer domains.}
For a point $s$ of a scheme (resp., for a prime ideal $\mathfrak{p}$ of a ring), we let $\kappa(s)$ (resp., $\kappa(\mathfrak{p})$) denote its residue
field. For a global section $s$ of a scheme $S$, we write $S[\frac{1}{s}]$ for the open locus where $s$ does not
vanish. For a ring $A$, we let $\text{Frac}\, A$ denote its total ring of fractions.
For a morphism of algebraic spaces $S'\to S$, we let $(-)_{S'}$ denote the base change functor from $S$ to $S'$; if $S=\text{Spec}\,R$ and $S'=\text{Spec}\,R'$ are both affine schemes, we will often write $(-)_{R'}$ for $(-)_{S'}$.

 Let $S$ be an algebraic space, and let $G$ be an $S$-group algebraic space. For an $S$-algebraic space $T$, by a $G$-torsor over $T$ we shall mean a $G_T\ce G\times_RT$-torsor (see \Cref{def of torsos}).  Denote by $\textbf{Tors}(S_{\fppf},G)$ (resp., $\textbf{Tors}(S_{\et},G)$) the groupoid of $G$-torsors on $S$ that are fppf locally (resp., \'etale locally) trivial; specifically, if $G$ is $S$-smooth (\emph{e.g.}, $G$ is $S$-reductive, see below), then every fppf locally trivial $G$-torsor is \'etale locally trivial, so we have
$$\textbf{Tors}(S_{\fppf},G)= \textbf{Tors}(S_{\et},G).$$

For a scheme $X$, let $\mathbf{Pic}(X)$ denote the category of invertible $\sO_X$-modules. When $X$ is locally coherent (\S \ref{coh rings & schemes}), let $\sO_X\x{-}\mathbf{Rflx}$ denote the category of reflexive $\sO_X$-modules.

Let $S$ be an algebraic space. By a reductive $S$-group algebraic space we mean a smooth affine $S$-group algebraic space whose geometric $S$-fibres are (connected) reductive algebraic groups. When $S$ is a scheme, this aligns with the definition of reductive $S$-group schemes presented in \cite{SGA3IIInew}.
\epp

\subsection*{Acknowledgements}
The authors wish to express their gratitude to K\k{e}stutis \v{C}esnavi\v{c}ius and Ivan Panin for their unwavering support throughout the development of this work.
Additionally, the authors extend their appreciation to Matthew Morrow and Colliot-Thélène for proposing the Grothendieck--Serre conjecture on smooth schemes over semilocal Pr\"ufer rings during the defense of the first author.
The authors also thank K\k{e}stutis \v{C}esnavi\v{c}ius, Arnab Kundu, Shang Li, Ivan Panin, and Jiandi Zou for insightful discussions on various aspects of the article over the past few months.
The authors are especially grateful to Laurent Moret-Bailly for providing useful remarks and correspondence, which have helped to improve the quality of the work.
The authors would like to thank Enlin Yang from Peking University and Weizhe Zheng from the Morningside Center of Mathematics (Chinese Academy of Sciences) for their hospitality.
Finally, the authors would like to express their gratitude to the referees for their careful reading, insightful remarks and suggestions, which significantly enhanced the quality of this paper.
This project has received funding from the European Research Council (ERC) under the European Union's Horizon 2020 research and innovation programme (grant agreement No.~851146).
This work was done under support of the grant №075-15-2022-289 for
creation and development of Euler International Mathematical Institute.

\section{Reflexive sheaves and depths}
\bpp[Locally coherent schemes]\label{coh rings & schemes}
A module $M$ finitely generated over a ring $A$ is \emph{coherent} if its   finitely generated $A$-submodules are all finitely presented.
A ring $A$ is \emph{coherent} if it is coherent as an $A$-module.
It is worth noting that one can determine the coherence of either a module or a ring Zariski locally. Coherent rings characteristically encompass Noetherian rings, but of primary significance to our discussion are finitely generated flat algebras over Prüfer domains, as referenced in \Cref{coh of V-flat ft alg}.

{On a scheme $X$, a quasi-coherent $\sO_X$-module $\sF$  is \emph{coherent} if there exists an affine open cover $X=\cup_i U_i$ such that $\sF(U_i)$ is a coherent $\sO_X(U_i)$-module for every $i$.} In such instances, this property holds true for all affine open covers of $X$.  A scheme $X$ is \emph{locally coherent} if $\sO_X$ is coherent as an $\sO_X$-module. A locally coherent scheme is \emph{coherent} if it is quasi-compact quasi-separated.

Given a scheme $X$, the \emph{dual} of an $\sO_X$-module $\sF$ is defined as
\[
  \sF^{\vee}\ce \sHom_{\sO_X}(\sF,\sO_X).
\]
\epp
\blem\label{lem-coh}
Let $X$ be a scheme and let $\sF$ and $\sG$ be coherent $\sO_X$-modules.
\benumr
\item\label{ker-coker-coh} {If $\sF\overset{f}{\ra} \sG$ is a morphism of $\sO_X$-modules, then $\Ker f$ and $\Coker f$ are coherent.}
\item\label{hom-tf} Assume that $X$ is integral. If $\sG$ is $\sO_X$-torsion-free, so is $\sHom_{\sO_X}(\sF,\sG)$. 
\[
   \x{In particular,\q $\sF^{\vee}$ is $\sO_X$-torsion-free.}
\]

Now, assume that $X$ is locally coherent.

\item\label{coh modules on coh rings} An $\sO_X$-module is coherent if and only if it is Zariski locally finitely presented.
\item\label{hom-coh} {$\sHom_{\sO_X}(\sF,\sG)$ is coherent. In particular, $\sF^{\vee}$ is coherent.}
\eenum
\elem
\bpf
We will argue Zariski locally on $X$.
For \ref{ker-coker-coh}, see \SP{01BY}.
For \ref{hom-tf}, we apply $\sHom_{\sO_X}(-,\sG)$ to $\sO_X^I\surjects \sF$ and get an embedding $\sHom_{\sO_X}(\sF,\sG)\hra \prod_{I}\sG$. As $\prod_{I}\sG$ is $\sO_X$-torsion-free, so is $\sHom_{\sO_X}(\sF,\sG)$.
For \ref{coh modules on coh rings}, see \SP{05CX}.
For \ref{hom-coh}, we apply \ref{coh modules on coh rings} to choose a finite presentation $\sO_X^{\oplus m}\ra \sO_X^{\oplus n}\ra \sF\ra 0$ and take $\sHom_{\sO_X}(-,\sG)$ of it.
Then $\sHom_{\sO_X}(\sF,\sG)$ is a kernel of a map of coherent modules, so by \ref{ker-coker-coh} is coherent.
\epf
\bppt[Reflexive sheaves]
For a locally coherent scheme $X$ and a coherent $\sO_X$-module $\sF$, consider the canonical map
\[
   \text{$\phi_{\sF}\colon \sF\ra \sF^{\ddual}$}.
\]
If $\phi_{\sF}$ is an isomorphism, then $\sF$ is \emph{reflexive}.
We let $\sO_X\x{-}\mathbf{Rflx}$ denote the category of reflexive $\sO_X$-modules.
For instance, every vector bundle on $X$ is reflexive.
A $\sO_X$-module $\sG$ is \emph{Zariski locally finitely copresented} if it Zariski locally fits into an exact sequence $0\ra \sG\ra \sO_X^{\oplus m}\ra \sO_{X}^{\oplus n}$ for some integers $m$ and $n$.
 If $X$ is further assumed to be integral, then the following \Cref{lem-rf}~\ref{hom-rf} shows that, for every coherent $\sO_X$-module $\sG$, the double dual $\sG^{\vee\!\vee}$ is $\sO_X$-reflexive, hence $\sG^{\vee\!\vee}$ is called the \emph{reflexive hull} of $\sG$.
\eppt

\blem\label{lem-rf}
Let $X$ be an integral locally coherent scheme and let $\sF$ and $\sG$ be coherent $\sO_X$-modules.
\benumr
\item\label{ddual-tf} The double dual $\phi_{\sF}\colon \sF\hra \sF^{\ddual}$ is injective if and only if $\sF$ is $\sO_X$-torsion-free. In particular, all reflexive $\sO_X$-modules are $\sO_X$-torsion-free.
\item\label{hom-rf} If $\sG$ is reflexive, then so is $\sHom_{\sO_X}(\sF,\sG)$. 
Therefore,
\[
   \text{the dual of any coherent $\sO_X$-module is reflexive.}
\]
\item\label{rf-fcpr} $\sF$ is reflexive if and only if it is Zariski locally finitely copresented, if and only if there are a finite locally free $\sO_X$-module $\sL$ and a torsion-free $\sO_X$-module $\sN$ fitting into the short exact sequence
\[
   0\ra \sF\ra \sL\ra \sN\ra 0.
\]
\item\label{dual-rf-iso} If $\sG$ is reflexive, then the following natural map $\sHom_{\sO_X}(\phi_{\sF},\sG)$
\[
   \text{$\sHom_{\sO_X}(\sF^{\ddual},\sG)\isoto \sHom_{\sO_X}(\sF,\sG)$\q is an isomorphism.}
\]
\eenum
\elem
\bpf
We may prove all the assertions Zariski locally.
\benumr
\item For \ref{ddual-tf}, by \Cref{lem-coh}~\ref{hom-tf}, the injectivity of $\phi_{\sF}$ implies that $\sF$ is $\sO_X$-torsion-free.
Conversely, if $\sF$ is torsion-free, then locally $\sF\subset \sO_X^{\oplus n}$ for some integer $n$; taking double dual, we find that the composite $\sF\ra \sF^{\ddual} \ra \sO_X^{\oplus n}$ is injective, so $\phi_{\sF}$ is also injective.
\item For \ref{hom-rf}, we apply \Cref{lem-coh}~\ref{coh modules on coh rings} to choose a finite presentation and take its $\sHom_{\sO_X}(-,\sG)$.
So we get an exact sequence $0\ra \sHom_{\sO_X}(\sF,\sG)\overset{u}{\ra} \sG^{\oplus m}\ra \sG^{\oplus n}$ for some integers $m$ and $n$.
Taking double dual of this exact sequence, we obtain the following commutative diagram
\[
\begin{tikzcd}
0 \arrow[r] & \sHom_{\sO_X}(\sF,\sG) \arrow[r, "u"] \arrow[d] & \sG^{\oplus m} \arrow[r] \arrow[d] & \sG^{\oplus n} \arrow[d] \\
            & \sHom_{\sO_{X}}(\sF,\sG)^{\ddual} \arrow[r, "u^{\ddual}"]    & (\sG^{\ddual})^{\oplus m} \arrow[r]           &  (\sG^{\ddual})^{\oplus n}
\end{tikzcd}
\]
the reflexivity of $\sG$ and diagram chase reduces us to showing that $u^{\ddual}$ is injective.
Take $(-)\otimes_{\sO_X}K(X)$ of the exact sequence
\[
(\sG^{\vee})^{\oplus m}\overset{u^{\vee}}{\ra}\sHom_{\sO_X}(\sF,\sG)^{\vee}\ra \Coker(u^{\vee})\ra 0.
\]
Since $\sF$ is finitely presented, by \SP{0583}, we have $\Coker(u^{\vee})_{K(X)}=0$, namely, $\Coker(u^{\vee})$ is torsion. 
Thus $\Ker(u^{\ddual})=\Coker(u^{\vee})^{\vee}=0$, as desired.
\item If $\sF$ is reflexive, we choose a finite presentation of $\sF^{\vee}$ (which is coherent by \Cref{lem-coh}~\ref{hom-coh}) and then taking its dual yields a finite copresentation of $\sF^{\ddual} \simeq \sF$.
Conversely, if $\sF$ is finitely copresented as $0\ra \sF \ra \sO_X^{\oplus m}\overset{v}{\ra} \sO_X^{\oplus n}$, then $\sF\simeq \Coker(v^{\vee})^{\vee}$, which is reflexive by \ref{hom-rf}.
Taking $\sL\ce \sO_{X}^{\oplus m}$ and $\sN\ce \mathrm{Im(v)}$, we obtain the desired short exact sequence.
\item Pick a finite copresentation $\sO_X^{\oplus n}\ra \sO_X^{\oplus m}\ra \sG\ra 0$ and take $\sHom(\phi_{\sF},-)$ on it.
We have the commutative diagram of $\sO_X$-modules with exact rows
\[
      \tikz {
      \node (E) at (-7.5,0) {$0$};
      \node (F) at (-7.5,-1.5) {$0$};
      \node  (C) at (-5.4,0) {\small$\sHom_{\sO_X}(\sF^{\ddual},\sG)$};
      \node  (A) at (-1.7,0) {$\sHom_{\sO_X}(\sF^{\ddual},\sO_X^{\oplus m})$};
      \node  (B) at (2.2,0) {$\sHom_{\sO_X}(\sF^{\ddual},\sO_X^{\oplus n}) $};
      \node  (c) at (-5.4,-1.5) {$\sHom_{\sO_X}(\sF,\sG)$};
      \node  (a) at (-1.7,-1.5) {$\sHom_{\sO_X}(\sF,\sO_X^{\oplus m})$};
      \node  (b) at (2.2,-1.5)   {$\sHom_{\sO_X}(\sF,\sO_X^{\oplus n}).$};
      \node  (L) at (-4.2,-0.73) {\small$\sHom_{\sO_X}(\phi_{\sF},\sG)$};
       \draw [draw = black, thin,
      arrows={
      - Stealth }]
      (A) edge  (B);
      \draw [draw = black, thin,
      arrows={
      - Stealth }]
      (a) edge (b);
       \draw [draw = black, thin,
      arrows={
      - Stealth }]
      (F) edge (c);
       \draw [draw = black, thin,
      arrows={
      - Stealth }]
      (E) edge (C);
      \draw [draw = black, thin,
      arrows={
      - Stealth }]
      (B) edge (b);
      \draw [draw = black, thin,
      arrows={
      - Stealth }]
      (A) edge (a);
      \draw [draw = black, thin,
      arrows={
      - Stealth }]
      (C) edge  (c);
      \draw [draw = black, thin,
      arrows={
      - Stealth }]
      (C) edge  (A);
      \draw [draw = black, thin,
      arrows={
      - Stealth }]
      (c) edge  (a);
      }
      \]
By the assertion \ref{hom-rf}, $\sF^{\vee}$ is reflexive, hence $\phi_{\sF}^{\vee,m}, \phi_{\sF}^{\vee,n}$ are bijective and so is $\sHom_{\sO_X}(\phi_{\sF},\sG)$. \qedhere
\eenum
\epf

By reflexive hull, reflexive sheaves extend from quasi-compact open (\emph{cf.}~\cite{GR18}*{Proposition~11.3.8~(i)}).
\bcor\label{ext-rflx}
For a quasi-compact open $U$ in a coherent integral scheme $X$,
\[
   \text{the restriction \q $\sO_X\x{-}\mathbf{Rflx}\ra \sO_U\x{-}\mathbf{Rflx}$ \q is essentially surjective.}
\]
\ecor
\bpf
A reflexive $\sO_U$-module $\sF$, by \SP{0G41} and \Cref{lem-coh}~\ref{coh modules on coh rings}, extends to a coherent $\sO_X$-module $\wt{\sF}$. 
By \Cref{lem-rf}~\ref{hom-rf}, the reflexive  $\wt{\sF}^{\ddual}$ extends $\sF$.
\epf
Now, we recall the notion of depth in terms of local cohomology and use it to describe reflexive sheaves.
\bpp[Depths] \label{setup of depth & proj dim}

For a scheme $X$ and an open immersion $j\colon U\hookrightarrow X$ with closed complement $i\colon Z\ce X\backslash U\hookrightarrow X$, consider the functor
\[
\un{\GG}_Z\colon \mathbf{Z}_X\x{-}\mathbf{Mod} \to \mathbf{Z}_X\x{-}\mathbf{Mod} \qq \sF \mapsto \Ker(\sF\to j_*j^*\sF)
\]
which sends an abelian sheaf $\sF$ on $X$ to its largest subsheaf supported on $Z$.
The functor $\un{\GG}_Z$ is left exact, giving rise to a right derived functor
\[
R\underline{\GG}_Z:D^+(\mathbf{Z}_X\x{-}\mathbf{Mod})\to D^+(\mathbf{Z}_X\x{-}\mathbf{Mod}).
\]
Explicitly, for an object $K\in D^+(\mathbf{Z}_X\x{-}\mathbf{Mod})$, taking an injective resolution $K\ra I^{\bullet}$, one computes $R\un{\GG}_Z(K)$ via $R\un{\GG}_Z(I^{\bullet})$.
Thus, $R\un{\GG}_Z$ actually factors over the canonical functor $D^+(\mathbf{Z}_Z\x{-}\mathbf{Mod})\to D^+(\mathbf{Z}_X\x{-}\mathbf{Mod})$ and satisfies $\un{\GG}_Z\simeq R^0\un{\GG}_Z$. Moreover, we have the following exact sequence
$
0\to \un{\GG}_Z(I^{\bullet})\to I^{\bullet} \to j_*j^*(I^{\bullet})\to 0,
$
giving rise to the functorial distinguished triangle
\[
R\un{\GG}_Z(K)\to K \to Rj_*j^*(K)\to R\un{\GG}_Z(K)[1].
\]
{Assume now that $j\colon U\hra X$ is \emph{quasi-compact}.} Then, for any quasi-coherent sheaf $\sF$ on $X$, from the above triangle we see that $R\underline{\GG}_Z(\sF)$ has quasi-coherent cohomology sheaves as $Rj_*j^*(\sF)$ does so.
Consider the functor $\GG_Z\ce \GG(X,-)\circ \underline{\GG}_Z$.
Since the functor $\underline{\GG}_Z$ sends injective sheaves to injective sheaves ($\underline{\GG}_Z$ admits an exact left adjoint $i_*i^*$), the derived functor $R\GG_Z$ of $\GG_Z$ is canonically isomorphic to $ R\GG(X,-)\circ R\underline{\GG}_Z$.
Therefore, if $X=\Spec A$ is affine and $M$ is an $A$-module, by the quasi-coherence of $R^i\underline{\GG}_Z(\wt{M})$ we get $R^i\underline{\GG}_Z(\wt{M}) \simeq (R^i\GG_Z(\wt{M}))^{\sim}$ for all $i\ge 0$; in particular, all the $A$-modules $R^i\GG_Z(\wt{M})$ are supported on $Z$. Moreover, we have the following distinguished triangle
\[
R\GG_Z(\wt{M})\to M \to R\GG(U,\wt{M})\to R\GG_Z(\wt{M})[1],
\]
which gives an exact sequence $0\to \GG_Z(\wt{M})\to M \to \GG(U,\wt{M})\to R^1\GG_Z(\wt{M})\to 0$ and isomorphisms $R^i\GG_Z(\wt{M}) \simeq H^{i-1}(U,\wt{M})$ for all $i>1$.

In what follows, we will mainly consider the case where $X=\Spec A$ for a local ring $(A,x\ce \fm_A)$ such that the punctured spectrum $U_A\ce \Spec A\backslash \{x\}$ is \emph{quasi-compact}.
Denote by $j\colon U_A\hra X$ the open immersion. In this case, we have seen that the cohomology modules $R^i\GG_x(\wt{M})$ are all supported on $\{x\}$, where $M$ is an $A$-module.
For simplicity, we shall often write $R^i\GG_x(M)$ for $R^i\GG_x(\wt{M})$. We warn readers that the derived functor $M\mapsto R^i\GG_x(M)$ here is obtained from an injective resolution of the associated quasi-coherent sheaf $\wt{M}$, not from an injective resolution of the $A$-module $M$. These two approaches do not, in general, produce equivalent theories, unless $A$ is Noetherian (see \SP{0A6P}).

As defined in \cite{GR18}*{Definition~10.4.14}, the \emph{depth} of an $A$-module $M$ is given by
\[
\depth_A(M)\ce \sup\{n\in \b{Z}\,|\, R^i\GG_{x} M=0\q \x{for all $i<n$}\}\in \b{Z}_{\ge 0}\cup \{+\infty\}.
\]
For a \emph{finitely generated} $A$-module $N $ supported on $\{x\}$, consider a closely related quantity
\[
\tau_N(M)\ce  \sup\{n\in \b{Z}\,|\, \Ext_A^i(N,M)=0\, \x{for all $i<n$}\}\in \b{Z}_{\ge 0}\cup \{+\infty\}.
\]
We note that by the finite generation of $N$, all the $A$-modules $\Ext_A^i(N,M)$ are supported on $\{x\}$.
\epp

\blem \label{depth and tau}
Let $(A, \fm_A)$ be a local ring with quasi-compact punctured spectrum.
Let $M$ be an $A$-module and $N$ a finitely generated $A$-module supported on $\{\fm_A\}$.
Let $(f_1,\cdots, f_r)$ be an $M$-regular sequence in $\fm_A$.
We have (notations as in \S\ref{setup of depth & proj dim})
\[
\tst \depth_A(M)=\depth_A(M/\sum_{i=1}^rf_iM)+r \, \text{ and }\,  \tau_N(M)=\tau_N(M/\sum_{i=1}^r f_iM)+r.
\]
\elem
\bpf
Denote $x\ce \fm_A$ as the closed point of $\Spec A$. The two equalities are proved similarly, so we only treat the one concerning depths.
By induction on $r$, it suffices to consider for a nonzero $f\in \fm_A$ the following short exact sequence $$0\to M \xrightarrow{f} M \to M/fM \to 0$$ to show that $\depth_A(M)=\depth_A(M/fM)+1$.
We have a long exact sequence
\[
\cdots \to  R^{i-1}\GG_{x}M \xrightarrow{f} R^{i-1}\GG_{x}M \to R^{i-1}\GG_{x} (M/fM) \to R^i\GG_{x}M \xrightarrow{f} R^i\GG_{x} M \to \cdots.
\]
Assume that for an integer $n$ we have $\depth_A(M)\geq n+1$, so $R^i\GG_xM=0$ for $0\leq i\leq n$.
Then the displayed exact sequence implies that $R^i\GG_x(M/fM)=0$ for $0\leq i\leq n-1$, so $\depth_A(M/fM)\geq n$.
Conversely, if $\depth_A(M/fM)\geq n$, then $R^i\GG_x(M/fM)=0$ for $0\leq i\leq n-1$.
The displayed exact sequence implies that $R^i\GG_xM=0$ for $0\leq i\leq n-1$, and there is an injection $R^n\GG_{x}M \overset{\times f}{\hra} R^n\GG_{x} M$.
However, since the $A$-module $R^n\GG_{x}M$ is supported on $\{x\}$ and $f\in \mathfrak{m}_A$, we deduce that $R^n\GG_{x}M=0$, that is, $\depth_A(M)\geq n+1$.
\epf
\begs
\egi \label{exa-val-depth>=2}
Consider a valuation ring $V$ that is not a field. If the punctured spectrum $U_V$ is quasi-compact (\emph{e.g.}, when $0<\dim V<\infty$), then there exists an $f\in \fm_V$ such that $\dim (V/fV)=0$, which implies depth$_V(V/fV)=0$.
From the formula in \Cref{depth and tau}, we deduce that depth$_V(V)=1$. Conversely, if $U_V$ is not quasi-compact, then one can show that depth$_V(V)\ge 2$.

\egi Assume that $A$ is a Noetherian local ring and $N=A/I$ for an ideal {$I\subsetneq A$} (e.g., $N=A/\mathfrak{m}_A$).
Then for any finitely generated $A$-module $M$, we have
\[
\depth_A(M)=\tau_N(M).
\]
Indeed, utilizing \Cref{depth and tau}, one can check that both of them equal the length of any maximal $M$-regular sequence in $\mathfrak{m}_A$ (so the length is independent of all choices).
However, this may be false when $A$ is non-Noetherian.
 For instance, we let $A\ce V$ be a \emph{non-discrete} valuation ring of finite rank and let $N\ce V/\mathfrak{m}_V$ be its residue field.
Take $M=V/fV$ for a nonzero $f\in \mathfrak{m}_V$. Then,
\begin{itemize}
  \item $\depth_V(V/fV)=0$.
  Let $\fp\subsetneq \fm_V$ be the second largest prime ideal of $V$. The non-discreteness of $V$ implies that $\fm_V\neq fV$, so we can pick $g\in \fm_V\backslash (\fp\cup fV)$. Then $f_1\ce f/g$ is in $\fm_V$, and its image $\overline{f}_1\in V/fV$ has annihilator $gV$ which {strictly} contains $\fp$, that is, $0\neq \overline{f}_1\in \GG_{\fm_V}(V/fV)$.
  \item For any nonzero element $\overline{h}\in V/fV$, we have $\fm_V\overline{h}\neq 0$, that is, $\Hom_V(V/\mathfrak{m}_V,V/fV)=0$. Indeed, any such non-invertible $\overline{h}$ is represented by some $h\in \fm_V$ with $fV\subsetneq hV$. Then $\overline{h}$ has annihilator $(f/h)V$, which is strictly contained in $\fm_V$ because $V$ is non-discrete valued.
\end{itemize}
Therefore, we have $\tau_{V/\mathfrak{m}_V}(V/fV)\ge 1 >0=\depth_V(V/fV)$.
\eegs

\blem\label{dep-inh}
Let $A$ be a local ring with quasi-compact punctured spectrum, and let $M,N$ be $A$-modules.
\benumr
\item If $\depth_A(N)\geq 1$, then $\depth_A(\Hom_A(M,N))\geq 1$.
\item If $\depth_A(N)\geq 2$ and $M$ is finite over $A$, then $\depth_A(\Hom_A(M,N))\geq 2$.
\item \label{depth>=2 of refx} {If $A$ is coherent and $M$ is a reflexive $A$-module, then }\[
   \depth_A(M)\ge \min(2, \depth_A(A)).
\]
\eenum
\elem
\bpf
Pick a presentation $A^{\oplus J}\ra A^{\oplus I}\ra M\ra 0$ and take $\Hom_A(-,N)$, we get a short exact sequence
\[
 \tst  0\ra \Hom_A(M,N)\ra \prod_I N\ra N\pr\ra 0,
\]
where $N\pr\subset \prod_J N$ for some index set $J$.
 Since $\depth_A(N)\geq 1$, by definition, we have $\depth_A(\prod_J N)\geq 1$ and $\depth_A(N\pr)\geq 1$.
For (ii), take $R\GG_{x}$ of this short exact sequence, where $x$ is the closed point of $\Spec A$.
Since $\depth_A(N)\geq 2$ and $I$ is finite, we have $R\GG^i_{x}(\prod_I N)=\bigoplus_I R\GG^i_{x}N=0$.
Combine this with $R\GG^0_{x}N\pr=0$, both $R^0\GG_{x}$ and $R^1\GG_{x}$ of $\Hom_A(M,N)$ vanish, so we get $\depth_A(\Hom_A(M,N))\geq 2$. The last assertion follows from (i) and (ii).
\epf
\bcor\label{ref-dep-2}
Let $X$ be an integral, coherent, {and topologically locally Noetherian} scheme and $j\colon U\hookrightarrow X$ an open immersion such that every $z\in X\backslash U$ satisfies $\depth_{\sO_{X,z}}(\sO_{X,z})\geq 2$.
Then, every reflexive $\sO_X$-module $\sF$ satisfies
\[
\sF\isoto j_{\ast}j^{\ast}\sF.
\]
\ecor
\bpf
As $U$ contains the generic point of $X$, the injectivity follows because $\sF$ and $j_{\ast}j^{\ast}\sF$ are subsheaves of their common generic stalk.
To show the surjectivity, we show that every section $s\in \sF(U)$ extends over $X$.
Let $U_1\subset X$ be the domain of definition of $s$; it is open and contains $U$. If $U_1\neq X$, then every maximal point $z\in X\backslash U_1$ is contained in $X\backslash U$. In particular, we have $\depth_{\sO_{X,z}}(\sO_{X,z})\geq 2$ by assumption.
By \Cref{dep-inh}~\ref{depth>=2 of refx}, the localization $\sF_z\ce \sF|_{\sO_{X,z}}$ satisfies $\depth_{\sO_{X,z}}(\sF_z)\geq 2$, hence the restriction of $s$ to $\Spec \sO_{X,z}\cap U_1=\Spec \sO_{X,z}\backslash \{z\}$ extends to a section $s_1\in \GG(V_1, \sF)$ for an affine open neighborhood $V_1$ of $z$ in $X$ (by the quasi-compactness of the open immersion $U_1\hookrightarrow X$).
 Since $X$ is integral and $\sF$ is $\sO_X$-torsion-free, the two local sections $s_1$ and $s$ agree on $V_1\cap U$, and thus can be glued to a section over $V_1\cup U$.
 This implies that $z$ is already in $U_1$, a contradiction.
\epf

\bcor\label{rf-equiv}
Let $X$ be an integral, coherent, {and topologically locally Noetherian} scheme. Let $j:U\hookrightarrow X$ be an open subscheme such that every $z\in X\backslash U$ satisfies $\depth_{\sO_{X,z}}(\sO_{X,z})\geq 2$. Then, taking $j^{\ast}$ and $j_{\ast}$ induce an equivalence of categories
\[
   \sO_X\x{-}\mathbf{Rflx}\isoto \sO_U\x{-}\mathbf{Rflx}.
\]
\ecor
\bpf
By \Cref{ext-rflx}, the restriction functor is essentially surjective.
To show the full faithfulness, we pick two reflexive $\sO_X$-modules $\sF,\sG$ and consider $\sH\ce \sHom_{\sO_X}(\sF,\sG)$.
By \Cref{lem-rf}~\ref{hom-rf}, $\sH$ is a reflexive $\sO_X$-module, so \Cref{ref-dep-2} implies that $\sH(X)\simeq \sH(U)$, as desired.
\epf

\bpp[(Weakly) {associated} primes]\label{par-(Weakly) associated primes}

Let $R$ be a ring and $M$ an $R$-module. We define the set of \emph{associated primes} of $M$, denoted as $\operatorname{Ass}_R(M)$, to be the collection of prime ideals $\fp \subset R$ for which $R^0\GG_{\fp}M_{\fp}\neq 0$.
{Note that in the literature, the associated primes referred to in this context are often termed as “weakly associated primes” (see \SP{0546} and \cite{Bou98}*{Chapitre~IV, \S1, exercice~17)}).} It then follows that Ass$_R(M)\neq \emptyset$ whenever $M\neq 0$: for any nonzero $m\in M$, Ass$_R(M)$ contains all the minimal elements of
\[
\text{Supp}(m)\ce \{\fp\in \Spec\, R\colon m_{\fp}\neq 0\in M_{\fp}\}.
\]
\epp
Now, we present \Cref{HomAss}, whose argument is pointed out by L. Moret-Bailly, as a generalization of its Noetherian case \cite{SGA2new}*{Exposé~III, Proposition~1.3}.
In addition to its generality, this result leads to an alternative proof of \Cref{a lem on linear algebra} concerning the ``matrix of direct sums of modules''.
\bprop\label{HomAss}
Let $A$ be a ring, $M$ a finitely presented $A$-module, and $N$ an $A$-module.
We have
\[
   \Ass_A(\Hom_A(M,N))=\Supp(M)\cap \Ass_A(N).
\]
\eprop
\bpf
It is clear that $ \Ass_A(\Hom_A(M,N))\subset \Supp(M)\cap \Ass_A(N)$, so it suffices to prove the converse inclusion.
For a prime ideal $\fp\in \Supp(M)\cap \Ass_A(N)$, we have $\fp\in \Supp(M_{\fp})\cap \Ass_{A_{\fp}}(N_{\fp})$.
As $M$ is finitely presented, it suffices to show that there is an $f\in \Hom_A(M,N)$ such that $\fp$ the minimal one in
 \[
   \Supp(f)=\{\fq\in \Spec A\,|\, 0\neq f_{\fq}\in \Hom_{A_{\fq}}(M_{\fq},N_{\fq})\}.
\]
Hence we are reduced to showing the following local case \Cref{HomAss-local}.
\epf
\blem\label{HomAss-local}
Let $(A,\fm)$ be a local ring, $M$ a finitely presented nonzero $A$-module, and $N$ an $A$-module with $\Supp(N)=\{\fm\}$.
Then we have $\Hom_A(M,N)\neq 0$.
\elem
\bpf
As $N$ has an $A$-submodule of the form $A/J$ for an ideal $J\subset A$ with $\Supp(J)=\{\fm\}$, we may replace both $A$ and $N$ by $A/J$ and replace $M$ by $C/JC$ iteratively to assume that $A$ is a local ring of dimension zero.
It suffices to prove that $M^{\vee}\neq 0$.
As $M$ has a nonzero quotient of the form $A/I$ for a finitely generated ideal $I\subset A$, it suffices to note that $(A/I)^{\vee}=\Ann(I)\neq 0$, because $I$ is nilpotent.
\epf

\blem\label{lem-iso}
Let $R$ be {a} domain with topologically Noetherian spectrum and $\alpha\colon M\ra N$ a morphism of $R$-modules.
Then $\alpha$ is an isomorphism, provided that $N$ is torsion-free and every prime $\fp\subset R$ satisfies
\[
   \x{either  $\alpha_{\fp}\colon M_{\fp}\isoto N_{\fp}$  is an isomorphism,\q or $\depth_{R_{\fp}}(M_{\fp})\geq 2$.}
\]
\elem
\bpf
Take $K\ce \Ker(\alpha)$. If $K\neq 0$, we choose an associated prime $\fp\in \text{Ass}_{R}(K)$, then $0\neq R^0\GG_{\fp}K_{\fp}\subset R^0\GG_{\fp}M_{\fp}=0$ (since depth$_{R_{\fp}}(M_{\fp})\ge 2>0$), a contradiction. This proves that $\alpha$ is injective.
Take $Q\ce \text{Coker}(\alpha)$. If $Q\neq 0$, we choose an associated prime $\fp\in \text{Ass}_{R}(Q)$. Since $R^0\GG_{\fp}M_{\fp}=R^1\GG_{\fp}M_{\fp}=R^0\GG_{\fp}N_{\fp}=0$, taking $R\GG_{\fp}$ of the short exact sequence
$0\ra M_{\fp}\ra N_{\fp}\ra Q_{\fp}\to 0$ yields that $R^0\GG_{\fp}Q_{\fp}=0$. a contradiction. Thus $\alpha$ is also surjective.
\epf

\bppt[Weak dimensions]\label{par-Weak dimesions}
Recall the \emph{weak dimension} of a ring $B$, denoted by $\wdim(B)$, is defined as
\[
   \wdim(B)=\sup\{\fldim_B(M)\,|\, \x{$M$ is a $B$-module}\},
\]
the supremum of flat dimensions of all $B$-modules.
   For example, we know from \SP{092S} that wdim$(B)\le 1$ if and only if all local rings of $B$ are valuation rings. {Note that if a $B$-module $M$ has a (possibly infinite) resolution by finite free modules, then $\fldim_{B}(M)=\pd_B(M)$. This holds in particular if $B$ is coherent and $M$ is finitely presented. Since any such $M$ is an iterated extension of cyclic, finitely presented $B$-modules, we conclude that (for $B$ coherent)
   \[
   \wdim(B)=\sup\{\pd(B/J)\,|\, \x{$J\subset B$ is a finitely generated ideal}\}.
   \]}
\eppt

\bpp[Coherent regular rings]\label{reg-coh}
A coherent ring $R$ is termed \emph{regular} if every finitely generated ideal of $R$ has finite projective dimension. By Serre's homological characterization, one recovers the classical regularity for Noetherian rings. However, the primary focus of this paper is on the class of non-Noetherian regular rings (e.g., a flat finite type ring over a valuation ring with regular fibres, see \Cref{coh of V-flat ft alg}).  For readers' convenience, we briefly summarize the key properties of general regular coherent rings, which closely resemble their Noetherian counterparts. Other useful properties can be found in \Cref{apdxA}.
\epp

\bthm \label{ring property coh reg}
For a coherent regular local ring $A$, the following assertions hold.
\benumr
\item\label{ring property coh reg-i} $A$ is a normal domain; more precisely
\item\label{ring property coh reg-ii} $A$ is the intersection (in $\Frac A$) of its local rings which are valuation rings;
\item\label{ring property coh reg-iii} Any two nonzero elements of $A$ have a greatest common divisor and a least common multiple;
\item\label{qc coh reg has fin wdim} If the punctured spectrum $U_A$ is quasi-compact, then
\[
   \depth_A(A)=\wdim(A)<\infty.
\]
So, for every coherent local ring $R$ with quasi-compact punctured spectrum, 
\[
   \wdim(R)< \infty \q \x{if and only if}\q \x{$R$ is regular.}
\]
\eenum
\ethm

\bpf \hfill
\benumr
\item This follows from \cite{Ber72}*{Corollaire 4.3}.

\item The reference for this is \cite{Que71}*{Proposition 2.4}. Assuming the topological Noetherianness of $\Spec A$, this also follows from (i) and the more general \Cref{rf-s2} (or \Cref{normal-s2}).

\item See \cite{Que71}*{Proposition~3.2}. Assuming the topological Noetherianness of $X\ce\Spec A$, this result can also be deduced from \Cref{Rflx-Pic-Vect on reg} as follows: consider two nonzero elements, $f$ and $g$, in $A$.
Now, examine the two (finitely generated) ideals, $fA+gA$ and $fA\cap gA$, of $A$. These two ideals define two coherent sheaves on $X$, which are invertible around each point of $X$ whose local ring is a valuation ring, or, equivalently, has a weak dimension $\le 1$ (\SP{092S}). By \Cref{Rflx-Pic-Vect on reg}, there exist two invertible sheaves on $X$ that coincide with the two ideal sheaves at all local rings of $X$ having weak dimension $\le 1$. Since $X$ is local, these two invertible sheaves are free, and, in view of (ii), any generators provide a greatest common divisor and a least common multiple of $f$ and $g$.

\item The finiteness of $\wdim(A)$ under the quasi-compactness assumption on $U_A$ was originally established by Quentel \cite{Que71}*{Corollaire~1.1}. We refer to \Cref{wdim=depth} for additional details. \qedhere
\eenum
\epf
\blem\label{wdim-2}
Let $R$ be a coherent ring and $n\geq 2$ an integer.
Denote $(-)^{\vee}\ce \Hom_R(-,R)$.
Then
\[
   \x{$\wdim(R)\leq n$\q if and only if \q $\fldim_R(M^{\vee})\leq n-2$ for every $R$-module $M$.}
\]
\elem
\bpf
If $\wdim(R)\leq n$, then for any $R$-module $M$ one takes a presention $F_2\ra F_1\ra M\ra 0$ for free $R$-modules $F_1$ and $F_2$.
As $R$ is coherent, by \SP{05CZ} and Lazard's theorem \SP{058G}, $F_1^{\vee}$ and $F_2^{\vee}$ are flat, so the exact sequence $0\ra M^{\vee}\ra F_1^{\vee}\ra F_2^{\vee}\ra C\ra 0$ for $C\ce \Coker(F_1^{\vee}\ra F_2^{\vee})$ yields $\fldim_R(M^{\vee})\leq n-2$.
Conversely, to estimate $\wdim(R)$, it suffices to consider finitely presented $R$-modules and their flat dimensions.
Every finitely presented $R$-module $M$ fits into an exact sequence
\[
   0\ra K\ra R^{\oplus m}\ra R^{\oplus n}\ra M\ra 0,
\]
where $K\ce \Ker(R^{\oplus m}\ra R^{\oplus n})$, so by \Cref{lem-rf}~\ref{rf-fcpr}, $K$ is a reflexive $R$-module.
Thus, by assumption, we have $\fldim_R(M)=\fldim(K^{\ddual})+2\leq n$.
This gives the desired inequality $\wdim(R)\leq n$.
\epf

A locally coherent scheme is \emph{regular} if it is covered by spectra of coherent regular rings.
\bthm\label{wdim2-vect}
Let $X$ be a locally coherent, integral scheme and let $\sF$ be a coherent $\sO_X$-module.
\benumr
\item\label{wdim2-vect-i} If $\sF$ is reflexive at a point $x\in X$ and $\wdim\sO_{X,x}$ is finite\footnote{For instance, when $\sO_{X,x}$ is regular with a quasi-compact punctured spectrum; see \Cref{ring property coh reg}~\ref{qc coh reg has fin wdim}.}, then we have
\[
\pd_{\sO_{X,x}}\sF_x\leq \max(0,\wdim\sO_{X,x}-2).
\]
\item\label{wdim2-vect-ii} If $\wdim\sO_{X,x}\leq 2$ for all $x\in X$, every reflexive $\sO_X$-module is locally free.
\eenum
\ethm
\bpf
We first note that $\wdim\sO_{X,x}<\infty$ implies that $\sO_{X,x}$ is regular, and the regularity of $\sO_{X,x}$ guarantees that every coherent $\sO_{X,x}$-module is a perfect complex.
The assertion \ref{wdim2-vect-ii} follows from \ref{wdim2-vect-i}.
The assertion \ref{wdim2-vect-i} follows from \Cref{wdim-2} because the perfectness of $\sF_x$ leads to the following desired inequality
\[
   \pd_{\sO_{X,x}}(\sF_x)=\fldim_{\sO_{X,x}}(\sF_x^{\ddual})\leq \max(0,\wdim\sO_{X,x}-2). \qedhere
\]
\epf

\bthm\label{Rflx-Pic-Vect on reg}
Let $X$ be a topologically locally Noetherian, integral, locally coherent, regular scheme.
 For a closed subset $Z\subset X$ with $\wdim \sO_{X,z}\geq 2$ for all $z\in Z$ and the canonical open immersion $j\colon X\backslash Z\hra X$, the restriction $j^*$ and pushforward $j_*$ induce the following equivalences of categories:
\begin{equation*}\label{restriction}
\text{$\sO_{X}\text{-}\mathbf{Rflx}\isoto \sO_{X\backslash Z}\text{-}\mathbf{Rflx}$\qq $\mathbf{Pic}\,X\isoto \mathbf{Pic}\, X\backslash Z$.}
\end{equation*}
In particular, for every scheme $Y$ affine over $X$, we have a bijection of sets
\[
  \tst Y(X)\simeq Y(X\backslash Z).
\]
Moreover, if $\wdim\sO_{X,x}\leq 2$ for all $x\in X$, then we have an equivalence of categories:
\begin{equation*}
\mathbf{Vect}(X) \isoto \mathbf{Vect}(X\backslash Z).
\end{equation*}
\ethm
\bpf
We first note that the regularity of $X$ guarantees that every coherent sheaf on $X$ is a perfect complex.
The equivalence for vector bundles follows from that for reflexive sheaves and \Cref{wdim2-vect} \ref{wdim2-vect-ii}.
By \Cref{ring property coh reg} \ref{qc coh reg has fin wdim}, the depth and weak dimension agree for each local ring of $X$, so the equivalence of reflexive sheaves follows from \Cref{rf-equiv}.
For the assertion concerning $\mathbf{Pic}$, by the result for reflexive modules, it is enough to show that, for every invertible $\sO_{X\backslash Z}$-module $\sI$, its unique $\sO_X$-reflexive extension $\sF\ce j_*\sI$ is actually invertible.
Consider the \emph{category} $\mathbf{gr.Pic}(X)$ \emph{of graded invertible $\sO_X$-modules}, whose objects are pairs $(\sL,\alpha)$ consisting of an invertible $\sO_X$-module $\sL$ and a {locally constant} function $\alpha \colon X\ra \b{Z}$, and morphisms $h\colon (\sL,\alpha)\ra (\sM,\beta)$ satisfies $h_x=0$ if $\alpha(x)\neq \beta(x)$ for all $x\in X$.
By \cite{KM76}*{Theorem~1}, there is a unique determinant functor from $D(\sO_X\x{-}\mathbf{Mod})^{\ast}_{\mathrm{perf}}$ the groupoid of perfect complexes of $\sO_X$-modules
\[
   \mathbf{det}\colon D(\sO_X\x{-}\mathbf{Mod})^{\ast}_{\mathrm{perf}}\ra \mathbf{gr.Pic}(X).
\]
Also, this functor $\mathbf{det}$ commutes with arbitrary base change.
Thus, since the complex $\sF[0]$ is perfect, by \Cref{wdim2-vect} \ref{wdim2-vect-i}, the invertible $\sO_X$-module $\mathbf{det}(\sF[0])$ is well-defined, and we have isomorphisms
\[
\mathbf{det}(\sF[0])\xrightarrow{\sim}j_*j^*\mathbf{det}(\sF[0])\xrightarrow{\sim}j_*\mathbf{det}(j^*\sF[0])\xleftarrow{\sim}j_*\mathbf{det}(\sI[0])\xleftarrow{\sim}j_*\sI.
\]
Here, the first isomorphism follows from the assertion concerning $\mathbf{Rflx}$, and the second isomorphism used the base-change property of $\mathbf{det}$.
This shows that 
\[
   \text{$\sF\simeq \mathbf{det}(\sF[0])$ \q is invertible.}
\]

Finally, since $X\backslash Z$ is schematically-dense in $X$, the injectivity of the restriction map $Y(X)\to Y(X\backslash Z)$ follows from \SP{084N}. To prove the surjectivity, we (locally on $X$) write $Y$ as the relative spectrum of a quasi-coherent $\sO_X$-algebra $\sA$. It remains to observe that every morphism $\sA_U\to \sO_U$ gives rise to a morphism 
\[
   \sA\to j_*\sA_U\to j_*\sO_U\simeq \sO_X. \qedhere
\]
\epf

The following \Cref{normal-s2}  generalizes Serre's conditions $(R_1)+(S_2)$ to the coherent case.
The key ingredient can be traced back to the research of the ``multiplicative theory of ideals'' by M. Zafrullah, \emph{c.f.}~\cite{Zaf78}*{Lemma~7}.

\bthm\label{normal-s2}
For a coherent domain $R$, the following are equivalent:
\benumr
\item\label{normal-s2-i} $R$ is normal;
\item \label{normal-s2-ii} $R$ is the intersection (in $\Frac R$) of its local rings which are valuation rings;
\item\label{normal-s2-iii} every local ring $A$ of $R$ is either a valuation ring, or for any $f\in \fm_A$ we have $\fm_A\not\in \Ass(A/fA)$.
\eenum
 If local rings of $R$ have topologically Noetherian spectra, then these are equivalent to
\begin{equation}
\x{every local ring of $R$ is either a valuation ring or of $\depth\geq 2.$} \label{S2} \tag{$\dagger$}
\end{equation}
\ethm
\bpf
First, for \ref{normal-s2-iii}$\Leftrightarrow$(\ref{S2}) under the Noetherian assumption, note that for every $R$-module $M$ we have $\Ass(M)=\{\fp\in \Spec R\,|\, \depth_{R_{\fp}}M_{\fp}=0\}$, so it suffices to apply \Cref{depth and tau}, since $R$ is a domain.

The implication \ref{normal-s2-ii}$\Rightarrow$\ref{normal-s2-i} is clear, since an intersection of integrally closed subrings is integrally closed.

For \ref{normal-s2-iii}$\Rightarrow$\ref{normal-s2-ii}, {we} need to show that, for nonzero elements $f,g\in R$, if $g\in f R_{\fp}$ whenever $R_{\fp}$ is a valuation ring, $\fp\in \Spec R$, then already $g\in fR$. If not, we consider the annihilator $\Ann\subsetneq R$ of $\bar{g}\in R/fR$ in $R$. By assumption on $f,g$, we have $\Ann_{\fp}=R_{\fp}$ whenever $R_{\fp}$ is a valuation ring. For any minimal element $\fp$ in $V(\Ann) \subset \Spec R$, since $\Ann_{\fp}\neq R_{\fp}$, we deduce that $R_{\fp}$ is not a valuation ring. But since $\fp R_{\fp} \in \Ass(R_{\fp}/fR_{\fp})$, we have obtained a contradiction.

Finally, for \ref{normal-s2-i}$\Rightarrow$\ref{normal-s2-iii}, we may assume that $R$ is a normal, coherent local domain.
It suffices to show that, if there exists an $f\in \fm_R$ such that $\fm_R\in \Ass(R/fR)$, then $R$ is a valuation ring.

Note that since $R$ is coherent, the inverse $J^{-1}$ of every finitely generated ideal $J\subset R$ is again finitely generated; under the normality assumption on $R$, we also have $J\subset (J^{-1})^{-1} \subset R$. Indeed, the first inclusion is clear, and for the second, any $x\in (J^{-1})^{-1}$ satisfies $xJ^{-1}\subset R \subset J^{-1}$ which, by the Cayley--Hamilton, implies that $x$ is integral over $R$. Thus, $x\in R$, by normality.

For every ideal $I\subset R$ define $I^t$ as the union of $(J^{-1})^{-1}$ for all finitely generated subideals $J\subset I$, which is again an ideal of $R$ (unless it is equal to $R$, but this is shown to be impossible below). It is straightforward to check that $I_1^t\subset I_2^t$ whenever $I_1\subset I_2\subset R$, $(aI)^t=aI^t$ for every $a\in \Frac R\backslash \{0\}$, and $I\subset I^t=(I^t)^t$ (e.g., using that any finitely generated ideal of $R$ equals its triple inverse).
\bcl \label{cl-t}
We have $\fm_R=\fm_R^t$. In particular, every subideal $J\subset \fm_R$ satisfies $J^t\subset \fm_R$.
\ecl
\bpf[Proof of the claim]
It suffices to prove the first assertion. Consider the set of all ideals $I\subset R$ satisfying $I^t=I$. This set is stable under taking increasing unions (of ideals), so Zorn's lemma implies that every such $I$ is contained in a maximal one (ordered by inclusion); denote by $\cM$ the subset of such maximal elements. Then we have the following:
\begin{itemize}
\item every ideal $I\subset R$ not contained in any ideal in $\cM$ satisfies $I^t=R$: if not, since the ideal $I^t\subset R$ satisfies $I\subset I^t=(I^t)^t$, it would be contained in some element of $\cM$, a contradiction;
  \item every $\fq \in \cM$ is a prime ideal: if $x,y\in R\backslash \fq$ but $xy\in \fq$, then the inclusion $y(\fq+(x))\subset \fq $ yields $y(\fq+(x))^t\subset \fq^t=\fq$, but $(\fq+(x))^t=R$ by the above, hence $y\in \fq$, a contradiction;
  \item  $R=\bigcap_{\fq\in \cM}R_{\fq}$. It suffices to show that each $a\in \bigcap_{\fq\in \cM}R_{\fq}$ is contained in $R$.
By assumption, we have $a^{-1}R\cap R\not\subset \fq$ for all $\fq\in \cM$, so by the above $R=((a^{-1}R\cap R)^{-1})^{-1}$. Taking a further inverse, we get $R=(a^{-1}R\cap R)^{-1}$, but the latter contains $a$, as desired.
\end{itemize}

Now we will use the crucial assumption that there exists an $f\in \fm_R$ such that $\fm_R\in \Ass(R/fR)$.
This is equivalent to that, there exist $f,g\in R$ such that $g\notin fR$ and $\fm_R$ is the unique prime ideal of $R$ containing $\Ann$, the annihilator of $0\neq \bar{g}\in R/fR$ in $R$. If $\fm_R \notin \cM$ (i.e., $\fm_R\neq \fm_R^t$), then none of $\fq \in \cM$ (which are primes) contains $\Ann$, so that $\Ann_{\fq}=R_{\fq}$, i.e., $g\in fR_{\fq}$, for all $\fq \in \cM$.
But by the intersection formula $R=\bigcap_{\fq\in \cM}R_{\fq}$, we would have $g\in fR$, a contradiction.
\epf

Now we show that $R$ is a valuation ring. By \SP{090Q}, it suffices to show that $I\ce (x,y)$ is principal for arbitrary nonzero elements $x,y\in R$.
Observe that $(I\cdot I^{-1})^{-1}=R$: indeed, every $c\in (I\cdot I^{-1})^{-1}$ satisfies that $c\cdot I\cdot I^{-1}\subset R$, hence $c\cdot I^{-1}\subset I^{-1}$.
By Cayley--Hamilton, $c$ is integral over $R$, but $R$ is normal, thus $c\in R$, as desired.
Consequently, we have $(I\cdot I^{-1})^t=((I\cdot I^{-1})^{-1})^{-1}=R$, which, by \Cref{cl-t}, implies that $I\cdot I^{-1}=R$. Thus, $I$ is invertible and so principal because $R$ is local.
\epf

\bprop\label{rf-s2}
For an integral, locally coherent, and topologically locally Noetherian scheme $X$, a coherent $\sO_X$-module $\sF$, if $X$ is normal, then the following conditions are equivalent:
\benumr
\item\label{eq-i} $\sF$ is reflexive;
\item\label{eq-ii} $\sF$ is torsion-free and $\depth_{\sO_{x}}(\sF_{x})\geq \min(2,\depth_{\sO_X,x}(\sO_{x}))$ {for all $x\in X$};
\item\label{eq-iii} $\sF$ is torsion-free and we have the following equality 
\[
  \tst \text{$\sF=\bigcap_{x\in X_0}\sF_x$\q in the generic stalk $\sF_K\ce \sF\otimes_{\sO_X}\sK_X$}
\]
where $X_0\ce \{x\in X\,|\,\wdim\sO_{X,x}\leq 1\}=\{x\in X\,|\, \x{$\sO_{X,x}$ is a valuation ring}\}$.

\eenum
\eprop
\bpf
The assertions (i)--(iii) are local, we may assume that $X=\Spec A$ is affine and $\sF=\wt{M}$ for a coherent $A$-module $M$.
Thus, $A$ has quasi-compact punctured spectrum.
By \Cref{normal-s2}, every local ring of $A$ is either of depth $\geq 2$ or a valuation ring.
 It follows from Lemmata \ref{lem-rf}~\ref{ddual-tf} and \ref{dep-inh} that \ref{eq-i}$\Rightarrow$\ref{eq-ii}.
Assume that \ref{eq-ii} holds, then the inclusion map $M\ra \bigcap_{\fp\in X_0}M_{\fp}$ is an isomorphism at every prime $\fp\in X_0$.
For any prime ideal $\fp\not\in X_0$, the assumption on depths implies that $\depth_{A_{\fp}}(M_{\fp})\geq 2$, hence \Cref{lem-iso} gives that $M\isoto \bigcap_{\fp\in X_0}M_{\fp}$ is an isomorphism, namely, \ref{eq-iii} holds.
Finally, assume that \ref{eq-iii} holds.
By the assumption on $X$ (or $A$), for every prime $\fp\in X_0$, the local ring $A_{\fp}$ is a valuation ring with topologically Noetherian spectrum.
As a torsion-free finite $A_{\fp}$-module, $M_{\fp}$ is free, so it is reflexive.
Thus, the reflexive hull $\phi_M\colon M\ra M^{\ddual}$ is an isomorphism at all $\fp\in X_0$.
Moreover, since $M^{\ddual}$ is reflexive (\Cref{lem-rf}~\ref{hom-rf}) and \ref{eq-i}$\Rightarrow$\ref{eq-iii} is already proved, we have $M^{\ddual}=\bigcap_{\fp\in X_0}(M^{\ddual})_{\fp}$. Combining this with the assumption $M=\bigcap_{\fp\in X_0}M_{\fp}$, we see that $M^{\ddual}=\bigcap_{\fp\in X_0}(M^{\ddual})_{\fp}=\bigcap_{\fp\in X_0}M_{\fp}=M$ and $\phi_M$ is an isomorphism, that is, $M$ is reflexive.
\epf

\bcor
Let $A$ be a topologically locally Noetherian normal coherent domain with fraction field $K$.
Let $L/K$ be an extension of fields and $B$ the integral closure of $A$ in $L$.
If $B$ is a coherent $A$-module, then $B$ is reflexive over $A$.
\ecor
\bpf
By \Cref{rf-s2}, it suffices to prove that $B=\bigcap_{\fp}B_{\fp}$ where $\fp\subset A$ ranges over prime ideals such that $A_{\fp}$ is a valuation ring.
As $B$ is normal, it remains to show that $B\supset \bigcap_{\fp}B_{\fp}$.
Take $b\in \bigcap_{\fp}B_{\fp}$ and its minimal polynomial $F(x)=x^r+c_{r-1}x^{r-1}+\cdots +c_0\in K[x]$.
For each $\fp\subset A$, the extension of domains $B_{\fp}/A_{\fp}$ is finite, so $b\in B_{\fp}$ satisfies $F(x)\in \Frac(A_{\fp})[x]=K[x]$ and $c_i\in A_{\fp}$.
Consequently, by \Cref{rf-s2} again, we have $c_i\in A$ for each $i$.
In particular, $F(x)\in A[x]$ so $b\in B$, as desired.
\epf
By using the criterion \Cref{rf-s2}, we easily obtain the following result for deducing \Cref{purity for finite flat group schemes}.
\bcor\label{pp-ref}
Let $f\colon X\ra Y$ be a finite, finitely presented, surjective morphism of topologically locally Noetherian, integral, normal, coherent schemes.
Let $\sF$ be a reflexive $\sO_X$-module.
\benumr
\item $f_{\ast}\sF$ is a reflexive $\sO_Y$-module.
\item Let $j\colon Y\hra \ov{Y}$ be an open immersion of integral coherent topologically locally Noetherian schemes.
 If $\depth\sO_{\ov{Y},y}\geq 2$ for all $y\in \ov{Y}\backslash Y$, then $j_{\ast}f_{\ast}\sO_X$ is a reflexive $\sO_{\ov{Y}}$-module, and the morphism
\[
\x{$\ov{f}\colon \un{\Spec}_{\ov{Y}}(j_{\ast}f_{\ast}\sO_X)\ra \ov{Y}$\q is finite.}
\]
In particular, $\un{\Spec}_{\ov{Y}}(j_{\ast}f_{\ast}\sO_X)$ is the relative normalization of $\ov{Y}$ in $X$.
\eenum
\ecor

\section{Geometry of schemes over Pr\"ufer bases}
In this section, we recollect useful geometric properties of scheme over Pr\"ufer bases.

\csub[Geometric properties and reduction methods]
\vskip -0.5cm
\blem\label{geom}
For a Pr\"ufer domain $R$ with spectrum $S$, a finite type irreducible $S$-scheme $X$, a point $x\in X$ and its image $s\in S$, the following assertions hold
\benumr
\item\label{geo-i} all nonempty $S$-fibres have the same dimension;
\item\label{geo-iii} if $X$ is $S$-flat and $X_s$ is generically reduced, then for any maximal point $\xi\in X_s$,  we have an extension of valuation rings
\[
   \text{$\sO_{S,s}\hra \sO_{X,\xi}$ \q inducing an isomorphism of value groups;}
\]
\item\label{geo-iv} for $x\pr\in X$ with $x\pr \neq x$ whose image is denoted by $s\pr$, if $x\in \overline{\{x'\}}$, then
\begin{itemize}
  \item  {either $s=s'$}, and we have $\dim (\sO_{X_{s'},x'}) <\dim(\sO_{X_s,x})$;
  \item or {$s\in \overline{\{s'\}}, s\neq s'$}, and we have $\dim (\sO_{X_{s'},x'}) \le \dim(\sO_{X_s,x})$.
\end{itemize}
\eenum
\elem
\bpf
For \ref{geo-i}, see \cite{EGAIV3}*{Lemme~14.3.10}.
For \ref{geo-iii}, see \cite{MB22}*{Théorème~A}.
Now, to prove \ref{geo-iv}, we may assume that $X$ is affine and has a pure relative dimension, say, $n$, over $R$.
By assumption, we have $s\in \overline{\{s'\}}$. Assume that we are not in the first case. The schematic closure $\overline{\{x'\}}$ is a dominant scheme of finite type over $\overline{\{s'\}}$ (which corresponds to the spectrum of a valuation ring). Therefore, by \ref{geo-i}, all its non-empty fibres have the same dimension. Thus, we deduce from $\overline{\{x'\}} \supset \overline{\{x\}}$ that
\[
\dim (\overline{\{x'\}}_{s'}) = \dim (\overline{\{x'\}}_{s}) \ge \dim (\overline{\{x\}}_{s}).
\]
Hence, we have
\[
\dim (\sO_{X_{s'},x'})= n-\dim (\overline{\{x'\}}_{s'}) \le n-\dim (\overline{\{x\}}_{s}) = \dim(\sO_{X_s,x}). \qedhere
\]
\epf

The following \Cref{enlarge valuation rings} provides us a passage to the case when there is a section.
\blem\label{enlarge valuation rings}
For a valuation ring $V$, an essentially finitely presented (resp., essentially smooth\footnote{By definition, this means that $A$ is a local ring of a smooth $V$-algebra.}) $V$-local algebra $A$, there are an extension of valuation rings $V'/V$ with trivial extension of value groups, and an essentially finitely presented (resp., essentially smooth) $V$-map $V'\to A$ with finite extension of residue fields.
\elem
\bpf
Assume $A=\sO_{X,x}$ for an affine scheme $X$ finitely presented over $V$ and a point $x\in X$ lying over the closed point $s \in \Spec\,V$. Denote $t=\text{tr.deg}(\kappa(x)/\kappa(s))$.
As $\kappa(x)$ is a finite extension of $l\ce \kappa(s)(a_1,\cdots, a_t)$ for a transcendence basis $(a_i)_1^t$ of $\kappa(x)/\kappa(s)$, we have $t=\dim_l\Omega^1_{l/\kappa(s)}\leq \dim_{\kappa(x)}\Omega_{\kappa(x)/\kappa(s)}^1$.
Choose sections $b_1,\cdots,b_t\in \GG(X,\sO_X)$ such that $d\overline{b_1},\cdots,d\overline{b_t} \in \Omega^1_{\kappa(x)/\kappa(s)} $ are linearly independent over $\kappa(x)$, where the bar stands for their images in $\kappa(x)$.
Define $p:X\to \mathbb{A}_V^t$ by sending the standard coordinates $T_1,\cdots, T_t$ of $\mathbb{A}_V^t$ to $b_1,\cdots, b_t$, respectively.
Since $d\overline{b_1},\cdots,d\overline{b_t} \in \Omega^1_{\kappa(x)/\kappa(s)} $ are linearly independent, the image $\eta\ce p(x)$ is the generic point of $\bA^t_{\kappa(s)}$, so $V\pr\ce \sO_{\bA^t_{V},\eta}$ is a valuation ring whose value group is $\GG_{V\pr}\simeq \GG_V$.
Note that $\kappa(x)/\kappa(\eta)$ is finite, the map $V\pr\ra A$ induces a finite residue fields extension.

When $V\ra A$ is essentially smooth, the images of $db_1,\cdots,db_t $ under  the map $\Omega^1_{X/V}\otimes \kappa(x) \to \Omega^1_{\kappa(x)/\kappa(s)}$ are linearly independent, so are their images in $\Omega^1_{X/V}\otimes \kappa(x)$. Hence, $p$ is essentially smooth at $x$.
\epf
In the sequel, \Cref{approxm semi-local Prufer ring} (\emph{cf.}~\cite{GL23}*{Lemma 2.2}), combined with limit arguments, often allows us to only consider Pr\"{u}fer rings of finite Krull dimension.
\blem \label{approxm semi-local Prufer ring}
Every semilocal Pr\"{u}fer domain $R$ is a {filtered union} of its subrings $ R_i$ such that:
\benumr
\item for every $i$, $R_i$ is a semilocal Pr\"{u}fer domain of finite Krull dimension; and
\item for $i$ large enough, $R_i\to R$ induces a bijection on the sets of maximal ideals hence is fpqc.
\eenum
\elem

\bppt[Geometric presentation for the Grothendieck--Serre]
In both Fedorov's and $\check{\mathrm{C}}$esnavi$\check{\mathrm{c}}$ius' works on mixed characteristic Grothendieck--Serre, significant emphasis is placed on geometric results of a certain type, reminiscent of Gabber--Quillen, as demonstrated in \cite{Fed22b}*{Proposition~3.18} and \cite{Ces22a}*{Variant~3.7}, respectively.
Similarly, in our context, we observe an analogous result to \emph{loc.~cit.} and record it below.
\eppt

\blem \label{Ces's Variant 3.7}
Let
\benumr
  \item $R$ be a semilocal Pr\"{u}fer ring;
  \item $\overline{X}$ be a projective, flat $R$-scheme with fibres of pure dimension $d>0$;
  \item  $X\subset \overline{X}$ be an open subscheme, smooth over $R$;
  \item  $\textbf{x}\subset X$ be a finite subset;
  \item $Y\subset X$ be a closed subscheme which is
   $R$-fibrewise of codimension $\ge 1$; assume also that $\overline{Y}\backslash Y$ is $R$-fibrewise of codimension $\ge 2$.
\eenum
Then, there are
\benumr
\item an affine open $S\subset \mathbb{A}_R^{d-1}$ and an affine open neighbourhood $ U\subset X$ of $\textbf{x}$, and
\item a smooth morphism $\pi:U\to S$ of relative dimension 1
\eenum
such that $Y \cap U$ is $S$-finite.
\elem
\bpf
This can be proved similarly as \cite{Ces22a}*{Variant~3.7}.
\epf

\csub[Regularity and Reflexive sheaves over Pr\"ufer bases]\label{subsec-Regularity and Reflexive sheaves over Pr\"ufer bases}

Now, we turn to the special case of schemes over Pr\"ufer bases.
To begin with, we consider the coherence and calculation of depths over Pr\"ufer bases, as the following Lemmata \ref{coh of V-flat ft alg} and \ref{depth of O_{X,x}}.

\blem\label{coh of V-flat ft alg}
 Let $X$ be a scheme that is flat and locally of finite type over a Pr\"ufer domain $R$.
\benumr
  \item \label{coherence of X}$X$ is locally of finite presentation and locally coherent.
  \item \label{coherence of O}For every point $x\in X$, the local ring $\sO_{X,x}$ is coherent.
  \item \label{G-R 17.4.1 bound}If $x\in X$ lies over $s \in \Spec R$ and $\sO_{X_s,x}$ is regular, then every finitely presented $\sO_{X,x}$-module has a finite resolution by finite free $\sO_{X,x}$-modules of length $\le \dim \sO_{X_s,x}+1$. In particular,
  \[
  \x{$\sO_{X,x}$ is coherent regular.}
  \]
  \item \label{normality of O} The local ring $\sO_{X,x}$ in \ref{G-R 17.4.1 bound} is a normal domain.
 \eenum
\elem
\bpf
We may assume that $X=\Spec A$ is affine for a finite type flat $R$-algebra $A$.

For \ref{coherence of X}, by \cite{RG71}*{Corollaire~3.4.7} ({or \cite{Nag66}*{Theorem 3'} when $R$ is a valuation ring}), every finitely generated, flat algebra over a domain is finitely presented.
So $A$ is finitely presented over $R$. The coherence of $A$ thus follows from the following facts: 1) any polynomial ring over a Pr\"ufer domain $R$ is coherent (\cite{Gla89}*{Corollary 7.3.4}); 2) the quotient of a coherent ring by a finitely generated ideal is again coherent.

The assertion \ref{coherence of O} follows because each local ring of a coherent ring is coherent.
For \ref{G-R 17.4.1 bound}, see \cite{GR18}*{Proposition~11.4.1}.
Finally, \ref{normality of O} follows from a more general result \Cref{ring property coh reg}~\ref{ring property coh reg-i}.
\epf

By the above \Cref{coh of V-flat ft alg}, a flat, locally of finite type scheme over a valuation ring is locally coherent.
Further, when the fibres are Cohen-Macaulay, the depth of its local ring can be computed as follows.
\blem \label{depth of O_{X,x}}
Let $R$ be a Pr\"ufer domain and let $X$ be an $R$-flat finite type scheme.
Let $x\in X$ be a point with image $s\in \Spec R$ and local ring $A\ce \sO_{X,x}$.
If the following conditions hold
\benum
\item $s$ is not the generic point of $\Spec R$;
\item $\sO_{X_s,x}$ is Cohen-Macaulay; and
\item $A$ has quasi-compact punctured spectrum, or equivalently, so does $\sO_{\Spec R,s}$,
\eenum
{then, $A$ has a regular sequence in $\fm_A$ of length $d+1$ for $d=\dim \sO_{X_s,x}$.} In particular, 
\benumr
      \item \label{depth of A} $\depth_A(A)=d+1$; and
      \item $\tau_N(A)\ge d+1$ for any finitely generated $A$-module $N$ supported on $\{x\}$, 
       \[
       \x{therefore,\q $\Ext_A^i(N,A)=0$\q for all $i\le d$.}
       \]
\eenum
\elem
     \bpf
We may assume that $R$ is a valuation ring $V$ with quasi-compact punctured spectrum $U_V$ such that $x$ is in the closed fibre.
     In view of \Cref{depth and tau}, it suffices to show that $\fm_A$ contains a regular sequence of $A$ whose common vanishing locus is zero-dimensional.
     As $U_V$ is quasi-compact, we can pick an $f\in \mathfrak{m}_V$ such that $\dim (V/fV) =0$.
       Let $(g_1,\cdots, g_d)$ be a sequence in $\fm_A$ such that their images in the local ring $A/\mathfrak{m}_VA$ forms a regular sequence.
     By the flatness criterion \cite{EGAIV3}*{Théorème~11.3.8}, $(g_1,\cdots, g_d)$ is a regular sequence of $A$, and the quotient ring $\overline{A}\ce A/(g_1,\cdots, g_d)$ is $V$-flat with maximal ideal $\mathfrak{m}_V\overline{A}= \mathfrak{m}_{\overline{A}}$. Therefore, $(g_1,\cdots, g_d,f)$ is a regular sequence of $A$ for which $\dim A/(fA+\sum_{i=1}^dg_iA) =0$, as desired.
     \epf

Let $R$ be a Prüfer domain, $X$ a scheme flat and of finite type over $R$, and $A$ the local ring of $X$ at a point $x\in X$.
So $A$ is coherent.
Knaf proved \cite{Kna08}*{Theorem~1.1} that $A$ is coherent regular if and only if $\wdim(A)$ is finite; moreover, in this coherent regular case we have (\emph{cf.}~\Cref{fibre criterion coh reg}~(ii))
\[
\wdim(A)=\dim (A\otimes_R \kappa(\fq))+\wdim(R_{\fq}),
\]
where $\fq\subset R$ lies below $x$ and we have $\wdim(R_{\fq})=0$ whenever $x$ lies over $\Frac R$ otherwise $\wdim(R_{\fq})=1$.
Besides, if $A$ is coherent regular, then $A\otimes_R \kappa(\fq)$ is Cohen-Macaulay and even Noetherian regular if $R_{\fq}$, or equivalently, $A$ is non-Noetherian (\cite{Kna08}*{Theorem~1.3}). Conversely, we see that the regularity of $A\otimes_R \kappa(\fq)$ implies the regularity of $A$ (\emph{cf.} \Cref{coh of V-flat ft alg}~\ref{G-R 17.4.1 bound}).
We formulate these into the following.

\bprop\label{prufer-regularity}
Let $R$ be a Prüfer domain and $(A,\fm_A)$ a local ring that is $R$-flat and essentially finitely presented.
Denote $\fq\subset R$ the prime underlies $\fm_A$.
\benumr
\item\label{prufer-regularity-i} If $\fm_A$ is a minimal prime of $A\otimes_{R_{\fq}} \kappa(\fq)$, then $A$ is regular if and only if $A$ is a valuation ring.
\item\label{prufer-regularity-ii} If $A$ is non-Noetherian, then $A$ is regular if and only if $A\otimes_{R_{\fq}} \kappa(\fq)$ is regular,
in which case we have $\wdim(A)=\dim(A\otimes_{R}\kappa(\fq))+\wdim(R_{\fq})$.
\item\label{prufer-regularity-iii} Assume that $A$ is regular with Noetherian spectrum. Define $X_0\subset \Spec A$ as
\[
   X_0\ce \{\x{height-one primes of $A\otimes \Frac R$}\}\cup \{\x{minimal primes of $R$-fibres of $A$} \}
\]
Then $X_0=\{x\in \Spec A\,|\, \x{$\sO_{\Spec A,x}$ is a valuation ring}\}$, and all local rings at $x\in (\Spec A)\backslash X_0$ {have} $\depth \geq 2$, or equivalently, have weak dimensions $\geq 2$.
\eenum
\eprop
\bpf
It suffices to deal with \ref{prufer-regularity-iii}.
By \Cref{geom}~\ref{geo-iii}, $X_0$ is the set of all points of $X$ where the local rings are valuation rings.
The regularity of $A$ and \Cref{ring property coh reg}~\ref{ring property coh reg-i} imply that $A$ is normal, thus by \Cref{normal-s2}, all local rings at $x\in (\Spec A)\backslash X_0$ have depth $\geq 2$.
Finally, by \Cref{wdim=depth}, the regularity of $A$ yields $\depth=\wdim$ locally on $\Spec A$, thus all local rings at $x\in (\Spec A)\backslash X_0$ have $\wdim\geq 2$.
\epf

Finally, we present a Pr\"uferian variant of \Cref{Rflx-Pic-Vect on reg}, the purities of reflexive modules, of line bundles, and of vector bundles.
Note that our base now is a general Pr\"ufer domain so the local rings under consideration may have non-quasi-compact punctured spectra.
The crux is to carefully carry out a limit argument, which preserves the fibrewise codimensions of closed subsets.

\blem\label{lim-codim}
Let $S$ be a semilocal, affine, integral Pr\"ufer scheme and $\eta$ its generic point.
Given
\benumr
\item a flat, surjective, {finitely presented} morphism $f\colon X\ra S$  with regular fibres (resp., $f$ is smooth);
\item a coherent $\sO_X$-module $\sF$ that is reflexive at a point $x\in X$; and
\item a constructible closed $Z\subset X$ such that 
\[
   \text{$\codim(Z_s, X_s)\geq 1$ for every $s\in S$ \q and\q  $\codim(Z_{\eta},X_{\eta})\geq 2$.}
\]
\eenum
Then, there are
\benuma
\item a semilocal, affine, integral Pr\"ufer scheme $S_0$ of finite dimension with generic point $\eta_0$;
\item a flat, surjective, finite type morphism $f_0\colon X_{0}\ra S_{0}$ with regular fibres (resp., $f_0$ is smooth) such that $X_{0}\times_{S_{0}}S\simeq X$;
\item a coherent $\sO_{X_{0}}$-module $\sF_{0}$ as the inverse image of $\sF$ and is reflexive at the image $x_0$ of $x$; and
\item a constructible closed subset $Z_{0}\subset X_{0}$ such that $Z_{{0}}\times_{S_{0}}S\simeq Z$ and
\[
\q \text{$\codim(({Z_{0}})_s, ({X_{0}})_s)\geq 1$ for every $s\in S_{0}$\, and \, $\codim(({Z_{0}})_{\eta_0},({X_{0}})_{\eta_0})\geq 2$.}
\]
\eenum
\elem
\bpf

We apply \Cref{approxm semi-local Prufer ring} to $S$ then use limit arguments in \cite{EGAIV3}*{\S8}.
The condition that $X$ has regular $S$-fibres descends to $X_{\lambda_0}$ by \cite{EGAIV2}*{Proposition~6.5.3} (resp., the smoothness of $f$ descends by \cite{EGAIV4}*{Proposition~17.7.8}).
The reflexive $\sO_X$-module $\sF$ descends thanks to \cite{EGAIV3}*{Théorème~8.5.2} and by applying \cite{EGAIV3}*{Corollaire~8.5.2.5} to $\sF\isoto \sF^{\ddual}$.
Because $Z$ is contructible closed, by \cite{EGAIV3}*{Théorème~8.3.11}, it descends to $Z_{\gL}$ such that $p_{\gL}^{-1}(Z_{\gL})=Z$. For $f_{\gL}\colon X_{\gL}\ra S_{\gL}$, by the transversity of fibres and \cite{EGAIV2}*{Corollaire~4.2.6}, $Z_{\gL}$ does not contain any irreducible components of $f_{\gL}^{-1}(s_{\gL})$ for any $s_{\gL}\in S_{\gL}$.
  Finally, the image of the generic point $\eta\in S$ is the generic point $\eta_{\gL}\in S_{\gL}$. By \cite{EGAIV2}*{Corollaire~6.1.4}, we have $\codim((Z_{\gL})_{\eta_{\gL}}, (X_{\gL})_{\eta_{\gL}})=\codim(Z_{\eta}, X_{\eta})\geq 2$.
\epf

\bthm\label{equiv-cats}
Let $R$ be a semilocal Pr\"ufer domain with spectrum $S$ and $f\colon X\ra S$ a flat, {locally} of finite type morphism of {schemes} with regular fibres.
Let $Z\subset X$ be a constructible closed subset such that
\[
   \begin{cases}
    \codim(Z_s,X_s)\geq 1 & \text{for all $s\in S$, and}\\
    \codim(Z_{\eta},X_{\eta})\geq 2 & \text{for the generic point $\eta\in S$}.
   \end{cases}
\]
For the open immersion $j\colon X\backslash Z\hra X$, the restriction $j^{\ast}$ and pushforward $j_{\ast}$ induce
\[
   \x{the equivalences \q $\sO_X\x{-}\mathbf{Rflx}\isoto \sO_{X\backslash Z}\x{-}\mathbf{Rflx}$\q and \q $\mathbf{Pic}\,X\isoto \mathbf{Pic}\,X\backslash Z$.}
\]
In particular, {for} every scheme $Y$ affine over $X$, we have a bijection of sets
\[
   Y(X)\isoto Y(X\backslash Z).
\]
\ethm
\bpf
The problem is Zariski-local on $X$, so we may assume that $f$ is of finite type. By a limit argument involving \Cref{lim-codim}, we may further assume that $|S|$ is finite.
Then, $\abs{X}$ is the finite union of its $S$-fibres $\abs{X_s}$, which are Noetherian spaces, so $X$ is topologically Noetherian.
By \Cref{coh of V-flat ft alg}~\ref{G-R 17.4.1 bound}, $X$ is coherent regular (see also \Cref{reg-reg-reg}).
Since the local rings of $X$ are normal (\Cref{ring property coh reg}~\ref{ring property coh reg-i}) and the generic fibre $X_{\eta}$ is regular and schematic dense in $X$, $X$ is Zariski locally an integral scheme.
The assumption on $Z$ implies that $X\backslash Z$ contains $X_0\subset X$ defined in \Cref{prufer-regularity}~\ref{prufer-regularity-iii}, so for every $z\in Z$, we have $\wdim \sO_{X,z}\geq 2$.
Hence, the assertion follows from \Cref{Rflx-Pic-Vect on reg}.
\epf

\section{Auslander's flatness criterion}
\numberwithin{equation}{subsection}
The goal is to establish \Cref{Auslander's flatness criterion} as a counterpart of  Auslander's flatness criterion \cite{Aus62}*{Theorem~1.3} on schemes smooth over valuation rings.
As expected, our criterion leads to the analog of Zariski--Nagata purity (\Cref{fet cover}).
\begin{thm-tweak} \label{Auslander's flatness criterion}
    Let $V$ be a valuation ring, $X$ a $V$-smooth scheme, and $x\in X$ a point. Set $A\ce \sO_{X,x}$. Let $M$ be a reflexive $A$-module. If $\End_A(M)$ is isomorphic to a direct sum of copies of $M$, then $M$ is $A$-free.
\end{thm-tweak}

Similar to Auslander's proof, our strategy relies on an estimate of the length of cohomology groups of $M$.
To begin with, we introduce the length function on torsion modules over valuation rings.

\bpp[Lengths of torsion modules] \label{delta}
Let $V$ be a valuation ring, distinct from a field, with fraction field $K$, and a valuation map $\nu\colon K\ra \GG$. Every finitely presented torsion $V$-module $M$ can be expressed as
\[
\tst \text{$M\simeq \bigoplus_{i}V/a_i V$\q for finitely many $a_i\in V\backslash \{0\}$.}
\]
Define the \emph{length} of $M$ as $ \delta(M)=\sum_i\nu(a_i) \in \Gamma_{\ge 0}$.
The element $\delta(M)$ is well defined, and $\delta(M)=0$ if and only if $M=0$.
Every acyclic, bounded complex $M^{\bullet}$ of torsion, finitely presented $V$-modules satisfies
\[
  \tst\sum_j (-1)^j \delta(M^j)=0.
\]
\epp

\begin{lemma-tweak} \label{the map l}
Let $V$ be a valuation ring that is not a field. Let $(V,\fm_V)\to (A, \mathfrak{m}_A)$ be an essentially smooth, local map of local rings. Denote by $A\x{-}\mathbf{Mod}_{\mathrm{tor, fp}}$ the collection of all finitely presented $A$-modules $M$ supported on $\{\fm_A\}$. Then there exist a totally ordered abelian group $\Gamma$ and a map $l\colon A\x{-}\mathbf{Mod}_{\mathrm{tor, fp}}\ra \GG_{\geq 0}$ satisfying the following properties:
\begin{itemize}
\item for $A$-module $M\in A\x{-}\mathbf{Mod}_{\mathrm{tor, fp}}$, we have $l(M)=0$ if and only if $M=0$;
\item for every acyclic, bounded complex $M^{\bullet}$ such that $M^j\in A\x{-}\mathbf{Mod}_{\mathrm{tor, fp}}$ for each $j$, one has
    \[
    \tst \sum_j (-1)^j l(M^j)=0.
    \]
\end{itemize}
\end{lemma-tweak}

{It is worth mentioning that the set $A\x{-}\mathbf{Mod}_{\mathrm{tor, fp}}$ might be trivial, consisting only of the zero module. In fact, this occurs precisely when $\Spec V\backslash \{\fm_V\}$ is not quasi-compact.
This totally ordered abelian group $\GG$ is a value group, and in the sequel, we will only use its property of being partially ordered.}
\bpf
First we assume that the structural map $V\to A$ admits a retraction $A \twoheadrightarrow V$. In this case we claim that $M$ is finitely presented over $V$ and is $V$-torsion. So we can simply let $\Gamma$ be the valuation group of $V$ and set $l(M)\ce \delta(M)$, where $\delta$ is delivered from \S\S\Cref{delta}.
Indeed, it is clear that $M$ is $V$-torsion.
Any section $\Spec V\ra \Spec A$ is a regular immersion \SP{067R}, so there is a finitely generated ideal $J\subset A$ such that $V\simeq A/J$.
Hence, since $M\in A\x{-}\mathbf{Mod}_{\mathrm{tor, fp}}$, we see that $J^nM=0$ for a  large $n$. On the other hand, the essential smoothness of $A$ over $V$ implies that $J/J^2$ is a free $V\simeq A/J$-module whose rank equals the rank of the free $A$-module $\Omega^1_{A/V}$, and there is a natural isomorphism of graded $V\simeq A/J$-algebras
\[
\tst \bigoplus_{n\ge 0} J^n/J^{n+1} \simeq \text{Sym}^{\bullet}_{A/J}(J/J^2).
\]
In particular, $A/J^n$ is a finite free $V$-module for every $n\ge 1$. 
Therefore, by tensoring 
\[
   \text{a presentation\q  $A^N\to A^{N'} \to M\to 0$\q of $M$}
\]  
 with $A/J^n$ for a large enough $n$, we get a desired finite presentation of $M$ over $V$.

In the general case, we first use \Cref{enlarge valuation rings} to reduce to the case when the residue fields extension of $V\to A$ is finite. Then, if $B$ is the integral closure of $V$ in an algebraic closure of $\Frac V$, we let $V\pr$ be a valuation ring of $\Frac (B)$ centered at a maximal ideal of $B$. It is clear that $V\pr$ is absolutely integral closed, so it is strictly Henselian and there exists a $V$-map $\phi\colon A/\mathfrak{m}_{A}\to V\pr/\mathfrak{m}_{V\pr}$. Let $A\pr\ce A\otimes_VV\pr$. Then $\phi$ induces a $V'$-map $\phi\pr \colon A\pr\to V\pr/\mathfrak{m}_{V\pr}$; let $\mathfrak{p}\subset A\pr$ be its kernel. Then $A\pr_{\mathfrak{p}}$ is essentially smooth over $V\pr$ and $\phi'$ induces a $V\pr$-map $A\pr_{\mathfrak{p}}\to V\pr/\mathfrak{m}_{V\pr}$, which, by the Henselianity of $V\pr$, lifts to a $V\pr$-map $A\pr_{\mathfrak{p}}\to V\pr$. By the previous paragraph, the lemma is true for $A\pr_{\mathfrak{p}}$, say, with corresponding map $l\pr$ valued in $\Gamma$, where $\Gamma$ is the valuation group of $V\pr$. Since $A\to A\pr_{\mathfrak{p}}$ is faithfully flat, it suffices to define 
\[
   l(M)\ce l\pr(M\otimes_AA\pr_{\mathfrak{p}}). \qedhere
 \]
\epf

\bpp[Homological algebra Lemmata]
\epp

\begin{lemma-tweak} \label{lem on isom of Ext's}
  Let $(A,\fm_A)$ be a coherent local ring with a regular sequence of length $d\geq 1$. Let $M\overset{\phi}{\to} N$ be a morphism of coherent $A$-modules that induces an isomorphism over $\Spec A\backslash \{\fm_A\}$.
  We have isomorphisms $\Ext_{A}^i(N,A)\xrightarrow{\sim} \Ext_{A}^i(M,A)$ for all $i<d-1$ and a monomorphism $\Ext_{A}^{d-1}(N,A) \hookrightarrow \Ext_{A}^{d-1}(M,A)$.
\end{lemma-tweak}
This will be applied to \Cref{Auslander's flatness criterion} with $A\ce \sO_{X,x}$, $d=\dim(\sO_{X_s,x})+1$, and $s$ is not the generic point.
\bpf
 By \Cref{lem-coh}~\ref{ker-coker-coh}, $\Ker \phi$ and $\Coker \phi$ are coherent $A$-modules supported on $\{\fm_A\}$. By assumption and \Cref{depth and tau}, we have $\tau_{\Ker \phi}(A)\ge d$ and $\tau_{\Coker \phi}(A)\ge d$ (see \S\S ~\ref{setup of depth & proj dim} for the definition of $\tau_{-}(-)$). Consider the following short exact sequences
\[
\x{$0\to \Ker \phi \to M \to \im \phi \to 0$\q and \q $0\to \im \phi \to N \to \Coker \phi \to 0$.}
\]
By applying $\Hom_{A}(-,A)$, we get two long exact sequences concerning Ext's, and the lemma follows from $\Ext_{A}^i(\Ker \phi,A)=0$
and $\Ext_{A}^i(\Coker \phi,A)=0$ for $i< d$.
\epf

\begin{lemma-tweak}\label{coherence of Tor and Ext}
Let $V$ be a valuation ring, $X$ a $V$-smooth finite type scheme, and $x\in X$ a point that lies over a non-generic point $s\in \Spec V$. For finitely presented $A\ce \sO_{X,x}$-modules $M$ and $N$, $\Ext_A^i(M,N)$ and $\Tor_i^A(M,N)$ are finitely presented over $A$ for all $i\ge 0$ and are zero for $i> d+1$, where $d=\dim \sO_{X_s,x}$.
\end{lemma-tweak}
\bpf
By \Cref{coh of V-flat ft alg}~\ref{G-R 17.4.1 bound}, $M$ has a resolution by finite free $A$-modules of length $\le d+1$: $F_{\bullet}\to M$, $F_i=0$ for $i>d+1$.
Then the $A$-modules
\[
\Ext_A^i(M,N)=H^i(\Hom(F_{\bullet},N))\q \x{and} \q \Tor^A_i(M,N)=H_i(F_{\bullet}\otimes N)
\]
are all coherent, i.e., finitely presented $A$-modules, and are zero for $i>d+1$.
\epf

\blemt \label{lem on dual and hom}
Let $A$ be a coherent domain with a topological Noetherian spectrum. If $A$ is normal, then for every finitely presented $A$-module $M$, we have 
\[
\text{a natural isomorphism \q $\mathrm{End}_A(M)^{\ddual} \isoto \mathrm{End}_A(M^{\ddual})$.}
\]
\elemt
\bpf
By the functoriality of $(-)^{\ddual}$, there is a natural homomorphism
\[
\mathrm{End}_A(M) \to \mathrm{End}_A(M^{\ddual}).
\]
As the target is reflexive by \Cref{lem-rf}~\ref{hom-rf}, this map factors through $\mathrm{End}_A(M)^{\ddual}$, thus yielding a natural map $\mathrm{End}_A(M)^{\ddual} \to \mathrm{End}_A(M^{\ddual})$.
Since both the source and target of this map are reflexive, by \Cref{normal-s2} and \Cref{rf-s2}, it is enough to check this map is an isomorphism when $A$ is a valuation ring.
But if this is the case, then there are an $N\in \b{Z}_{\geq 0}$ and finitely many $a_i\in \fm_A\backslash \{0\}$ for which
\[
  \tst M\simeq  A^{\oplus N} \bigoplus \,\p{\oplus_i A/a_i A}.
\]
Consequently, we conclude by the following isomorphisms 
\[
   \End_A(M)^{\ddual} \simeq \End_A(A^{\oplus N}) \simeq  \End_A(M^{\ddual}). \qedhere
\]
\epf

\bpp[Proof of \Cref{Auslander's flatness criterion}]\label{proof of aus-flat}
The proof proceeds as the following steps.
\epp
\textbf{Preliminary cases and reductions.}
Firstly, since $X$ is locally of finite presentation over $S$ and $M$ is finitely presented over $A$, by a standard limit argument involving \Cref{approxm semi-local Prufer ring}, we are reduced to the case when $V$ is a finite-rank valuation ring. Secondly, if $V'$ is a valuation ring of an algebraic closure of $\Frac (V)$ that dominates $V$ and if $x'\in X'\ce X\times_VV'$ is a point lying over $x\in X$, then $M_{A'}\ce M\otimes_A A'$ is a finitely presented reflexive $A'$-module and $\End_{A'}(M_{A'}) \simeq \End_{A}(M)\otimes_AA'$ is isomorphic to a direct sum of copies of $M_{A'}$, where $A'\ce \sO_{X',x'}$ (because $A'$ is faithfully flat over $A$).
By faithfully flat descent \SPD{08XD}{00NX}, the freeness of $M$ over $A$ is equivalent to the freeness of $M_{A'}$ over $A'$. Therefore, by replacing $V$ by $V'$, $A$ by $A'$, and $M$ by $M_{A'}$, we are reduced to the case where $\Frac (V)$ is algebraically closed.
{It is important to emphasize that the smoothness of
$X$ plays a crucial role in this reduction, as the regularity of fibres is not maintained under base changes.
Furthermore, while the assumption that $\Frac V$ is algebraically closed will be invoked only towards the conclusion of the proof, the true necessity lies in ensuring the perfectness of all residue fields of $V$.}

Let $s\in \Spec V$ be the image of $x$. Set $d_x\ce \dim \sO_{X_s,x}$ and $r\ce \dim V$. The case $r=0$ and $d_x$ arbitrary is classical. The case $r$ arbitrary and $d_x=0$ is trivial, where $A$ is a valuation ring (\Cref{geom}~\ref{geo-iii}). The case $r$ arbitrary and $d_x=1$ follows from \Cref{wdim2-vect}~\ref{wdim2-vect-ii}.
Thus, we may assume $d_x\ge 2$ in the sequel.

\textbf{Case 1: $r$ is arbitrary and $d_x=2$.}
  We will proceed by induction on $r$. The induction hypothesis is that the assertion holds for $d_x=2$ and $r\pr\leq r-1$. Notice that, for any proper generalization $x'\in X$ of $x$ that lies over, say, $s'\in \Spec V$, by \Cref{geom}~\ref{geo-iv}, we have either $s'=s$ and $d_{x'}<2$, or the height of $s\pr$ is less than $r$ and $d_{x'}\le 2$.
 Hence, by induction hypothesis and the preliminary cases above, the assertion holds for $\sO_{X,x'}$.
 As $M_{x'}$ is a finitely presented reflexive $\sO_{X,x'}$-module and
\[
  \End_{\sO_{X,x'}}(M_{x'})=\End_{\sO_{X,x}}(M)\otimes_{\sO_{X,x}}\sO_{X,x\pr} \simeq (\bigoplus M)\otimes_{\sO_{X,x}}\sO_{X,x\pr}=\bigoplus M_{x\pr},
\]
the induction hypothesis applies to the $\sO_{X,x\pr}$-module $M_{x\pr}$, implying that $M_{x'}$ is $\sO_{X,x'}$-free. In other words, $\widetilde{M}$ is locally free over $\Spec A \backslash \{x\}$. Consider the following  
\[
\text{evaluation map \q $M^{\vee}\otimes_A M\ra \Hom_A(M,M), \q f\otimes m\mapsto [m\pr \mapsto f(m\pr )m]$.}
\]
By the local freeness of $\widetilde{M}$ over $\Spec A\backslash \{x\}$, it is an isomorphism over $\Spec A\backslash \{x\}$.
Since $d_x=2>1$, by \Cref{lem on isom of Ext's}, we apply $\Ext_A^1(-,A)$ to the above map to obtain
\begin{equation}\label{extmma-equal-extendma}
\Ext^1_A(M^{\vee}\otimes M, A) \simeq \Ext_A^1(\End_A(M),A)\simeq \Ext^1_A(M,A)^{\oplus \mathrm{rk}_M}
\end{equation}
as isomorphisms of $A$-modules that are supported on $\{x\}$ by the local freeness of $\wt{M}$ over $\Spec A\backslash \{x\}$, where $\mathrm{rk}_M=\dim_{\Frac A}M\otimes_A\Frac A$.
By \Cref{coherence of Tor and Ext}, the modules in (\ref{extmma-equal-extendma}) are also finitely presented over $A$.

For the adjunction $\Hom_A(M,\Hom_A(M^{\vee},-))\simeq \Hom_A(M\otimes M^{\vee},-)$, we take their derived functors valued at $A$, so the $E_2$-page of the associated Grothendieck spectral sequence yields a monomorphism
\[
 \Ext_A^1(M,M) \hra \Ext^1_A(M\otimes M^{\vee},A)\overset{(\ref{extmma-equal-extendma})}{\simeq} \Ext^1_A(M,A)^{\oplus \mathrm{rk}_M},
\]
where we have used $M^{\ddual}\simeq M$; again, by the local freeness of $\widetilde{M}$ over $\Spec A\backslash \{x\}$ and \Cref{coherence of Tor and Ext}, they are finitely presented supported on $\{x\}$. In particular, the map $l$ from \Cref{the map l} applies so we have
\begin{equation}\label{length-extmm-extmma}
l(\Ext_A^1(M,M)) \le \mathrm{rk}_M\cdot l(\Ext_A^1(M,A)).
\end{equation}
Since $M$ is reflexive, by \Cref{wdim2-vect}~\ref{wdim2-vect-i}, we have $\pd_A M\le d_x-1=1$.
We prove $\pd(M)=0$ by contradiction.
If $\pd(M)= 1$, then $M$ has a free resolution $0\to F_1\to F_0\to M\to 0$ by finite $A$-modules.
This sequence is nonsplit, corresponding to a nontrivial extension class in
\[
  \Ext^1_A(M,F_1)\simeq \Ext^1_A(M,A)^{\mathrm{rank}(F_1)}.
\]
In particular, $C\ce \Ext_A^1(M,A) \neq 0$.  Applying $\Hom_A(-,A)$ to $F_{\bullet}\ra M$ yields 
\[
   \text{an exact sequence \q $0\ra M^{\vee} \ra F_0^{\vee} \ra F_1^{\vee} \ra \Ext_A^1(M,A) \to 0$.}
\]
Tensoring it with $M$, we get an exact sequence
\[
F_0^{\vee}\otimes_A M \ra F_1^{\vee}\otimes_A M \ra \Ext_A^1(M,A) \otimes_A M \to 0.
\]
Since
\begin{equation*}
\begin{split}
\Coker\left(F_0^{\vee}\otimes M \to F_1^{\vee}\otimes M\right)  & \simeq  \Coker\left(\Hom_A(F_0,M)\to \Hom_A(F_1,M)\right)  \\
 & = \Ext_A^1(M,M),
\end{split}
\end{equation*}
we deduce that $\Ext_A^1(M,M) \simeq \Ext_A^1(M,A) \otimes_A M=C\otimes_A M$.

By tensoring $0\ra F_1\ra F_0\to M\to 0$ with $C=\Ext_A^1(M,A)$ (which is nonzero, finitely presented, and supported at $\{x\}$, by the local freeness of $\wt{M}$ over $\Spec A\backslash \{x\}$), we get an exact sequence of finitely presented $A$-modules supported on $\{x\}$:
\[
0\ra \Tor^A_1(C,M) \ra C\otimes_A F_1 \ra C\otimes_A F_0 \to C\otimes_A M \to 0.
\]
Denote $\mathrm{rk}_M\ce \mathrm{rank} F_0-\mathrm{rank} F_1 >0$.
Applying the map $l$ from \Cref{the map l}, we get
\begin{equation}\label{length-cm}
\begin{split}
l(C\otimes_A M) & = l(C\otimes_A F_0)-l(C\otimes_A F_1)+ l(\Tor^A_1(C,M) \\ 
& = \mathrm{rk}_M\cdot l(C)+l(\Tor^A_1(C,M)).
\end{split}
\end{equation}
On the other hand, since $C\otimes_A M\simeq \Ext^1_A(M,M) $, we deduce the following inequality
\begin{equation}\label{length-cm2}
l(C\otimes_A M)  \overset{(\ref{length-extmm-extmma})}\le  \mathrm{rk}_M\cdot  l(C).
\end{equation}
The combination of (\ref{length-cm}) and (\ref{length-cm2}) leads to $l(\Tor_1^A(C,M))=0$.
So, we obtain
\[
   \text{a short exact sequence \q $0\ra C\otimes_A F_1 \ra C\otimes_A F_0 \ra C\otimes_A M \to 0$,}
\]
which combined with \Cref{a lem on linear algebra} implies that the map $F_1\ra F_0$ splits, that is, $M$ is $A$-free, contradicting our assumption that $\pd_A(M)=1$. This completes the case when $r$ is arbitrary and $d_x=2$.

\textbf{Case 2: $r$ is arbitrary and $d_x>2$.}
We will proceed by double induction on the pair $(r=\text{ht}(s),d_x)$. By induction hypothesis, the assertion holds for {all smooth $V$-schemes} $X\pr$ and all points $x\pr\in X\pr$ such that $\text{ht}(s\pr)\le \text{ht}(s)$ and $d_{x\pr} \le d_x$, where $s\pr\in \Spec V$ lies below $x\pr$, and at least one of equalities is strict. 
In particular, by \Cref{geom}~\ref{geo-iv}, the induction hypothesis applies to $\sO_{X,x\pr}$ for all proper generalization $x\pr\in X$ of $x$.
As $M_{x\pr}$ is finitely presented reflexive over $\sO_{X,x\pr}$ and $$\End_{\sO_{X,x\pr}}(M_{x\pr})=\End_{\sO_{X,x}}(M)_{x\pr} \simeq \bigoplus M_{x\pr},$$ the induction hypothesis gives that $M_{x\pr}$ is $\sO_{X,x\pr}$-free. 
In other words, $\widetilde{M}$ is locally free over $\Spec A \backslash \{x\}$.

\bcl[\SP{057F}]
Assume that the residue field extension of $V\to A$ is separable (e.g., this holds if $\kappa(s)\ce V/\mathfrak{m}_V$ is perfect),
 then there exists an $a \in A$ such that $\overline{A}\ce A/(a)$ is essentially $V$-smooth and $$\dim(\overline{A}/\mathfrak{m}_V \overline{A})=d_x-1.$$
 \ecl

Since our $V$ has algebraically closed fraction field (by the first paragraph), all of its primes have algebraically closed residue fields, so we can choose $a\in A$ as in the above claim. Since $a$ is a nonzerodivisor in $A$ and $M=\Hom_A(M^{\vee},A)$, we see that $a$ is $M$-regular. Set $\overline{M}\ce M/aM$. Applying $\Hom _A(M,-)$ to the short exact sequence $0\to M \xrightarrow{a} M \to \ov{M}\to 0$, we get an exact sequence
\[
0\to \Hom_A(M,M)\xrightarrow{a}  \Hom_A(M,M) \to \Hom_A(M,\overline{M}) \to \Ext^1_A(M,M).
\]
Substituting our assumption $\Hom_A(M,M)\cong M^{\oplus \mathrm{rk}_M}$ into it yields an exact sequence 
\[
0\to \ov{M}^{\oplus \mathrm{rk}_M} \to \Hom _{\ov{A}}(\ov{M},\ov{M}) \to T\to 0
\]
of $\overline{A}$-modules,
where $T\subset \Ext^1_A(M,M)$ is a finitely presented $\overline{A}$-submodule (see \Cref{coherence of Tor and Ext}), which,
by the locally freeness of $\widetilde{M}$ over $\Spec A\backslash \{x\}$, is supported on $\{x\}$. Since $\dim(\overline{A}/\mathfrak{m}_V \overline{A})=d_x-1 \ge 2$, taking dual (as $\overline{A}$-modules) of the above short exact sequence and using \Cref{lem on isom of Ext's}, we see that
\[
(\ov{M}^{\vee})^{\oplus \mathrm{rk}_M} \simeq  \Hom _{\overline{A}}(\overline{M},\overline{M})^{\vee}.
\]
 Taking dual further and invoking \Cref{lem on dual and hom}, we get the following isomorphism
\[
(\ov{M}^{\ddual})^{\oplus \mathrm{rk}_M} \simeq  \Hom _{\ov{A}}(\ov{M}^{\ddual},\ov{M}^{\ddual}).
\]
Since the double dual $\ov{M}^{\ddual}$ is finitely presented and reflexive over $\ov{A}$  (\Cref{lem-rf} \ref{hom-rf}), we can apply the induction hypothesis to the $\ov{A}$-module $\ov{M}^{\ddual}$ and conclude that it is $\ov{A}$-free. The same lemma also implies that $\ov{M}^{\vee}$ is $\ov{A}$-reflexive, so that $\overline{M}^{\vee}\simeq \overline{M}^{\ddual\! \vee}$ is $\ov{A}$-free of rank $\mathrm{rk}_M$.

Finally, we show that $M$ is $A$-free. Since $\widetilde{M}$ is locally free over $\Spec A\backslash \{x\}$, the natural map $\overline{M} \to \overline{M}^{\ddual}$ is an isomorphism over $\Spec A\backslash \{x\}$. Since $\dim(\overline{A}/\mathfrak{m}_V\overline{A})=d_x-1>1$, we may apply \Cref{lem on isom of Ext's} to see that $\Ext_{\overline{A}}^1(\overline{M},\overline{A})   \simeq \Ext_{\overline{A}}^1(\overline{M}^{\ddual},\overline{A}) =0$. Since $a$ is $M$-regular, we have
\[
 \Ext_A^1(M,\overline{A}) \simeq
 \Ext_{\overline{A}}^1(M\otimes_A^{\textbf{L}}\overline{A},\overline{A})
 \simeq
 \Ext_{\overline{A}}^1(\overline{M},\overline{A})=0.
\]
 Applying Hom$_{A}(M,-)$ to the short exact sequence $0\to A\xrightarrow{a} A \to \overline{A}\to 0$, we get an exact sequence
 \[
   0 \to M^{\vee} \xrightarrow{a} M^{\vee} \to \Hom _A(M,\ov{A}) \to \Ext_A^1(M,A)
    \xrightarrow{a} \Ext_A^1(M,A) \to \Ext_A^1(M,\ov{A}).
\]
  All modules are finitely presented over $A$. Since $\Ext_A^1(M,\overline{A})=0$, Nakayama's lemma gives that $\Ext_A^1(M,A)=0$.
  {Therefore, $M^{\vee}/aM^{\vee}\simeq \Hom _A(M,\overline{A}) = \overline{M}^{\vee}$ is $\overline{A}$-free of rank $\mathrm{rk}_M$, implying that $\dim_{\kappa(x)}(M^{\vee}\otimes_A \kappa(x))=\mathrm{rk}_M=\mathrm{rk}_{M^{\vee}}$.} It follows that $M^{\vee}$, and equivalently $M\simeq M^{\ddual}$, is $A$-free.

\section{Generalities on torsors over algebraic spaces}

\bpp[Setup] \label{situation on S,X,G}
Throughout this section, we let $S$ denote a base scheme, $X$ an algebraic space over $S$, and $G$ an $X$-group algebraic space.
\epp
\begin{defn-tweak} \label{def of torsos}
\hfill
\begin{itemize}
\item [(1)]  A (right) \emph{$G$-torsor} (for the fppf topology) is an $X$-algebraic space $\cP$ equipped with a $G$-action $a:\cP\times_X G \to \cP$ such that the following conditions hold:
\benumr
\item the induced morphism $\cP\times_X G\isoto \cP \times_X \cP, (p,g)\mapsto (p,a(p,g))$, is an isomorphism; and
\item  there exists a fppf covering $\{X_i\to X\}_{i\in I}$ of algebraic spaces \SP{03Y8} such that $\cP(X_i)\neq \emptyset$ for every $i\in I$.
\eenum
\item [(2)] For $G$-torsors $\cP_1$ and $\cP_2$, a \emph{morphism} $\cP_1 \to \cP_2$ is a $G$-equivariant morphism $\cP_1\to \cP_2$ of $X$-algebraic spaces.
\item [(3)] By a trivialization of a $G$-torsor $\cP$ we mean a $G$-equivariant isomorphism $t:G \isoto \cP$, where $G$ acts on itself via right multiplication; this amounts to the choice of a section $t(1_G)\in \cP(X)$. A $G$-torsor $\cP$ is \emph{trivial} if there exists a trivialization, or, equivalently, if $\cP(X)\neq \emptyset$.
\end{itemize}
\end{defn-tweak}

Note that every morphism of two $G$-torsors is an isomorphism. To see this, one may pass to a fppf covering of $X$ to reduce to the case when both torsors are trivial; in this case the assertion is trivial.

\bremt
One can also define a sheaf torsor for an $X$-group algebraic space $G$. It is a sheaf $$\cP:(\Sch_{/S})_{\fppf}^{\mathrm{opp}}\to \Set$$ equipped with a map $\cP\to X$ of sheaves and a $G$-action $a:\cP\times_X G \to \cP$ such that the above two conditions (i) and (ii) in (1)  hold. However, it turns out that such a sheaf torsor is necessarily representable by an algebraic space, so working with sheaf torsors adds no more generality. To see this, let $\{X_i\to X\}_{i\in I}$ be a fppf covering as in (ii) that trivializes $\cP$. Then every $\cP\times_XX_i\simeq G\times_XX_i$ is an algebraic space, and the map
\[
  \tst \bigsqcup_i\cP\times_XX_i\to \cP
\]
is representable by algebraic spaces and is a fppf covering, because it is the base change of the fppf covering $\bigsqcup_i X_i\to X$ of algebraic spaces via $\cP\to X$. Here, all coproducts are taken in the category of sheaves on $(\Sch/S)_{\fppf}$. It follows from (3) of \SP{04S6} that $\cP$ is an algebraic space, as desired.
\eremt

 Let $\cP_1,\cP_2$ be two $G$-torsors. Define a functor
\[
\un{\mathrm{Isom}}_X(\cP_1,\cP_2):(\Sch_{/X})^{\mathrm{opp}}\to \Set
\]
which associates to any scheme $T$ over $X$ the set of $G_T$-equivariant isomorphisms $\cP_{1,T} \to \cP_{2,T}$ over $T$.
\blemt \label{lem on Isom}
For two $G$-torsors $\cP_1$ and $\cP_2$, $\un{\mathrm{Isom}}_X(\cP_1,\cP_2)$ is an algebraic space over $S$. Further,  $G\to X$ is quasi-compact (resp., \'etale, smooth, flat, separated, (locally) of finite type, (locally) of finite presentation, quasi-affine, affine, or finite) if and only if $\un{\mathrm{Isom}}_X(\cP_1,\cP_2)\to X$ is so.
\elemt
\bpf
Since $\un{\mathrm{Isom}}_{X}(\cP_1,\cP_2)$ is fppf locally on $X$ isomorphic to $G$, it admits a representable fppf covering by algebraic spaces, hence it is an algebraic space by \SP{04S6}.

The list properties of morphisms of algebraic spaces are all stable under base changes and are fppf local on the target, see \SP{03KG} (resp., \SPN{03XT}{03ZF}{03MM}{03KM}{040Y}{0410}{03WM}{03WG}{03ZQ}).
Consequently, since the functor $\un{\mathrm{Isom}}_{X}(\cP_1,\cP_2)$ is fppf locally on $X$ isomorphic to $G$, the properties of $G$ are inherited by and can be detected from $\un{\mathrm{Isom}}_{X}(\cP_1,\cP_2)$.
\epf

Since every $G$-torsor $\cP\to X$ trivializes over a fppf covering $\{X_i\to X\}$, one may try to obtain $\cP$ by gluing the trivial $G_{X_i}$-torsors $\cP_{X_i}$ using the canonical isomorphisms $$\phi_{ij}:(\cP_{X_i})_{X_{ij}} \simeq \cP_{X_{ij}} \simeq \cP_{X_j})_{X_{ji}}, \q \text{ where } \q X_{ij}=X_i\times_XX_j.
$$
It turns out that, unlike the case of schemes, this is always possible in the framework of algebraic spaces, see \Cref{descent lem for torsors}.
Note that, by taking $U\ce \bigsqcup_i X_i$, we may assume that $\cP_U$ is trivial for a fppf covering $U\to X$ with $U$ an algebraic space.

\bdt [Descent {data} for torsors] \label{def of descent datum for torsors}
Let $S$, $X$ and $G$ be as in \ref{situation on S,X,G}. Let $U \to X$ be a fppf covering of algebraic spaces over $S$. For every integer $n\ge 0$, denote by $U^{(n)}$ the $n$-fold fibre product of $U$ over $X$. The \emph{category of descent {data} } for $G$-torsors relative to $U\to X$, denoted $$\mathbf{Tors}\bigl((U^{(2)}\rightrightarrows U)_{\fppf}, G\bigr),$$ has pairs $(\cQ, \phi)$ as objects, where
\begin{itemize}
  \item $\cQ\to U$ is a $G_U$-torsor; and
  \item $\phi:\mathrm{pr}_1^*\cQ \isoto \mathrm{pr}_2^*\cQ$ is an isomorphism of $G_{U^{(2)}}$-torsors such that the following diagram commutes (i.e., the cocycle condition holds)
      \begin{equation*}
      \begin{tikzcd}
\mathrm{pr}_{12}^*\mathrm{pr}_{1}^*\cQ  \ar[r, "\mathrm{pr}_{12}^*(\phi)"] \ar[d, "\simeq"]
& \mathrm{pr}_{12}^*\mathrm{pr}_{2}^*\cQ
 \ar[r, "\simeq "]
& \mathrm{pr}_{23}^*\mathrm{pr}_{1}^*\cQ
\ar[r, "\mathrm{pr}_{23}^*(\phi)"]
& \mathrm{pr}_{23}^*\mathrm{pr}_{2}^*\cQ
 \ar[d,"\simeq"]
 \\
\mathrm{pr}_{13}^*\mathrm{pr}_{1}^*\cQ
\ar[rrr, "\mathrm{pr}_{13}^*(\phi)"]
&
&
& \mathrm{pr}_{13}^*\mathrm{pr}_{2}^*\cQ.
\end{tikzcd}
\end{equation*}

\end{itemize}
A \emph{morphism} from a pair $(\cQ, \phi)$ to another pair $(\cQ', \phi')$ is a morphism $\theta: \cQ\to \cQ'$ of $G_U$-torsors compatible with $\phi$ and $\phi'$, that is, $\mathrm{pr}_{2}^*(\theta)\phi=\phi' \mathrm{pr}_{1}^*(\theta)$.
\edt
To every $G$-torsor $\cP$ one can associate a pair $\Psi(\cP)\ce (\cP_U,\mathrm{can})$ via base changes, where $\mathrm{can}$ denotes the canonical isomorphism $\mathrm{pr}_{1}^*(\cP_U) \simeq \cP_{U^{(2)}} \simeq \mathrm{pr}_{2}^*(\cP_U)$. Thus we obtain a functor $$\Psi:\mathbf{Tors}(X_{\fppf},G)\to \mathbf{Tors}((U^{(2)}\rightrightarrows U)_{\fppf},G).$$

\blemt [Descent $G$-torsors]\label{descent lem for torsors}
$\Psi$ is an equivalence of {categories}.
\elemt
In other words, every descent data $(\cQ, \phi)$ for $G$-torsors are effective in the sense that there exists a $G$-torsor $\cP$ and an isomorphism $\cQ\simeq \cP_U$ compatible with $\theta$ and the canonical descent data for $\cP_U$.
\bpf
The full faithfulness of $\Psi$ follows from the sheaf property of the functor $\un{\mathrm{Isom}}_X(\cP_1,\cP_2)$ for any $G$-torsors $\cP_1$ and $\cP_2$. For the essential surjectivity, we pick a descent data $(\cQ, \phi)$, and need to show that there exists a $G$-torsor $\cP$ for which $(\cP_U,\mathrm{can})\simeq (\cQ, \phi)$.

 When both $X$ and $U$ are schemes, this is proven in \SP{04U1}. The case of algebraic spaces can be proved similarly, and we repeat the argument for convenience. First we view $\cQ$ as a sheaf on the site $(\mathbf{AS}/U)_{\fppf}$ (by the natural equivalence of the topoi associated to $(\mathbf{AS}/U)_{\fppf}$ and $(\Sch/U)_{\fppf}$). Since descent data for sheaves on any site are always effective \SP{04TR}, we may find a sheaf $\cP$ on the site $(\mathbf{AS}/X)_{\fppf}$ and an isomorphism of sheaves $\cP_U\simeq \cQ$ compatible with the descent data. Further, since maps of sheaves on any site can be glued \SP{04TQ}, the $G_U$-action on $\cQ$ descends to a $G$-action on $\cP$. All the assumptions (i) and (ii) of \Cref{def of torsos} hold, because they can be checked on the fppf covering $U\to X$. It remains to see that $\cP$ is representable by an algebraic space over $X$. However, this follows from (3) of \SP{04S6}, in view of the fact that the map $\cQ\to \cP$ is representable by algebraic spaces and is a fppf covering (being a base change of the fppf covering $U\to X$).
\epf

We end this section with the following result, which is used repeatedly in the sequel.
\blemt \label{An equivalence criterion}
Let $S$ be a scheme, $X$ an algebraic space over $S$, and $G$ an $X$-group algebraic space. Let $f\colon Y\ra X$ be a morphism of algebraic spaces over $S$.
Assume that the following conditions hold:
\benumr
\item {for every fppf covering $T\to X$ with $T$ a scheme, the map $G(T)\to G(Y_T)$ is bijective},
    where $Y_T\ce Y\times_X T$; and
\item for every $G_Y$-torsor $\cP$, there is an fppf covering $T\to X$ with $T$ a scheme such that $\cP_{Y_T}$ lies in the essential image of $f_T^*$, where $f_T\ce f\times_XT$.
\eenum
Then pullback induces an equivalence $f^*\colon \mathbf{Tors}(X_{\fppf}, {G} ) \isoto \mathbf{Tors}(Y_{\fppf},G_{Y})$.

Similarly, if $G\to X$ is smooth, then we have an equivalence $$f^*:\mathbf{Tors}(X_{\et}, {G} ) \isoto \mathbf{Tors}(Y_{\et},G_{Y}),$$ provided that one replaces `fppf' by `\'etale' everywhere in the above assumptions.
\elemt
\bpf
We prove the lemma for fppf torsors. The assumption (i) implies that the functor $f^*$ is fully faithful. It remains to check the essential surjectivity. Let $\cP$ be a $G_Y$-torsor. By assumption (ii) there is a fppf covering $T\to X$ with $T$ a scheme and a $G_T$-torsor $\cQ$ such that $f_{T}^*\cQ\simeq \cP_{Y_T}$. Using this isomorphism we can transform the canonical descent data on $\cP_{Y_T}$ to a descent data $$\theta:\mathrm{pr}_{1}^*f_{T}^*\cQ \isoto \mathrm{pr}_{2}^*f_{T}^*\cQ$$ on $f_{T}^*\cQ$ (relative to the covering $Y_T\to Y$). For every integer $n\ge 0$, denote by $T^{(n)}$ the $n$-fold fibre product of $T$ over $X$. Using the canonical identifications $$\mathrm{pr}_{1}^*f_{T}^*\cQ = f_{T^{(2)}}^*\mathrm{pr}_{1}^*\cQ \q \text{ and } \q \mathrm{pr}_{2}^*f_{T}^*\cQ = f_{T^{(2)}}^*\mathrm{pr}_{2}^*\cQ,$$ the full faithfulness of $f_{T^{(2)}}$ implies that there is a unique isomorphism $$
\tau:\mathrm{pr}_1^*\cQ \isoto \mathrm{pr}_2^*\cQ
$$
of $G_{T^{(2)}}$-torsors such that $f_{T^{(2)}}^*(\tau)=\theta$. Since
\[
\mathrm{pr}_{13}^*(\theta)=\mathrm{pr}_{13}^*(f_{T^{(2)}}^*(\tau))=f_{T^{(3)}}^*\mathrm{pr}_{13}^*(\tau)
\]
and
\begin{align*}
\mathrm{pr}_{13}^*(\theta) &= \mathrm{pr}_{23}^*(\theta)\mathrm{pr}_{12}^*(\theta) \\
&=\mathrm{pr}_{23}^*\left(f_{T^{(2)}}^*(\tau)\right)\mathrm{pr}_{12}^*\left(f_{T^{(2)}}^*(\tau)\right) \\
&=f_{T^{(3)}}^*\left(\mathrm{pr}_{23}^*(\tau)\right)f_{T^{(3)}}^*\left(\mathrm{pr}_{12}^*(\tau)\right) \\
&=f_{T^{(3)}}^*\left(\mathrm{pr}_{23}^*(\tau)\mathrm{pr}_{12}^*(\tau)\right),
\end{align*}
the full faithfulness of $f_{T^{(3)}}^*$ implies that $\mathrm{pr}_{13}^*(\tau)=\mathrm{pr}_{23}^*(\tau)\mathrm{pr}_{12}^*(\tau)$, that is, $\tau$ is a descent data on $\cQ$ relative to $T\to X$. By \Cref{descent lem for torsors}, there is a $G$-torsor $\cR$ and an isomorphism $(\cQ, \phi)\simeq (\cR_T, \mathrm{can})$ of descent datas. Pulling back to $Y_T$, we get an isomorphism of descent data $$(\cP_{Y_T}, \mathrm{can})\simeq f_T^*(\cQ, \tau)\simeq (\cR_{Y_T}, \mathrm{can}).$$ By \Cref{descent lem for torsors} again (applied to the covering $Y_T\to Y$), we see that $f^*(\cR)=\cR_Y\simeq \cP$, as desired.
\epf

\section{Purity for torsors and finite \'etale covers}
We begin with discussing generalities on linear groups that will be fundamental in multiple types of purities for torsors. The overall argument is bootstrapped from that for vector bundles, and controlling on the projective dimensions of extended reflexive sheaves leads to relative-dimensional constraints.
In particular, we obtain the purity for torsors on relative curves and its local variants in \S\ref{csub-purity-reductive}, where the constraints on dimensions in the local case are more flexible.
This allows us to shrink complements of domains of reductive torsors (\emph{cf}. \cite{GL23}*{Proposition~2.9 and Corollary~2.10}), which is crucial for our proof of the Grothendieck--Serre for constant reductive group schemes given in \cite{GL23}.
Finally, utilizing our Pr\"uferian analog of Auslander's flatness criterion (\emph{cf}. \Cref{Auslander's flatness criterion}), we establish in \S\ref{csub-purity-zg} the Purity \Cref{purity for finite flat group schemes} for torsors under finite locally free groups, and, consequently, we obtain a Pr\"uferian counterpart of the Zariski--Nagata Purity \Cref{fet cover}.

\bppt[Coaffine locally linear groups]\label{linear-coaffine}
Let $X$ be an algebraic space. An $X$-group algebraic space $G$ is \emph{linear} if there exists a group monomorphism $G\hra \GL(\sV)$ for a locally free $\sO_X$-module $\sV$ of finite rank.

An $X$-group algebraic space $G$ is \emph{fppf locally coaffine} (resp., \emph{\'etale locally coaffine}), if, fppf locally (resp., \'etale locally) on $X$, there is a monomorphism $G\hra \GL(\sV)$ such that the sheaf quotient $\GL(\sV)/G$ is representable by an $X$-affine algebraic space. Such a $G$ is necessarily affine over $X$ {and the monomorphism $G\hra \GL(\sV)$ is automatically a closed immersion}.
For instance, if $G$ is $X$-reductive (resp., $X$-finite locally free), then $G$ is \'etale locally coaffine\footnote{If $G$ is $X$-reductive, then \'etale locally $G$ splits,  so (by e.g., \cite{Gil21}*{Thm.~1.1 and Cor.~4.3}) there exists a closed immersion $G\hra \mathrm{GL}_{n,X}$ for some integer $n$, and \cite{Alp14}*{9.4.1} ensures that the quotient $\mathrm{GL}_{n,X}/G$ is $X$-affine of finite type.} (resp., fppf locally coaffine).
In the sequel, we will mainly consider (fppf or \'etale) locally coaffine $X$-group algebraic spaces $G$.

\eppt

\csub[Purity for torsors on relative curves]\label{csub-purity-reductive}

Now we study the extension behavior of torsors over relative curves. Motivated by \cite{EGAIV4}*{Proposition~21.9.4} that every invertible sheaf on a curve over a field extends across finitely many closed points, \Cref{ext-vb-curve} concerns relative curves over Pr\"ufer rings and generalizes \cite{Guo20b}*{Lemma~7.3}.
\bppt[Torsors on relative curves]\label{tors-curve}
Let $R$ be a semilocal Pr\"ufer domain with spectrum $S$ and $X$  an $S$-flat scheme of finite type with regular one-dimensional fibres.
Let $D\subset X$ be an $S$-quasi-finite closed subscheme contained in an affine open $\Spec A \subset X$, and {cut out by a nonzerodivisor $t\in A$.}
Assume that 
\benumr
\item  $\abs{S}$ is finite; or
\item  $D$ is $S$-finite.
\eenum
Note that these imply that $D$ is semilocal, and (ii) holds for instance when $X$ is $S$-proper. 
Consider
\begin{itemize}
\item $B_D\ce \Spec \wh{A}$, the formal neighborhood of $D$, where $\wh{A}\ce \varprojlim_{n} A/t^nA$;
\item $U_{D}\ce B_D\backslash D=\Spec \wh{A}[\frac{1}{t}]$, the punctured formal neighborhood.
\end{itemize}
Indeed, $B_D$ is semilocal, which follows from the semilocality of $D$.
As $D$ is $S$-quasi-finite, it is $S$-flat.
By \SP{0B9D}, $D$ is a relative effective Cartier divisor, so its each nonempty fibre has codimension one.
By Hensel's lemma, $(\wh{A},t\wh{A})$ is a Henselian pair, in particular, $t\wh{A}$ is contained in all maximal ideals of $\wh{A}$. 
Combining this with the fibrewise codimension-one property of $D$, we conclude that $B_D$ is semilocal.
\eppt

The following proposition specializes to \cite{Guo20b}*{Lem.~7.3} when $A=V[t]$ and $D=\Spec A/tA$ for a valuation ring $V$.
\bprop\label{ext-vb-curve}
With the setup \S\ref{tors-curve}, the restriction functor for vector bundles
\[
\text{$\mathbf{Vect}_X\ra \mathbf{Vect}_{X\backslash D}$\q is essentially surjective.}
\]
Furthermore,  we have
\[
H^1_{\Zar}(U_D,\GL_n)=H^1_{\et}(U_D,\GL_n)=\{\ast\}.
\]
\eprop
\bpf
Note that by \Cref{coh of V-flat ft alg}~\ref{G-R 17.4.1 bound} or \Cref{prufer-regularity}~\ref{prufer-regularity-ii}, all local rings of $X$ have weak dimension $\le 2$. By \Cref{ext-rflx}, every vector bundle on $X\backslash D$  extends to a reflexive sheaf on $X$, which, by \Cref{wdim2-vect}~\ref{wdim2-vect-ii}, is actually a vector bundle. This proves the first assertion.

Next, let $\sV$ be a vector bundle on $U_D=\Spec \wh{A}[\frac{1}{t}]$ and denote the Henselization of the pair $(X, D)$ by $(B_D^{\mathrm{h}},D)$. Write $B_D^{\mathrm{h}}=\Spec A^{\mathrm{h}}$ and set $U_D^{\mathrm{h}}\ce B_D^{\mathrm{h}}\backslash D=\Spec A^{\mathrm{h}}[\frac{1}{t}]$. Then, since the $t$-adic completion of $A^{\mathrm{h}}$ is isomorphic to $\wh{A}$,
\cite{BC22}*{Corollary~2.1.22~(c)}  implies that
  $\sV$ descends to a vector bundle $\sV^{\mathrm{h}}$ on $U_D^{\mathrm{h}}$.
Since $B_D^{\mathrm{h}}$ is the limit of elementary \'etale neighbourhoods $D\subset X'$ of $D\subset X$, by a limit argument, $\sV^{\mathrm{h}}$ descends to a vector bundle $\sV\pr$ on $X'\backslash D$ for some $D\subset X'$.
Since $\mathbf{Vect}_{X\pr}\ra \mathbf{Vect}_{X\pr\backslash D}$ is essentially surjective, $\sV\pr$ extends to a vector bundle $\wt{\sV}\pr$ on $X\pr$ whose pullback gives a vector bundle $\wt{\sV}$ on $B_D$ extending $\sV$.
Finally, as in \S\ref{tors-curve}, $\wh{A}$ is semilocal, so over which the bundle $\wt{\sV}$ and thus $\sV$ is trivial.
\epf

The following purity result in the case when $G$ is a reductive $X$-group scheme was proven in \cite{GL23}*{Theorem 2.7}. It turns out that the same argument works for any $X$-group algebraic space $G$ that is \'etale locally coaffine (\emph{cf}. \S\ref{linear-coaffine}) in the setting of algebraic spaces.
\bthm\label{purity for rel. dim 1}
Let $S$ be a semilocal affine Pr\"{u}fer scheme and $X$ an $S$-flat, locally of finite type algebraic space with regular one-dimensional $S$-fibres.
Let $G$ be an $X$-affine group algebraic space that is \'etale locally coaffine\footnote{For example, $G$ could be $X$-reductive, or $X$-finite locally free.}.
Given a closed subspace $Z$ of $X$ such that the inclusion $j\colon X\backslash Z\hra X$ is quasi-compact, and
\[
\text{$Z_{\eta}= \emptyset$\q for each generic point $\eta\in S$ \q and \q  $\codim(Z_s,X_s)\ge 1$ \q for all $s\in S$.}
\]
Then, the restriction gives an equivalence of categories of  $G$-torsors
\begin{equation} \label{restriction of G-torsors}
\mathbf{Tors}(X_{\fppf},G) \isoto \mathbf{Tors}((X\backslash Z)_{\fppf},G).
\end{equation}
Consequently, when considering isomorphism classes, there exists a bijection:
 $$
 H^1_{\fppf}(X,G)\simeq H^1_{\fppf}(X\backslash Z,G).
 $$
 Furthermore, if $S$ is allowed to be a general Pr\"ufer algebraic space\footnote{Namely, it admits an \'etale cover by a disjoint union of spectra of Pr\"ufer domains.} {(not necessarily semilocal)}, the above conclusions remain valid as long as $G$ is finitely presented over $X$ (and \'etale locally coaffine).
\ethm
\bpf
First, consider the case when $S$ is a semilocal affine Pr\"{u}fer scheme. We show that (\ref{restriction of G-torsors}) is an equivalence.
    Since $G$ is $X$-affine, by checking \'etale locally and using \Cref{Rflx-Pic-Vect on reg}, we see that $G(X) \simeq  G(X\backslash Z)$, which proves the full faithfulness (this uses the analogous bijection when we base change everything to a scheme \'etale over $X$ and the condition on fibres remains valid).

   For essential surjectivity, we choose a $G$-torsor $\cP$ over $X\backslash Z$ and wish to extend $\cP$ to a $G$-torsor over $X$. In the special case where $G=\GL_n$, $G$-torsors correspond to vector bundles of rank $n$.
{Due to the condition on fibres imposed on $Z$, by \Cref{prufer-regularity}~\ref{prufer-regularity-iii} and \Cref{Rflx-Pic-Vect on reg}, the pushforward $j_*\cP$ is reflexive.
Now, we leverage the essential assumption that $X$ has one-dimensional regular $S$-fibres: by \Cref{prufer-regularity}~\ref{prufer-regularity-ii}, every local ring of $X$ is coherent regular with $\wdim \leq 2$.
   Therefore, \Cref{wdim2-vect}~\ref{wdim2-vect-ii} implies that the reflexive $\sO_X$-module  $j_*\cP$ is locally free.}
   For a general $G$, by glueing in the \'etale topology, it suffices to demonstrate that $\cP$ extends, at least \emph{\'etale locally} on $X$, to a $G$-torsor over $X$ (see \Cref{An equivalence criterion}).

To prove this, we may assume that $X$ is affine, $G\subset \GL_{n,X}$, and that $\mathrm{GL}_{n,X}/G$ is affine over $X$. We exploit the following commutative diagram of pointed sets, where {the two rows are exact sequences}:
\[
\begin{tikzcd}
{(\mathrm{GL}_{n,X}/G)(X)} \arrow[r] \arrow[d,"\simeq" labl, ] & {H^1_{\fppf}(X,G)} \arrow[r] \arrow[d] & {H^1_{\fppf}(X,\mathrm{GL}_{n,X})} \arrow[d] \\
{(\mathrm{GL}_{n,X}/G)(X\backslash Z)} \arrow[r]           & {H^1_{\fppf}(X\backslash Z,G)} \arrow[r]           & {H^1_{\fppf}(X\backslash Z,\mathrm{GL}_{n,X})}.
\end{tikzcd}
\]
The bijectivity of the left vertical arrow follows from \Cref{Rflx-Pic-Vect on reg}.
By the case of vector bundles, we may replace $X$ by an affine open cover to assume that the induced $\GL_{n,X\backslash Z}$-torsor $\cP\times ^{G_{X\backslash Z}}\GL_{n,X\backslash Z}$ is already trivial.
Then, a diagram chase yields a $G$-torsor $\cQ$ over $X$ such that $\cQ|_{X\backslash Z}\simeq \cP$.

Now,  let us assume that $S$ is a Pr\"ufer algebraic space, and that $G$ is finitely presented over $X$. We may assume that $X$ is of finite type over $S$. By \Cref{coh of V-flat ft alg}~\ref{coherence of X}, the $S$-flatness of $X$ ensures its finite presentation over $S$. So, $G$ and hence all $G$-torsors are also finitely presented over $S$. At this point, we can use a standard argument to deduce the conclusion in the present case from the previous semilocal case. For example, to establish full faithfulness, we prove that $G(X) \simeq G(X\backslash Z)$.  The question is \'etale local on $X$, allowing us to assume that $S$ is an affine Pr\"ufer scheme and $X$ is an affine scheme. We consider both sides as Zariski sheaves on $S$, given by $T\mapsto G(X_T)$ and $T\mapsto G(X_T\backslash Z_T)$, respectively. The previous semilocal case implies that $G(X(s)) \simeq G(X(s)\backslash Z(s))$ for each point $s\in S$, where $X(s)\ce X\times_{S}\Spec(\sO_{S,s})$, etc. Since $j$ is quasi-compact and $G\to X$ is finitely presented, both sides of the last bijection can be identified with the stalks of the two Zariski sheaves at $s$. This implies $G(X)\simeq G(X\backslash Z)$, as desired.
\epf

\brem \label{rem on purity fo rel curves}
In higher relative dimensions, even in the classical Noetherian setting, the purity \Cref{purity for rel. dim 1} is inapplicable, even for the simplest group $G=\GL_{n}$. For instance, there exists a vector bundle over $\Spec(R)\backslash \{\mathfrak{m}_R\}$ that cannot be extended to $\Spec\, R$, where $(R,\mathfrak{m}_R)$ is any Noetherian regular local ring of Krull dimension at least three.
\erem

 The following is a local version of  \Cref{purity for rel. dim 1}.
 {Its proof closely mirrors that of \Cref{purity for rel. dim 1} (also refer to \cite{GL23}*{Theorem~2.8}, or can be straightforwardly deduced from it.}

\bthm[Local variant of purity] \label{extends across codim-2 points}
Let $V$ be a valuation ring of finite rank, with spectrum $S$, and {let $\eta\in S$ be the generic point}.
Let $X$ be an $S$-flat, finite type scheme with regular fibres. Let $G$ be an $X$-group scheme that is coaffine \'etale locally. If $x\in X$ satisfies one of the following
\benumr
\item\label{g-tor-local-gen} $x\in X_{\eta}$ with $\dim \sO_{X_{\eta},x} =2$, or
\item\label{g-tor-local-ngen} $x\in X_s$ (where $s \neq \eta$) with $\dim \sO_{X_s,x} =1$,
\eenum
then every $G$-torsor over $\Spec \sO_{X,x} \backslash \{x\}$ extends uniquely to a $G$-torsor over $ \sO_{X,x}$.
\ethm
It is worth noting that the stipulation of $V$ possessing finite rank guarantees that any finite type $V$-scheme will be topologically Noetherian. This ensures that the punctured spectrum $\Spec \sO_{X,x} \backslash \{x\}$ remains quasi-compact. This property is fundamental when establishing results concerning vector bundles.

As a corollary of \Cref{extends across codim-2 points}, one proves the following result about extending generically trivial torsors.

\bcor [Extending generically trivial torsors, \cite{GL23}*{Corollary 2.10}]\label{extend generically trivial torsors}
Fix
\benumr
\item $R$  a semilocal Pr\"{u}fer domain with spectrum $S$;
\item $X$ an $S$-flat finite type quasi-separated scheme with regular fibres;
\item $Y$ the spectrum of a local ring of an affine open subset of $X$;
\item $r\in R$ a nonzero element; and
\item $G$ a reductive $X$-group scheme.
\eenum
Then, every generically trivial $G$-torsor over $Y$ (resp., over an open subset of $X[\frac{1}{r}]\ce X_{R[1/r]}$) extends to a $G$-torsor over an open subset $U\subset X$. Here, the complementary closed $Z\ce X\backslash U$ satisfies the condition:
\[
\begin{cases}
\codim(Z_{\eta}, X_{\eta})\geq 3 & \text{for each generic point $\eta\in S$} \\
\codim(Z_s,X_s)\geq 2 & \text{for all $s\in S$.}
\end{cases}
\]
\ecor

\csub[Purity for finite locally free torsors and the Zariski--Nagata]\label{csub-purity-zg}
By combining \Cref{Rflx-Pic-Vect on reg} on the purity of reflexive sheaves and Auslander's flatness criterion \Cref{Auslander's flatness criterion}, we are able to establish the following Pr\"uferian analog of a result of Moret-Bailly, detailed in \cite{Mar16}.
\bthm [Purity for torsors under finite locally free groups] \label{purity for finite flat group schemes}
\hfill
\benumr
\item \label{purity-G-gl} Let $S$ be a {Pr\"ufer algebraic space} and $X$ an $S$-smooth algebraic space.
Let $G$ be an $X$-finite, locally free group algebraic space.
Given a closed $Z\subset X$ such that $j\colon X\backslash Z\hra X$ is quasi-compact, and
\[
   \begin{cases}
   \codim(Z_{\eta}, X_{\eta})\geq 2 & \text{for each generic point $\eta\in S$}\\
   \codim(Z_s,X_s)\geq 1 & \text{for all $s\in S$,}
   \end{cases}
\]
then the restriction induces an equivalence of categories of $G$-torsors:
\[
\mathbf{Tors}(X_{\fppf},G) \isoto \mathbf{Tors}((X\backslash Z)_{\fppf},G).
\]
Consequently, when considering isomorphism classes, there exists a bijection:
\[
H^1_{\fppf}(X,G)\simeq H^1_{\fppf}(X\backslash Z,G).
\]
\item \label{purity-G-loc} Let $V$ be a finite-rank valuation ring, $X$ a $V$-smooth scheme, and $G$ an $X$-finite locally free group scheme.
Let $x\in X$ be a point such that $\dim \sO_{X,x} \geq 2$.
If $x$ is not a maximal point in the $S$-fibres of $X$, then the restriction functor establishes an equivalence of {categories} of $G$-torsors:
\[
\mathbf{Tors}((\Spec \sO_{X,x})_{\fppf},G) \isoto \mathbf{Tors}((\Spec \sO_{X,x}\backslash \{x\})_{\fppf},G).
\]
Consequently, when considering isomorphism classes, there exists a {bijection}:
 $$
 H^1_{\fppf}(\Spec \sO_{X,x},G)\simeq H^1_{\fppf}(\Spec \sO_{X,x}\backslash \{x\},G).
 $$
\eenum
\ethm
We anticipate that the theorem is valid for all $X$ which are flat, locally finitely presented over $S$ (or over $V$), with regular fibres. However, our approach does not confirm this due to its reliance on Auslander's criterion for flatness (see \Cref{Auslander's flatness criterion}), which necessitates the smoothness of $X$ over $S$ (or over $V$).

\bpf
(i) We may assume that $X\to S$ is quasi-compact. The assumption implies that $X$ and $G$ are both finitely presented over $S$. Thus, similar to the proof of \Cref{purity for rel. dim 1}, we are reduced to the case where $S$ is the spectrum of a valuation ring.

It remains to verify the assumptions of \Cref{An equivalence criterion}. All assumptions of \Cref{An equivalence criterion} are \'etale local on $X$, so we may assume that $X$ is a scheme, finitely presented over $S$.
Employing the limit argument \Cref{lim-codim} involving \Cref{approxm semi-local Prufer ring}, we can further restrict our scenario to instances where $S$ has a finite Krull dimension. Furthermore, since $|S|$ is finite and each $R$-fibre of $X$ is Noetherian, $|X|$ is Noetherian.

 The condition (i) of \Cref{An equivalence criterion} can be deduced from \Cref{prufer-regularity}~\ref{prufer-regularity-iii} and \Cref{Rflx-Pic-Vect on reg}, which needs the condition on fibres of $Z$.

To verify the condition (ii) of \Cref{An equivalence criterion}, we will check that, \'etale locally on $X$, every $G$-torsor over $X\backslash Z$ extends to a $G$-torsor over $X$.
Let $\cP$ be a $G_{X\backslash Z}$-torsor.
By \Cref{prufer-regularity}~\ref{prufer-regularity-iii} and \Cref{pp-ref}, $j_{\ast}\sO_{\cP}$ is a reflexive $\sO_X$-module.
First, we prove the $\sO_X$-flatness of $j_{\ast}\sO_{\cP}$.
 We can use Noetherian induction to reduce to the case where $X$ is local, essentially smooth over $R$, and $Z=\{x\}$ is its closed point.
Then, Auslander's criterion \Cref{Auslander's flatness criterion} reduces us to showing  the isomorphism
\[
\sHom_{\sO_X}(j_{\ast}\sO_{\cP},j_{\ast}\sO_{\cP}) \simeq  \left(j_{\ast}\sO_{\cP}\right)^{\oplus r}, \q \x{where}\q r=\mathrm{rank} _{\sO_X}\sO_G.
\]
Note that in such local case, we have $\sO_G\simeq \sO_{X}^{\oplus r}$. Consider the following map
\begin{align*}
\sHom_{\sO_X}(\sO_G,j_{\ast}\sO_{\cP}) & \to \sHom_{\sO_X}(j_{\ast}\sO_{\cP},j_{\ast}\sO_{\cP}), \\
f & \mapsto \Bigl( j_{\ast}\sO_{\cP} \xrightarrow{j_{\ast}\rho} \sO_G\otimes_{\sO_X}j_{\ast}\sO_{\cP} \xrightarrow{(f,\text{id})} j_{\ast}\sO_{\cP}   \Bigr)
\end{align*}
of reflexive $\sO_X$-modules. This is an isomorphism: by \Cref{Rflx-Pic-Vect on reg}, it suffices to argue over $X\backslash Z$, then its explicit inverse is
\[
g \mapsto  \Bigl( \sO_{G_{X\backslash Z}} \xrightarrow{\text{id}\otimes 1}\sO_{G_{X\backslash Z}}\otimes_{\sO_{X\backslash Z}}\sO_{\cP}   \xrightarrow{(\rho, \text{id})^{-1}} \sO_{\cP}\otimes_{\sO_{X\backslash Z}}\sO_{\cP} \xrightarrow{(g, \text{id})}\sO_{\cP} \Bigr).
\]

We now prove that the $G$-torsor structure of $\cP$ extends uniquely to that of the scheme $\un{\Spec}_X(j_{\ast}\sO_{\cP})$.
As $G$ is finite locally free, by projection formula \SP{01E8}, taking $j_{\ast}$ of the co-action $\rho:\sO_{\cP}\to j^*\sO_{G}\otimes_{\sO_{X\backslash Z}}\sO_{\cP}$ yields
\[
j_{\ast}\rho\colon  \q j_{\ast}\sO_{\cP}\to \sO_G\otimes_{\sO_X} j_{\ast}\sO_{\cP}.
\]
To check that $j_{\ast}\rho$ is a co-action, we verify the associativity, the commutativity of the following diagram
\begin{equation*}
\begin{tikzcd}
j_{\ast}\sO_{\cP} \arrow[r, "j_{\ast}(\rho)"]  \arrow[d,"j_{\ast}(\rho)"]
 & \sO_G\otimes_{\sO_X}j_{\ast}\sO_{\cP} \arrow[d,"\text{id}\otimes j_{\ast}(\rho)"] \\
\sO_G\otimes_{\sO_X}j_{\ast}\sO_{\cP} \arrow[r, "\mu_G \otimes \text{id}"]
& \sO_G\otimes_{\sO_X}\sO_G\otimes_{\sO_X}j_{\ast}\sO_{\cP},
\end{tikzcd}
\end{equation*}
where $\mu_G:\sO_G\to \sO_G\otimes_{\sO_X}\sO_G$ is the co-multiplication of $G$.
Since all sheaves involved are $\sO_X$-reflexive, the commutativity over $X\backslash Z$ by \Cref{Rflx-Pic-Vect on reg} extends over $X$.
Finally, the following map
\[
(j_{\ast}\rho, 1\otimes \text{id})\colon\q  j_{\ast}\sO_{\cP}\otimes_{\sO_X}j_{\ast}\sO_{\cP} \to \sO_G\otimes_{\sO_X}j_{\ast}\sO_{\cP},
\]
by the $\sO_X$-flatness of $j_{\ast}\sO_{\cP}$ and \Cref{Rflx-Pic-Vect on reg}, is an isomorphism since so is its restriction on $X\backslash Z$.

(ii) A similar argument can be employed for proof. Specifically, to establish the essential surjectivity of the restriction functor, the finite-rank assumption on $V$ ensures that $j:\Spec \sO_{X,x}\backslash \{x\} \hookrightarrow \Spec \sO_{X,x}$ is quasi-compact quasi-separated. Consequently, for any $G$-torsor $\cP$ over $\Spec \sO_{X,x}\backslash \{x\}$, $j_{\ast}\sO_{\cP}$ is a reflexive $\sO_{X,x}$-module (by \Cref{Rflx-Pic-Vect on reg}). By invoking Auslander's criterion (see \Cref{Auslander's flatness criterion}), we deduce that $j_{\ast}\sO_{\cP}$ is free over $\sO_{X,x}$ and retains the $G$-torsor structure from $\cP$.
 This subsequently yields the sought-after extension of $\cP$ to $\Spec \sO_{X,x}$.
\qedhere
\epf

From \Cref{purity for finite flat group schemes} one can now derive an analog of the classical Zariski--Nagata.
\bthm [Zariski--Nagata purity for finite \'etale covers] \label{fet cover}
\hfill
\benumr
\item \label{fet-gl version} Let $S$ be a Pr\"ufer algebraic space, and $X$ an $S$-smooth algebraic space. Given a closed subspace $Z$ of $X$ such that $j\colon X\backslash Z\hra X$ is quasi-compact, and that
\[
   \begin{cases}
    \codim(Z_{\eta}, X_{\eta})\geq 2 & \text{for each generic point $\eta\in S$}\\
    \codim(Z_s,X_s)\geq 1 & \text{for all $s\in S$},
   \end{cases}
\]
then the restriction gives an equivalence of categories of finite \'etale objects
\[
\mathrm{F\'Et}_X\isoto \mathrm{F\'Et}_{X\backslash Z}.
\]

\item \label{fet-loc version} Let $R$ be a finite-rank valuation ring with spectrum $S$, and consider an $S$-smooth scheme $X$. Given a point $x \in X$ with properties such that $\dim \sO_{X,x} \geq 2$ and $x$ is not a maximal point in the $S$-fibres of $X$. Then the restriction functor establishes an equivalence of categories of finite \'etale covers
\[
\mathrm{F\'Et}_{\Spec \sO_{X,x}}\isoto \mathrm{F\'Et}_{\Spec \sO_{X,x}\backslash \{x\}}.
\]
\eenum
\ethm
Again, we anticipate that the theorem is valid for all $X$ which are flat, locally finitely presented over $S$, with regular fibres.

\bpf
(i) Full faithfulness. For finite \'etale covers $\pi_i:X_i\to X$, $i=1,2$, consider 
\[
   \text{the $X$-functor\q $Y\ce \un{\Hom}_X(X_1,X_2)$}
\]
 that parameterizes $X$-morphisms from $X_1$ to $X_2$; it is a subfunctor of $$\un{\Hom}_X(\pi_{2,*}\sO_{X_2},\pi_{1,*}\sO_{X_1})$$ consisting of sections compatible with algebraic structures of $\pi_{2,*}\sO_{X_2}$ and $\pi_{1,*}\sO_{X_1}$, which amount to the commutativity of a certain diagram of $\sO_X$-modules. It follows that $Y\subset \un{\Hom}_X(\pi_{2,*}\sO_{X_2},\pi_{1,*}\sO_{X_1})$ is a closed subfunctor.
Thus, $Y$ is an algebraic space finite over $X$. (Using the infinitesimal criterion for formal smoothness, one can check that $Y\to X$ is even finite \'etale, but we will not need this in the sequel.)  By \Cref{Rflx-Pic-Vect on reg}, we have $Y(X)\simeq Y(X\backslash Z)$, thereby proving the full faithfulness.

Essential surjectivity.
Given the full faithfulness established in the preceding paragraph, we can first pass to an \'etale cover and then to a connected component, allowing us to assume that $X$ is an integral affine scheme. To apply \Cref{purity for finite flat group schemes}~\ref{purity-G-gl} and conclude, {it suffices to observe that the category of finite \'etale covers of $X$ of degree $n$ (with isomorphisms) is equivalent to the category of $\fS_{n,X}$-torsors. (Here $\fS_n$ denotes the $n$-th symmetric group.)}

(ii) This is proved in the same manner as (i), using \Cref{purity for finite flat group schemes}~\ref{purity-G-loc} instead of \Cref{purity for finite flat group schemes}~\ref{purity-G-gl}.
\epf

\section{Cohomology of groups of multiplicative type}
Inspired by the purity results in \cite{CS21}*{Theorem~7.2.8}, we investigate the parafactoriality over Pr\"ufer bases and then present the purity for cohomology of group algebraic spaces of multiplicative type.

\csub[Geometrically parafactorial pairs]
\vskip -0.5cm

\bppt[Parafactorial pairs]\label{parafact}
Let $(X,\sO_X)$ be a ringed space with a closed subspace $Z\subset X$ and the canonical open immersion $j\colon V\ce X\backslash Z\hra X$, if for every open subspace $U\subset X$ the following restriction is an equivalence of categories
\[
  \text{$\mathbf{Pic}\,U\isoto \mathbf{Pic}\,(U\cap V)$\qq $\sL\mapsto \sL|_{U\cap V}$,}
\]
then the pair $(X,Z)$ is \emph{parafactorial}.
In particular, for an invertible $\sO_X$-module $\sL$,
\[
\text{$\sL(U)=\Hom_{\sO_U}(\sO_U,\sL|_U)\simeq \Hom_{\sO_{U\cap V}}(\sO_{U\cap V},\sL|_{U\cap V})=\sL(U\cap V)$}
\]
for all open $U\subset X$;
in other words, $\sL\simeq j_{\ast}j^{\ast}\sL$.
A local ring $A$ is \emph{parafactorial} if the pair
 $(\Spec A, \{\fm_A\})$ is parafactorial.
We list several parafactorial pairs $(X,Z)$ and local rings.
\benumr
\item \label{para-i}By \cite{EGAIV4}*{Proposition~21.13.8}, a local ring $A$ is parafactorial amounts~to
\[
  \text{$\Pic \p{\Spec (A)\backslash \{\fm_A\}}=0$\q and \q $A\simeq \GG(\Spec (A)\backslash \{\fm_A\},\wt{A})$;}
\]
\item  By \cite{EGAIV4}*{Exemples~21.13.9~(ii)}, a Noetherian factorial local ring is parafactorial if and only if its Krull dimension is at least 2;
\item When $X$ is locally Noetherian, $Z$ satisfies $\codim(Z,X)\geq 4$, and  $\sO_{X,z}$ are locally complete intersection\footnote{This means that its completion is the quotient of complete regular local ring by an ideal generated by a regular sequence.} for all $z\in Z$, by \cite{SGA2new}*{Exposé~XI, Proposition~3.3 and Théorème~3.13~(ii)}, the pair $(X,Z)$ is parafactorial;
\item\label{para-iv}  \label{parafactorial pair} For a normal scheme $S$, an $S$-smooth scheme $X$ and a closed $Z\subset X$ satisfying
\[
   \begin{cases}
   \codim(Z_{\eta},X_{\eta})\geq 2 & \text{for each generic point $\eta\in S$}\\
   \codim(Z_s,X_s)\geq 1 & \text{for every $s\in S$},
   \end{cases}
\]
by \cite{EGAIV4}*{Proposition~21.14.3}, the pair $(X,Z)$ is parafactorial.
\eenum
\eppt

\bppt[Geometrically parafactorial pairs]\label{geom-parafact}
Below, we are mainly interested in the case when $X$ is a scheme, and $Z\subset X$ is a closed subscheme such that the canonical open immersion $j\colon V\ce X\backslash Z\hra X$ is quasi-compact. By \cite{EGAIV2}*{Lemme~2.3.1}, the quasi-compactness of $j$ guarantees that the pushforward by $j$ of a quasi-coherent $\sO_{V}$-module is quasi-coherent and its formation commutes with arbitrary flat base changes (in particular, localizations).

A pair $(X,Z)$ is \emph{geometrically parafactorial} if $j:V=X\backslash Z\hra X$ is quasi-compact and if, for every $X$-\'etale $X\pr$ with base change $Z\pr\ce Z\times_X X\pr$, the pair $(X\pr, Z\pr)$ is parafactorial.
A local ring $A$ with a quasi-compact punctured spectrum is \emph{geometrically parafactorial} if its strict Henselization $A^{\sh}$ is parafactorial. By the following \Cref{geom paraf vs paraf}, one can see that a local ring $A$ with a quasi-compact punctured spectrum is geometrically parafactorial if and only if the pair $(\Spec A, \{\fm_A\})$ is geometrically parafactorial.
\eppt

\blem \label{geom paraf vs paraf}
Let $A$ be a local ring with a quasi-compact punctured spectrum. Then $A$ is geometrically parafactorial if and only if for any {local and essentially \'etale\footnote{Recall that by definition an essentially \'etale ring map is a localization of an \'etale ring map.}} map $A\to B$ of local rings, $B$ is parafactorial.
\elem
\bpf
Assume that $A$ is geometrically parafactorial, that is, $A^{\sh}$ is parafactorial. Let $A\to B$ be a local and {essentially} \'etale map; we will show that $B$ is parafactorial. For any local ring $C$, denote by $j_C:U_C^{\circ}\ce \Spec (C)\backslash \{\fm_C\}\hra U_C\ce \Spec \,(C)$ the canonical open immersion. Choose an $A$-map $B\to A^{\sh}$. Let $\sL$ be an invertible $\sO_{U_B^{\circ}}$-module. By the quasi-compactness of $j_B$ (inherited from that of $j=j_A$) and the faithful flatness of $U_{A^{\sh}}\to U_B$, the $\sO_{U_B}$-module $j_{B,*}\sL$ is quasi-coherent and its pullback to $U_{A^{\sh}}$ is isomorphic to $j_{A^{\sh},*}(\sL|_{U_{A^{\sh}}^{\circ}})$. Since $A^{\sh}$ is parafactorial, we have $j_{A^{\sh},*}(\sL|_{U_{A^{\sh}}^{\circ}})\simeq \sO_{U_{A^{\sh}}}$. Descent theory implies that $j_{B,*}\sL\simeq \sO_{U_{B}}$ and thus $\sL\simeq \sO_{U_B^{\circ}}$. Similarly, by considering the pullback to $U_{A^{\sh}}$, we see that the natural map $\sO_{U_B}\to j_{B,*}\sO_{U_B^{\circ}}$ is bijective, that is, $B\simeq \Gamma(U_B^{\circ},\sO_{U_B^{\circ}})$. This proves that $B$ is parafactorial (\emph{cf}. \S\ref{parafact}~\ref{para-i}).

For the other side, fixing a separable closure $\overline{\kappa}_A$ of $\kappa_A=A/\fm_A$ and a geometric point $\overline{t}:A\to \overline{\kappa}_A$, then $A^{\sh}$ is the filtered colimit of all $B$ for {essentially \'etale local ring maps $A\to B$ along with a geometric point $\overline{t}_B:B\to \overline{\kappa}_A$ lifting $\overline{t}$}. Consequently, if all such $B$ are parafactorial, then we have the following equivalences
\[
\mathbf{Pic}\,U_{A^{\sh}}\xleftarrow{\sim} 2\text{-colim}_{(B,\overline{t}_B)}\mathbf{Pic}\,U_B \xrightarrow{\sim}
2\text{-colim}_{(B,\overline{t}_B)}\mathbf{Pic}\,U_B^{\circ}
\xrightarrow{\sim} \mathbf{Pic}\,U_{A^{\sh}}^{\circ},
\]
where the rightmost equivalence follows from \SP{0B8W}, because $U_A^{\circ}$ is quasi-compact. This proves that $A$ is geometrically parafactorial.
\epf

The following is a generalization of \cite{EGAIV4}*{Proposition~21.13.10} to the case of topologically locally Noetherian schemes\footnote{A scheme is topologically locally Noetherian if it admits a cover by open subschemes whose underlying topological spaces are Noetherian (i.e., any descending sequence of closed subsets is eventually constant); it is topologically Noetherian if its underlying topological space is Noetherian.}.
\blem\label{lim-parafactorial}
For a topologically locally Noetherian scheme $X$ and a closed subscheme $Z\subset X$,
\benumr
\item the pair $(X,Z)$ is parafactorial if and only if $\sO_{X,z}$ is parafactorial for every $z\in Z$;
\item the pair $(X,Z)$ is geometrically parafactorial if and only if $\sO_{X,\overline{z}}^{\sh}$ is parafactorial for every geometric point $\overline{z}\to Z$.
\eenum
\elem
\bpf
By \cite{EGAIV4}*{Corollaire~21.13.6~(i)}, we can work Zariski locally on $X$, so we may assume throughout that $X$ is topologically Noetherian. By \cite{BS15}*{Lemma~6.6.10~(3)}, any quasi-compact \'etale cover of $X$ is also topologically Noetherian. Therefore, taking into account of \Cref{geom paraf vs paraf}, we see that (ii) follows from (i) because $\Spec \sO_{X,\overline{z}}^{\sh}$ is the inverse limit of \'etale neighborhoods of $\overline{z}\to X$.

It remains to prove (i).
Assume that $(X,Z)$ is a parafactorial pair and denote by $j\colon V\ce X\backslash Z\hra X$ the canonical open immersion. For each $z\in Z$, denote $U_z\ce \Spec \sO_{X,z}$ and $U_z^{\circ}\ce U_z\backslash \{z\}$.
To show that $\sO_{X,z}$ is parafactorial, we prove that every invertible $\sO_{U_z^{\circ}}$-module $\sL_{0}$ is isomorphic to $\sO_{U_z^{\circ}}$.
By \cite{EGAIV3}*{Proposition~8.2.13} and \cite{EGAI}*{Proposition~2.4.2}, $U_z^{\circ}$ is the inverse limit of $B^{\circ}\ce B\backslash (B\cap \ov{\{z\}})$, where $B$ ranges over all open affine neighborhoods of $z\in X$.
Since $B^{\circ}\hookrightarrow B$ is quasi-compact (by topological Noetherianness), the limit argument \SP{0B8W} implies that there exists such a $B$ and an invertible $\sO_{B^{\circ}}$-module $\sL_{B^{\circ}}$ for which $\sL_0\simeq \sL_{B^{\circ}}|_{U_z^{\circ}}$.
By assumption and \cite{EGAIV4}*{Corollaire~21.13.6~(i)(ii)}, the pair $(B,B\cap \ov{\{z\}})$ is parafactorial.
Thus, there exists an invertible $\sO_B$-module $\wt{\sL}_B$ such that $\wt{\sL}_B|_{B^{\circ}}\simeq \sL_{B^{\circ}}$. Shrinking $B$ if necessary, we have $\wt{\sL}_B\simeq \sO_B$ hence $\sL_0\simeq \sO_{U_z^{\circ}}$.

Conversely, assuming that $\sO_{X,z}$ are parafactorial for all $z\in Z$, we will prove that the pair $(X,Z)$ is parafactorial.
We first show that for any open $U\subset X$ with the open immersion $j_U:U\cap V\hra U$, the canonical map $\sE\to j_{U,*}(\sE|_{U\cap V})$ is bijective for any finite-rank locally free $\sO_U$-module $\sE$. By \cite{EGAIV4}*{Corollaire~21.13.3~(b)}, this is equivalent to saying that $\sO_{X}\xrightarrow{\sim}j_*\sO_{V}$ is bijective.
The Noetherian assumption implies that $j$
is quasi-compact, so $j_*\sO_V$ is quasi-coherent and commutes with flat base changes. It suffices to prove that we have an isomorphism after localizing at points
$z\in X$. This is trivial if $z\in V$, and if
$z\in Z$ it follows from parafactoriality of $\sO_{X,z}$ (see \S\ref{parafact}~\ref{para-i}).

 To finish the proof, by \cite{EGAIV4}*{Proposition~21.13.5}, it remains to show that for every invertible $\sO_{V}$-module $\sL$, the pushforward $j_{\ast}\sL$ is an invertible $\sO_{X}$-module. 
 {For this, we consider the subset}
\[
   \Omega\ce \{x\in X\,|\, \x{$j_{\ast}\sL$ is invertible on an open neighborhood of $x$}\};
\]
it is open in $X$ and contains $V=X\backslash Z$.
Denote $  Y\ce X\backslash \Omega\subset Z$.
If $Y\neq \emptyset$, we choose a maximal point $y\in Y$ (every non-empty scheme has maximal points).
Then $\Omega\cap U_y=U_y^{\circ}$, and so $\sL_0\ce (j_{\ast}\sL)|_{U_y^{\circ}}$ is an invertible $\sO_{U_y^{\circ}}$-module.
The parafactoriality of $\sO_{X,y}$ yields an extension of $\sL_0$ to an invertible $\sO_{U_y}$-module $\wt{\sL}_0$, which, by the limit argument \SP{0B8W} again, descends to an invertible $\sO_{W}$-module $\wt{\sL}_W$ for an open affine neighborhood $W$ of $y\in X$. As $W$ shrinks, it becomes $U_y$ and $\Omega \cap W$ becomes $U_y^{\circ}$. Since $\Omega\cap W$ is quasi-compact, \emph{loc. cit.} implies that we may shrink $W$ to assume that the restrictions of $j_{\ast}\sL$ and $\wt{\sL}_W$ to $\Omega \cap W$ are equal. Set $\Omega\pr\ce \Omega\cup W$. By Zariski glueing, we obtain an invertible $\sO_{\Omega\pr}$-module $\sL\pr$ such that $\sL\pr|_W=\wt{\sL}_W$ and $\sL\pr|_{\Omega}=(j_{\ast}\sL)|_{\Omega}$.
Since $V\subset \Omega$, we have $\sL\pr|_{V}=\sL$, and so $\sL\pr \xrightarrow{\sim} (j_{\ast}(\sL\pr|_{V}))|_{\Omega\pr}=(j_{\ast}\sL)|_{\Omega\pr} $ (by the previous paragraph), which leads to a desired contradiction with the definition of $\Omega$.
\epf

\bprop \label{geom parafact}
Let $S$ be a normal scheme and $X$ an $S$-scheme. Assume that one of the following holds:
\benumr
\item either $X\ra S$ is a smooth morphism of topologically locally Noetherian schemes; or
\item $S$ is semilocal Pr\"ufer of finite dimension and $X$ is $S$-flat locally of finite type with regular $S$-fibres.
\eenum
Then, if $x\in X$ is not a maximal point of $S$-fibres of $X$ and satisfies $\dim \sO_{X,x}\geq 2$, it holds that
\[
   \text{$\sO_{X,x}$ is geometrically parafactorial,\q  namely, \q $\sO^{\sh}_{X,x}$ is parafactorial. }
\]
\eprop
\bpf
Notice that, in both cases (i)-(ii), the scheme $X$ is topologically locally Noetherian. By \Cref{geom paraf vs paraf}, the parafactoriality of $\sO_{X,x}^{\sh}$ is equivalent to those of $ \sO_{X\pr,x\pr}$ for all $X$-\'etale $X\pr$ and $x\pr\in X\pr$ lying over $x$. Moreover, since all $X\pr$ and $x\pr$ satisfy the conditions in the statement above (in the case (i), the topologically locally Noetherianness of $X\pr$ follows from \cite{BS15}*{Lemma~6.6.10~(3)}), thus it suffices to show that $\sO_{X,x}$ is parafactorial.
For the Zariski closure $Z\ce \ov{\{x\}}$, by \Cref{lim-parafactorial}, we are reduced to finding a small open neighborhood $U$ of $x\in X$ such that $(U, Z\cap U)$ is a parafactorial pair.
Now, take an arbitrary open neighborhood $U$ of $x\in X$, by \cite{EGAIV3}*{Proposition~9.5.3} applied to $Z\subset X$, shrinking $U$ if needed, we may assume that $U\cap Z$ does not contain any irreducible components of $S$-fibres of $X$.
In addition, if some $z\in Z$ lies over a maximal point $\eta\in S$, since $x$ specializes to $z$, then we have $\dim \sO_{X_{\eta},z}= \dim \sO_{X,z}\geq 2$. Consequently, we have $\codim(X_{\eta}\cap Z, X_{\eta})\geq 2$ and, by \S\ref{parafact}~\ref{para-iv} (in case (i)) and \Cref{Rflx-Pic-Vect on reg} (in case (ii)), the desired parafactoriality of $(U,Z\cap U)$ follows.
\epf

\csub[Purity for groups of multiplicative type]\label{subsection-purity of mult}
   The goal of this subsection is to study purity for groups of multiplicative type over algebraic spaces that are \emph{topologically locally Noetherian} in the following sense: they admit  \'etale covers by topologically locally Noetherian schemes. By \cite{BS15}*{Lemma~6.6.10~(3)}, we see that any scheme \'etale over a topologically locally Noetherian algebraic space is again topologically locally Noetherian.
   
   \bppt[The fppf local cohomology]\label{par-The fppf local cohomology}
Let $X$ be an algebraic space. Let $(\mathbf{Sch}_{/X})_{\fppf}$ denote the site of schemes over $X$ with fppf covers. {Recall that the sections of a sheaf on $(\mathbf{Sch}_{/X})_{\fppf}$ over an $X$-algebraic space $U$ are defined as the set of morphisms from the fppf sheaf
\[
\mathbf{Sch}_{/X} \ni V  \mapsto \Hom_{X}(V, U).
\]}
Given a closed subspace $Z\subset X$, denote by $j: X\backslash Z \hookrightarrow X$ the open immersion. For an abelian sheaf $\sF$ on $(\mathbf{Sch}_{/X})_{\fppf}$, we define
\[
\cH_Z^0(\sF)(V)\ce \text{Ker}(\sF(V)\to \sF(V\backslash (V\times_XZ))) ,\qq V\in \mathbf{Sch}_{/X};
\]
it is the largest subsheaf of $\sF$ supported on $Z$.
For an algebraic space $U$ over $X$, we denote
$$
H_Z^0(U,\sF)\ce \GG(U, \cH_Z^0(\sF)).
$$
Both the functors $\sF\mapsto \cH_Z^0(\sF)$ and $\sF \mapsto H_Z^0(U,\sF)$ are left exact. Let $\cH_Z^i$ and $H_Z^i(U,-)$ denote their $i$-th right derived functors. From the definition we see that $\cH_Z^i(\sF)$ can be identified with the sheafification of the presheaf $V\mapsto H^i_Z(V,\sF)$. Moreover, as the functor $\cH_Z^0$ sends injective sheaves to injective sheaves (because it admits an exact left adjoint), we have the local-to-global $E_2$-spectral sequence:
\[
E_2^{pq}=H^p(U,\cH^q_Z(\sF))\Rightarrow H^{p+q}_Z(U,\sF).
\]
{By employing an injective resolution, one also derives the long exact sequence for fppf local cohomology, analogous to the approach used for Zariski local cohomology as discussed in \S\S \ref{setup of depth & proj dim}. The key lies in the surjectivity of the restriction map $\sF(X)\to \sF(X\backslash Z)$ for an injective abelian sheaf $\sF$, a result that naturally follows by applying $\Hom(-,\sF)$ to the monomorphism $j_!\mathbf{Z}_{X\backslash Z}\to \mathbf{Z}_X$.  
 
 Although not directly pertinent to our discussion, let's briefly elucidate the connection with Zariski local cohomology presented in \S\S \ref{setup of depth & proj dim}. Consider the natural morphism $\mu: (\mathbf{Sch}{/X})_{\fppf} \to X_{\mathrm{Zar}}$ to the small Zariski site. For any abelian sheaf $\sF$ on $(\mathbf{Sch}_{/X})_{\fppf}$, there exists a natural comparison map $$R\Gamma_{\mathrm{Zar}}(X,\mu_*\sF)\to R\Gamma_{\mathrm{Zar}}(X,R\mu_*\sF)=R\Gamma_{\fppf}(X,\sF).$$ 
 While this is not an isomorphism in general, it holds true if $X$ is a scheme and $\sF$ is quasi-coherent, as $\mu_*\sF=R\mu_*\sF$  (quasi-coherent sheaves have vanishing higher fppf-cohomology on affines).  Considering the long exact sequence for local cohomology, a similar comparison result holds for local cohomology.}

On occasion, we will use the \'etale local cohomology. Its definition mirrors its counterpart, using either the big or small \'etale site of $X$ in lieu of the fppf site.
\eppt

The following result concerning \'etale descent of fppf local cohomology will be needed.
\blem [\emph{cf.}~ \cite{CS21}*{Lemma~7.1.1}]\label{et-desent-of-fppf-local}
For an algebraic space $X$, a closed subspace $Z \subset X$, and an abelian sheaf $\sF$ on $(\mathbf{Sch}_{/X})_{\fppf}$, if for any integer $q\ge 0$, $\wt{\cH}^q_Z(\sF)$ denotes the \'etale-sheafification of the presheaf
$(V\to X)\mapsto H^q_{Z}(V, \sF)$ where $V\to X$ is \'etale, then we have a convergent spectral sequence
\[
E_2^{pq}=H_{\et}^p(X,\wt{\cH}^q_Z(\sF))\Rightarrow H^{p+q}_Z(X,\sF).
\]
\elem
\bpf
There is a trouble with recursive references, so we give a proof for the convenience of the readers.

Let \textbf{Ab}$(-)$ denote the category of abelian sheaves on a site. Consider the following sequence of functors:
\[
\textbf{Ab}((\mathbf{Sch}_{/X})_{\fppf})
\xrightarrow{\cH_Z^0} \textbf{Ab}((\mathbf{Sch}_{/X})_{\fppf})
\xrightarrow{\nu_*} \textbf{Ab}((\textbf{\'EtSch}_{/X})_{\et})
\xrightarrow{\GG(X,-)}
\textbf{Ab},
\]
where $(\textbf{\'EtSch}_{/X})_{\et}$ is the site of \'etale schemes over $X$ with \'etale covers, and $\nu:(\mathbf{Sch}_{/X})_{\fppf}\to (\textbf{\'EtSch}_{/X})_{\et}$ denotes the natural morphism of sites.
The first two functors {send} injectives to injectives because they admit exact left adjoints. For any abelian sheaf $\sF$ on $(\mathbf{Sch}_{/X})_{\fppf}$, we have $(\nu_*\circ \cH_Z^0)(\sF)(V)=H_Z^0(V,\sF)$ where $V$ is a scheme \'etale over $X$. This implies that the $q$-th right derived functor of $\nu_*\circ \cH_Z^0$ is given by the \'etale sheafification of the presheaf $(V\to X)\mapsto H^q_{Z}(V, \sF)$.
Now the lemma follows from the Grothendieck spectral sequence applied to the functors $\nu_*\circ \cH_Z^0$ and $\GG(X,-)$.
\epf

  \bppt[Setup]\label{setup for purity over mult} From now on we assume the following, unless stated otherwise:
  \begin{itemize}
    \item let $X$ be a topologically locally Noetherian algebraic space, $Z\subset X$ a closed subspace, and $j:X\backslash Z\hra X$ the canonical open immersion;
    \item for every geometric point $\ov{z}\ra Z$, the strict local ring $\sO_{X,\ov{z}}^{\sh}$\footnote{ {By a geometric point of an algebraic space we refer to a map from the spectrum of a separably closed field. For a given geometric point $\ov{t}: \Spec(\Omega)\ra X$, the strict local ring $\sO_{X,\ov{t}}^{\sh}$ of $X$ at $\ov{t}$ is defined as the filtered colimit of all $B$, where $\Spec(B)\to X$ is an \'etale map, along with a geometric point $\ov{t}_B:\Spec(\Omega)\ra \Spec(B)$ lifting $\ov{t}$. The strict local ring $\sO_{X,\ov{t}}^{\sh}$, up to isomorphism, depends only on the equivalence class of the geometric point $\ov{t}$: indeed, if $\ov{t}':\Spec(\Omega')\ra \Spec(\Omega)\xrightarrow{\ov{t}} X$ is another geometric point and $\Spec(B)\to X$ is an \'etale map, there exists a natural bijection between liftings of $\ov{t}$ and $\ov{t}'$ to $\Spec(B)$. Clearly, this definition aligns with the classical definition of strict local rings for schemes.}}
     is (geometrically) parafactorial \footnote{Note that the topologically locally Noetherian assumption on $X$ implies that $\sO_{X,\ov{z}}^{\sh}$ has a quasi-compact punctured spectrum, thus it is geometrically parafactorial if and only if it is parafactorial.};
    \item let $M$ be an $X$-group algebraic space of multiplicative type, that is, its base change $M_{X'}$ is a group scheme of multiplicative type, where $X'\to X$ is some \'etale cover by a scheme $X'$.
  \end{itemize}
\eppt

\brem
In the case when $X$ is a scheme, \Cref{lim-parafactorial}~(ii) implies that the first two assumptions above are equivalent to saying that the pair $(X,Z)$ is geometrically parafactorial in the sense of \Cref{geom-parafact}. Moreover, \Cref{geom parafact}~(i)--(ii) gives examples of such pairs $(X,Z)$, where $Z$ satisfies the conditions
\[
   \begin{cases}
   \codim(Z_{\eta}, X_{\eta})\geq 2 & \text{for every generic point $\eta\in S$} \\
   \codim(Z_s,X_s)\geq 1 & \text{for all $s\in S$.}
   \end{cases}
\]

\erem

\bprop\label{local cohomology of tori}
In the Setup \ref{setup for purity over   mult}, assume that $X$ is a scheme. For a point $z\in Z$ and an $\sO_{X,z}$-torus $T$, we have
    \[
    \text{$H^i_{\{z\},\fppf}(\Spec(\sO_{X,z}),T)\simeq H^i_{\{z\},\et}(\Spec(\sO_{X,z}),T)=0$ \q for \, $0\le i \le 2$. }
    \]
\eprop
\bpf
The fppf and \'etale cohomology of a smooth group scheme coincides. So we may work with the \'etale site.
  By the local-to-global $E_2$-spectral sequence \cite{SGA4II}*{Exposé~V, Proposition~6.4},
\[
 H^p_{\et}(\Spec(\sO_{X,z}),\cH^q_{\{z\}}(T))\Rightarrow H^{p+q}_{\{z\}}(\Spec(\sO_{X,z}),T).
\]
Therefore, it suffices to prove that $\cH_{\{z\}}^q(T)=0$ for $0\le q \le 2$. Since $\sO_{X,z}$ has a quasi-compact punctured spectrum, we can identify their stalks at a geometric point $\overline{z}$ lying over $z$:
\[
\cH^q_{\{z\}}(T)_{\overline{z}}=H^q_{\{\overline{z}\}}(\Spec(\sO_{X,\overline{z}}^{\sh}),T).
\]
  Now, since $T_{\sO_{X,\overline{z}}^{\sh}} \simeq \mathbb{G}_{m,\sO_{X,\overline{z}}^{\sh}}^{\dim\,T}$ and $\sO_{X,\overline{z}}^{\sh}$ is parafactorial, we have
\[
H^q_{\et}(\Spec(\sO_{X,\overline{z}}^{\sh}),T)\simeq H^q_{\et}(\Spec(\sO_{X,\overline{z}}^{\sh})\backslash \{\ov{z}\},T) \q \text{ for } 0\le q\le 1.
\]
  Moreover, as $\sO_{X,\overline{z}}^{\sh}$ is strictly Henselian, we have
\[
H^2_{\et}(\Spec(\sO_{X,\overline{z}}^{\sh}),T)=0.
\]
Looking at the local cohomology exact sequence for the pair $(\Spec(\sO_{X,\overline{z}}^{\sh}),\overline{z})$ and $T$, we see that
$$
H_{\{\overline{z}\}}^q(\Spec(\sO_{X,\overline{z}}^{\sh}),T)=0 \q \text{ for } 0\le q\le 2,
$$
giving that $\cH_{\{z\}}^q(T)=0$ for $0\le q\le 2$, as desired.
\epf

  {The following result is a variant of \cite{CS21}*{Theorem~7.2.8~(a)}, where $X$ is assumed topologically locally Noetherian but the local rings $\sO_{X,z}$ are not supposed to be Noetherian for $z\in Z$.}
\bthm \label{purity for gp of mult type}
In the Setup \ref{setup for purity over   mult},
 we have
\[
\text{$H^i_{\fppf}(X,M)\isoto H^i_{\fppf}(X\backslash Z,M)$ for $i=0,1$ and $H^2_{\fppf}(X,M)\hra H^2_{\fppf}(X\backslash Z,M)$.}
\]
\ethm
\bpf
By the local cohomology exact sequence for the pair $(X,Z)$ and the sheaf $M$,
everything reduces to showing the vanishings $H^q_Z(X,M)=0$ for $0\le q \le 2$. By the spectral sequence in \Cref{et-desent-of-fppf-local},
it suffices to show the vanishings of $\wt{\cH}^q_Z(M)$, the \'etale-sheafification of the presheaf
\[
\text{$(V\to X)\mapsto H^q_{Z}(V, M)$, \qq where $V\to X$ is \'etale.}
\]
 In particular, the problem is \'etale-local on $X$, so we may pass to an \'etale cover to assume that $X$ is a scheme (noting that the assumptions of the Setup \Cref{setup for purity over   mult} still hold) and that $M$ splits as $\mu_n$ or $\bG_m$. Now,
since $\mu_n=\ker(\bG_m\overset{\times n}{\ra}\bG_m)$, it suffices to show that $\wt{\cH}^q_Z(\bG_m)=0$ for $0\le q\le 2$.

For $q=0,1$ this follows from the fact that the pair $(X, Z)$ is parafactorial (using \Cref{lim-parafactorial}). (This is where the Noetherian hypothesis is used.)

For $q=2$, by the case $q=0,1$ already proven and the long exact sequence for local cohomology, we have
\begin{align*}
  H^2_Z(X,\mathbb{G}_m) & \simeq \Ker(H^2_{\fppf}(X ,\bG_m) \to H^2_{\fppf}(X\backslash Z,\bG_m)) \\
   & \simeq \Ker(H^2_{\et}(X ,\bG_m) \to H^2_{\et}(X\backslash Z,\bG_m)).
\end{align*}

 The same is true for every scheme that is \'etale over $X$. Consequently, every class in  $H^2_Z(X,\mathbb{G}_m)$ vanishes in an \'etale cover of $X$, since this property holds for  $H^2_{\et}(X ,\bG_m)$. This implies that $\wt{\cH}_Z^2(\bG_m)=0$.
\epf

\csub[Grothendieck--Serre type results for groups of multiplicative type]

We record the following result from \cite{GL23}*{Proposition 3.6} for later use.

\bprop \label{G-S type results for mult type}
Let $R$ be a Pr\"ufer domain.
Consider an irreducible scheme $X$ essentially smooth over $R$ having function field $K(X)$, and an $X$-group scheme $M$ of multiplicative type. If there exists a connected finite \'etale Galois covering $X'\to X$ that splits $M$\footnote{Such a covering always exists, because $X$ is normal and so $M$ is isotrivial.}, then the restriction maps
 \[
 H^1_{\fppf}(X,M) \ra H^1_{\fppf}(K(X),M)\q \text{ and } \q H^2_{\fppf}(X,M) \ra H^2_{\fppf}(K(X),M)
 \]
  are injective in the following scenarios:
 \benumr
 \item \label{G-S for mult type gp}$X=\Spec (A)$ where $A$ is a semilocal ring that's essentially smooth over $R$;
 \item \label{B-Q for mult type} There exists an essentially smooth semilocal $R$-algebra $A$ such that $X$ embeds into $\overline{X}$ via a quasi-compact open immersion, with $\overline{X}$ being a smooth projective $A$-scheme having geometrically integral fibres. Furthermore, $\mathrm{Pic}(X_L)=0$ for any finite separable fields extension $L/\mathrm{Frac}(A)$, and $M=N_X$ where $N$ is an $A$-group of multiplicative type (e.g., $X$ could be a quasi-compact open subscheme of $\mathbb{P}_A^N$);
 \item For any \'etale covering $X''\to X$ {dominating} $X'\to X$, we have $\mathrm{Pic}(X'')=0$.
 \eenum
 Moreover, if $M$ is a flasque $X$-torus, then in all cases from $\mathrm{(i)}$--$\mathrm{(iii)}$, the restriction
 \[
 H^1_{\et}(X,M) \isoto H^1_{\et}(K(X),M) \q \text{ is bijective.}
 \]
\eprop

\section{Generically trivial torsors under a quasi-split group}
In this section, we study generically trivial torsors under quasi-split reductive group schemes. The main result is \Cref{qs-torsors}, which consists of \ref{qs-Nis-conj} a version of Nisnevich conjecture inspired by the recent preprint of $\check{\mathrm{C}}$esnavi$\check{\mathrm{c}}$ius \cite{Ces22b}*{Theorem~1.3~(2)}, who proved it in the case when $R$ is a Dedekind domain, and \ref{qs-GS-rk1} the Grothendieck--Serre conjecture over one-dimensional Pr\"{u}fer bases. The proof follows the strategy of \cite{Ces22a} (with its earlier version given by Fedorov \cite{Fed22b}), which is possible due to the availability of the main tools in our Pr\"uferian context, such as the toral version of purity (\emph{cf}. \cite{GL23}*{Theorem 3.3} and \Cref{purity for gp of mult type}) and the toral version of the Grothendieck--Serre conjecture (\emph{cf}. \Cref{G-S type results for mult type}~\ref{G-S for mult type      gp}).

\bthmt \label{qs-torsors}
For a semilocal Pr\"{u}fer domain $R$ with fraction field $K$, an {integral}, semilocal, and essentially smooth $R$-algebra $A$, and a quasi-split reductive $A$-group scheme $G$,
\benumr
\item \label{qs-Nis-conj} every generically trivial $G$-torsor over $A\otimes_RK$ is trivial, that is,
    \[
    \Ker \left(H^1_{\et}(A\otimes_RK,G)\to H^1_{\et}(\Frac A,G)\right) =\{*\};
    \]
\item  \label{qs-GS-rk1} if $R$ has Krull dimension 1, then every generically trivial $G$-torsor is trivial:
    \[
    \Ker \left(H^1_{\et}(A,G)\to H^1_{\et}(\Frac A,G)\right) =\{*\}.
    \]
\eenum
\ethmt

We start with the following consequence of \Cref{Ces's Variant 3.7}, which is the key geometric input permitting a series of reductions that eventually lead to \Cref{qs-torsors}.

\blemt [\emph{cf.}~\cite{Ces22a}*{Proposition~4.1}] \label{geome. input}
For
\benumr
  \item a semilocal Pr\"{u}fer domain $R$;
  \item a smooth, faithfully flat $R$-algebra $A$ of pure relative dimension $d\ge 1$;
  \item a finite subset $\mathbf{x}\subset X\ce \Spec\, A$;
  \item a closed subscheme $Y\subset X$ that satisfies
  \[
  \begin{cases}
   \codim(Y_s, X_s)\geq 1 & \text{for all closed points }s\in \Spec R, \text{ and} \\
   \codim(Y_s,X_s)\geq 2 & \text{otherwise}.
  \end{cases}
  \]
\eenum
there are an affine open $U\subset \Spec\, A$ containing $\mathbf{x}$, an affine open $S\subset \mathbb{A}_R^{d-1}$, and a smooth $R$-morphism $\pi:U\to S$ of relative dimension 1 such that $Y\cap U$ is $S$-finite.
\elemt
\bpf
 Choosing an embedding of $X$ into some affine space over $R$ and taking schematic closure in the corresponding projective space, we get a projective compactification $\overline{X}$ of $X$. Since $\overline{X}$ is flat and projective over $R$, by \Cref{geom}~\ref{geo-i}, all its $R$-fibres have the same dimension $d$. Denoting by $\overline{Y} \subset \overline{X}$ the schematic closure of $Y$,
to apply {\Cref{Ces's Variant 3.7}} and conclude (in which $X$ is $\overline{X}$ here, $W$ is $X$ here, and $Y$ is $\overline{Y}$ here), we need to check that
the boundary $\overline{Y}\backslash Y$ is $R$-fibrewisely of codimension $\ge 2$ in $\overline{X}$.

By \SP{01RG}, for a quasi-compact immersion of schemes, the schematic closure has the underlying space the topological closure. Thus, set-theoretically we have $\overline{Y}=\bigcup_y \overline{\{y\}}$, where $y$ runs through the (finitely many) generic points of $Y$.

Let $y\in Y$ be a generic point, lying over $s\in \Spec R$.
 By \Cref{geom}~\ref{geo-i}, $\overline{X}$ has equal $R$-fibre dimension $d$ and all non-empty $R$-fibres of $\overline{\{y\}}$ have the same dimension. If $s$ is not a closed point, then for any specialization $s\specializes s'\in \Spec R$, we have
\[
\codim(\ov{\{y\}}_{s\pr}, \overline{X}_{s\pr})= \codim(\ov{\{y\}}_{s}, \overline{X}_{s})\geq 2;
\]
\emph{a fortiori}, the contribution of such a generic point $y$ to the $R$-fibre codimension of $\overline{Y}\backslash Y$ in $\overline{X}$ is $\ge 2$.

Otherwise, $s$ is a closed point, then $\overline{\{y\}}_{s} =\overline{\{y\}} \subset \overline{Y_{s}}$. 
As assumed, $\codim(\overline{Y_{s}},\overline{X}_{s})=\codim(Y_{s},X_{s})\geq 1$, so we have $\codim(\overline{\{y\}}_{s}, \overline{X}_{s}) \ge 1$. But since the generic point $y$ of $\overline{\{y\}}_{s}$ is not contained in $ \overline{Y}\backslash Y$, we deduce that the contribution of such a generic point $y$ to the $s$-fibre codimension of $\overline{Y}\backslash Y$ in $\overline{X}$ is again $\ge 2$.
\epf

\blemt [Lifting the torsor to a smooth relative curve; \emph{cf.}~\cite{Ces22a}*{Proposition~4.2}]\label{red q-s torsor to rel curve}
Let $R$ be a semilocal Pr\"{u}fer domain with fraction field $K$.
Let $A$ be the semilocalization of an irreducible, $R$-smooth algebra $A'$ at a finite subset $\textbf{x}\subset \Spec A\pr$ and let $G$ be a quasi-split reductive $A$-group scheme with a Borel subgroup $B$.
\begin{itemize}
\item [(1)] Given a generically trivial $G$-torsor $P_K$ over $A_K\ce A\otimes_RK$, there are
\benumr
\item a smooth, affine relative $A$-curve $C$ with a section $s\in C(A)$;
\item an $A$-finite closed subscheme $Z\subset C$;
\item a quasi-split reductive $C$-group scheme $\sG$ with a Borel subgroup $\sB\subset \sG$ whose $s$-pullback is $B\subset G$, compatible with the quasi-pinnings;
\item a $\sG$-torsor $\cP_K$ over $C_K\ce C\times_RK$ whose $s_{A_K}$-pullback is $P_K$ such that $\cP_K$ reduces to a {$\rad^u(\sB)$}-torsor over $C_K\backslash Z_K$ (here $s_{A_K}$ denotes the image of $s$ in $C(A_K)$).
\eenum
\item [(2)] If $R$ has Krull dimension $1$ and $P$ is a generically trivial $G$-torsor, there are
\benumr
\item a smooth, affine relative $A$-curve $C$ with a section $s\in C(A)$;
\item an $A$-finite closed subscheme $Z\subset C$;
\item a quasi-split reductive $C$-group scheme $\sG$ with a Borel subgroup $\sB\subset \sG$ whose $s$-pullback is $B\subset G$, compatible with the quasi-pinnings;
\item a $\sG$-torsor $\cP$ whose $s$-pullback is $P$ such that $\cP$ reduces to a {$\rad^u(\sB)$}-torsor over $C\backslash Z$.
\eenum
\end{itemize}
\elemt

\bpf
In the case (1) we can first use a limit argument involving \Cref{approxm semi-local Prufer ring} to reduce to the case when $R$ has finite Krull dimension.

If $A'$ is of relative dimension 0 over $R$, then $A_K=\Frac (A)$ and $A$ is a semilocal Pr\"{u}fer domain. Thus, $P_K$ is trivial, and, by the Grothendieck--Serre conjecture on semilocal Pr\"{u}fer schemes (\emph{cf}. \cite{GL23}*{Theorem A.0.1}), $P$ is also trivial. In this case we simply take $C=\mathbb{A}_A^1$, $s=0\in \mathbb{A}_A^1(A)$, $Z=\emptyset$, $(\sG,\sB)=(G_{\mathbb{A}_A^1},B_{\mathbb{A}_A^1})$, and $\cP_K=(P_K)_{\mathbb{A}_{A_K}^1}$ (resp., $\cP=P_{\mathbb{A}_A^1}$). Thus, for what follows, we can assume that the relative dimension of $A'$ over $R$ is $d>0$.

By spreading out and localizing $A'$, we may assume that our quasi-split $G$ (in particular, the Borel $B$) and torsor $P$ all live over $A'$, and $P_K$ {lives} over $A'_K$. By \cite{SGA3IIInew}*{Exposé~XXVI, Corollaire~3.6 and Lemme~3.20}, the quotient $P_K/B_K$ (resp., $P/B$) is representable by a smooth projective scheme over $A'_K$ (resp., over $A'$). Now we treat the cases (1)--(2) separately.

(1) Since $P_K$ is generically trivial, the valuative criterion of properness applies to $P_K/B_K\to \Spec A\pr_K$, we find a closed subset $Y_K\subset \Spec A\pr_K$ of codimension $\ge 2$ such that $P_K/B_K \to \Spec A\pr_K$ has a section over $\Spec A\pr_K\backslash Y_K$ that lifts to a generic section of $P_K$.
  In other words, $(P_K)_{\Spec A\pr_K\backslash Y_K}$ reduces to a generically trivial $B_{\Spec A\pr_K\backslash Y_K}$-torsor $P_K^B$.
  Consider the $A'$-torus $T\ce B/\rad^u(B)$ and the induced
 \[
 \text{$T$-torsor \q $P_K^T\ce P_K^B/\rad^u(B)_K \q$  over  \q $\Spec (A'_K)\backslash Y_K$.}
 \]
Since $P_K^T$ is generically trivial, by \Cref{extend generically trivial torsors}, it extends to a $T$-torsor $\widetilde{P_K^T}$ over $\Spec  A\pr \backslash F$ for a closed subscheme $F\subset \Spec  A\pr$ satisfying
\[
\text{$\codim(F_K, \Spec A\pr_K)\geq 3$ \q and\q  $\codim(F_s,\Spec A\pr_s)\geq 2$ for all $s\in \Spec(R)$;}
\]
by purity for tori (\emph{cf}. \cite{GL23}*{Theorem 3.3} or \Cref{purity for gp of mult type}), this torsor further extends to the whole $\Spec A\pr$. As $\widetilde{P_K^T}$ is generically trivial, by the Grothendieck--Serre conjecture for tori (\Cref{G-S type results for mult type}~\ref{G-S for mult type      gp}), we may localize $A'$ around $\textbf{x}$ to assume that $\widetilde{P_K^T}$, and hence also $P_K^T$, is already trivial. In other words, $(P_K)_{\Spec (A'_K)\backslash Y_K}$ reduces to a $\rad^u(B)$-torsor over $\Spec (A'_K)\backslash Y_K$.

Denote by $Y$ the schematic closure of $Y_K$ in $\Spec A\pr$; by \Cref{geom}~\ref{geo-i}, it is $R$-fibrewisely of codimension $\ge 2$ in $\Spec  A\pr$.
Applying \Cref{geome. input} to the $R$-smooth algebra $A'$ and the closed subscheme $Y\subset \Spec  A\pr$, we obtain an affine open $U\subset \Spec  A\pr$ containing $\mathbf{x}$, an affine open $S\subset \mathbb{A}_R^{d-1}$, and a smooth $R$-morphism $\pi:U\to S$ of relative dimension 1 such that $Y\cap U$ is $S$-finite.

Recall that $A$ is the semilocal ring of $U$ at $\textbf{x}$. Denote
\[
   C\ce U\times_S\Spec A \q \text{ and } \q Z\ce (Y\cap U)\times_S\Spec A.
\]
Then $C$ is a smooth affine relative $A$-curve, the diagonal in $C$ induces a section $s\in C(A)$, and the closed subscheme $Z\subset C$ is $A$-finite. Thus (1)(i) and (1)(ii) hold. Let $\sB\subset \sG$ be the pullback of $B_U\subset G_U$ under the first projection $\text{pr}_1:C\to U$, and let $\cP_K$ be the pullback of $(P_K)_{U_K}$ under the first projection $\text{pr}_1:C_K\to U_K$.
Then, $\cP_K$ is a $\sG$-torsor over $C_K$, and, by construction, the $s$-pullback (resp., $s_{A_K}$-pullback) of $\sB\subset\sG$ (resp., of $\cP_K$) is $B\subset G$ (resp., $P_K$). Finally, since $P_K$ reduces to a $\rad^u(B)$-torsor over $\Spec (A'_K)\backslash Y_K$, $\cP_K$ reduces to a $\rad^u(\sB)$-torsor over $C_K\backslash Z_K$. Thus (1)(iii) and (1)(iv) also hold.

 (2) Recall that, by \Cref{geom}~\ref{geo-iii}, the local rings of all maximal points of $R$-fibres of $\Spec  A\pr$ are valuation rings. As $P$ is generically {trivial}, applying the valuative criterion of properness to $P/B\to \Spec  A\pr$ yields a closed subscheme $Y\subset \Spec  A\pr$, which avoids all the codimension $1$ points of the generic fibre $\Spec A\pr_K$ and all the maximal points of $R$-fibres of $\Spec  A\pr$. Moreover, $P/B \to \Spec  A\pr$ has a section over $\Spec  A\pr\backslash Y$ that lifts to a generic section of $P$. In other words, $Y$ satisfies
\[
\text{$\codim(Y_K, \Spec A\pr_K)\geq 2$ \q and\q  $\codim(Y_s,\Spec A\pr_s)\geq 1$ for all $s\in \Spec(R)$.}
\]
  Therefore, $P_{\Spec  A\pr\backslash Y}$ reduces to a generically trivial $B_{\Spec  A\pr\backslash Y}$-torsor $P^B$. Consider the $A'$-torus $T\ce B/\rad^u(B)$ and the induced $T$-torsor
 \[
 P^T\ce P^B/\rad^u(B)  \q \text{ over } \q \Spec  A\pr\backslash Y.
 \]
By purity for tori (\emph{cf}. \cite{GL23}*{Theorem 3.3} or \Cref{purity for gp of mult type}), $P^T$ extends to {a $T$-torsor $\widetilde{P^T}$ over $\Spec  A\pr$}. As $\widetilde{P^T}$ is generically trivial, by the Grothendieck--Serre conjecture for tori (\Cref{G-S type results for mult type}~\ref{G-S for mult type      gp}), we may localize $A'$ around $\textbf{x}$ to assume that $\widetilde{P^T}$, and hence also $P^T$, is already trivial. In other words, $P_{\Spec  A\pr\backslash Y}$ reduces to a $\rad^u(B)_{\Spec  A\pr\backslash Y}$-torsor.

Now, applying \Cref{geome. input} to the $R$-smooth algebra $A'$ and the closed subscheme $Y\subset \Spec  A\pr$, we obtain an affine open $U\subset \Spec  A\pr$ containing $\mathbf{x}$, an affine open $S\subset \mathbb{A}_R^{d-1}$, and a smooth $R$-morphism $\pi:U\to S$ of relative dimension 1 such that $Y\cap U$ is $S$-finite.

Recall that $A$ is the semilocal ring of $U$ at $\textbf{x}$. Denote $$C\ce U\times_S\Spec A \q \text{ and } \q Z\ce (Y\cap U)\times_S\Spec A.$$
Then $C$ is a smooth affine relative $A$-curve, the diagonal in $C$ induces a section $s\in C(A)$, and the closed subscheme $Z\subset C$ is $A$-finite. So (2)(i) and (2)(ii) hold. Let $\sB\subset \sG$ and $\cP$ be the pullback of $B_U\subset G_U$ and $P_U$ under the first projection $\text{pr}_1:C\to U$, respectively. Then, $\cP$ is a $\sG$-torsor over $C$, and, by construction, the $s$-pullback of $\sB\subset\sG$ and $\cP$ are $B\subset G$ and $P$, respectively. Finally, since $P$ reduces to a $\rad^u(B)$-torsor over $\Spec  A\pr\backslash Y$, $\cP$ reduces to a $\rad^u(\sB)$-torsor over $C\backslash Z$. Thus, (2)(iii) and (2)(iv) also hold.
\epf

\blemt [{\cite{Ces22a}*{Lemma~5.1}}] \label{equating q-s group}
For a semilocal ring $A$ whose local rings are geometrically unibranch, an ideal $I\subset A$, reductive $A$-groups $G$ and $G'$ that on geometric $A$-fibres have the same type, fixed quasi-pinnings of $G$ and $G'$ extending Borel $A$-subgroups $B\subset G$ and $B'\subset G'$ and an $A/I$-group isomorphism
\[
\text{$\iota\colon G_{A/I} \isoto G_{A/I}\pr$ \, respecting the quasi-pinnings; in particular,  \, $\iota(B_{A/I})=B_{A/I}'$,}
\]
there are
\benumr
\item a faithfully flat, finite, \'etale $A$-algebra $\wt{A}$ equipped with an $A/I$-point 
\[
   a\colon \wt{A} \twoheadrightarrow A/I; \text{}
\]
\item an $\wt{A}$-group isomorphism $\wt{\iota}\colon G_{\wt{A}} \isoto G_{\wt{A}}\pr$ respecting the quasi-pinnings such that $a^*(\wt{\iota})=\iota$.
\eenum
\elemt

Note that the original formulation \cite{Ces22a}*{Proposition~5.1} assumed that $A$ is Noetherian, though the Noetherianness of $A$ was not used in the proof.

\blemt [Changing the relative curve $C$ to equate $\sG$ and $G_C$; \emph{cf.}~\cite{Ces22a}*{Proposition~5.2}]
\label{Changing the rel cur to equ.}
In the setting of \Cref{red q-s torsor to rel curve}, for both cases (1) and (2), we may replace $C$ by an \'etale neighborhood of $\im(s) $ to further achieve that 
\[
   (\sG,\sB)=(G_C,B_C).
\]
\elemt
\bpf
Consider the semilocalization $\Spec (D)$ of $C$ at the closed points of $\im(s)\cup Z$. Since $C$ is normal, all the local rings of $D$ are geometrically unibranch. The image of the section $s:\Spec A\to \Spec (D)$ gives rise to a closed subscheme $\Spec (D/I) \subset \Spec (D)$. By the conclusion of \Cref{red q-s torsor to rel curve}, the restrictions of $\sB_D \subset \sG_D$ and $B_D\subset G_D$ to $\Spec (D/I)$ agree with each other in a way respecting their quasi-pinnings. Thus, by \Cref{equating q-s group}, there is a faithfully flat, finite, \'etale $D$-algebra $\wt{D}$, a point $\wt{s}:\wt{D} \twoheadrightarrow D/I \simeq A$ lifting $s:D \twoheadrightarrow D/I \simeq A$ such that $\sB_{\wt{D}} \subset \sG_{\wt{D}}$ is isomorphic to $B_{\wt{D}}\subset G_{\wt{D}}$ compatibly with the fixed identification of $\wt{s}$-pullbacks. We then spread out the finite \'etale morphism $\Spec (\wt{D})\to \Spec (D)$ to obtain a finite \'etale morphism $\wt{C}\to C'$ for an open $C'\subset C$ that contains $\im(s)\cup Z$, while preserving an $\wt{s}\in \wt{C}(A)$, and an isomorphism between $\sB_{\wt{C}} \subset \sG_{\wt{C}}$ and $B_{\wt{C}}\subset G_{\wt{C}}$. It remains to replace $C,s,Z$ and $\cP_K$ (resp., $\cP$) by $\wt{C}$, $\wt{s}$, $Z\times_C\wt{C}$ and $(\cP_K)_{\wt{C}_K}$ (resp., $\cP_{\wt{C}}$).
\epf

\blemt [Changing the smooth relative curve $C$ for descending to $\mathbb{A}_A^1$; \cite{Ces22a}*{Proposition~6.5}]
\label{prepare for des. to A^1}
In the setting of \Cref{red q-s torsor to rel curve}, for both cases (1) and (2), in addition to $(\sG,\sB)=(G_C,B_C)$,
we may change $C$ to further achieve that there is a flat $A$-map $C\to \mathbb{A}_A^1$ that maps $Z$ isomorphically onto a closed $Z'\subset \mathbb{A}_A^1$ such that
    \[
    Z\simeq Z' \times_{\mathbb{A}_A^1}C.
    \]
\elemt
\bpf
Assume that, in both cases (1) and (2) of \Cref{red q-s torsor to rel curve}, we have achieved the conclusion of \Cref{Changing the rel cur to equ.}. We have the data of a smooth affine relative $A$-curve $C$, a section $s\in C(A)$, and an $A$-finite closed subscheme $Z\subset C$.
Replacing $Z$ by $Z\cup \text{im}(s)$, we may assume that $s$ factors through $Z$. Unfortunately, in general, the $A$-finite scheme $Z$ may be too large to embed into $\mathbb{A}_A^1$.
(For instance, if $R=k$ is a finite field, then $Z$ cannot be embedded into $\mathbb{A}_k^1$ as soon as $\# \,Z(k) > \# \,k$.)
To overcome this difficulty, we first apply {Panin's `finite fields tricks' \cite{Ces22a}*{Lemma~6.1}} to obtain a finite morphism $\wt{C}\to C$ that is \'etale at the points in $\tilde{Z}\ce \wt{C}\times_CZ$ such that $s$ lifts to $\wt{s}\in \wt{C}(A)$, and there are no finite fields obstruction to embedding $\wt{Z}$ into $\mathbb{A}_A^1$ in the following sense: for every maximal ideal $\mathfrak{m}\subset A$ and every $d\geq 1$,
\[
\# \bigl\{
z\in \wt{Z}_{\kappa(\mathfrak{m})}\colon [\kappa(z)\colon \kappa(\mathfrak{m})]=d \bigr\}< \# \bigl\{z\in \mathbb{A}_{\kappa(\mathfrak{m})}^1\colon [\kappa(z)\colon\kappa(\mathfrak{m})]=d \bigr\}.
\]
Then, by \cite{Ces22a}*{Lemma~6.3}, there are an affine open $C'\subset \wt{C}$ containing $\text{im}(\wt{s}) $, a quasi-finite, flat $A$-map $C'\to \mathbb{A}_A^1$ that maps $Z$ isomorphically onto a closed subscheme $Z'\subset \mathbb{A}_A^1$ with
    \[
    Z\simeq Z' \times_{\mathbb{A}_A^1}C'.
    \]
It remains to replace $C$ by $C'$, $Z$ by $\wt{Z}$, $s$ by $\wt{s}$, and $\cP_K$ by $(\cP_K)_{C'_K}$ (resp., $\cP$ by $\cP_{C'}$).
\epf

\blemt [Descend to $\mathbb{A}_A^1$ via patching; \emph{cf.}~\cite{Ces22a}*{Proposition~7.4}]
\label{Desc. to A^1 via patching}
In the setting of \Cref{red q-s torsor to rel curve}, for both cases (1) and (2), we may further achieve that
\[
   \text{$(\sG,\sB)=(G_C,B_C)$, \q $C=\mathbb{A}_A^1$, \q  and $s=0\in \mathbb{A}_A^1(A)$.}
\]
\elemt
\bpf
By the reduction given in \Cref{prepare for des. to A^1}, we have a flat $A$-curve $C$, a section $s\in C(A)$, an $A$-finite closed subscheme $Z\subset C$, a quasi-finite, affine, flat $A$-map $C\to \mathbb{A}_A^1$ that maps $Z$ isomorphically onto a closed subscheme $Z'\subset \mathbb{A}_A^1$ such that $Z= Z' \times_{\mathbb{A}_A^1}C$, and
a $G$-torsor $\cP_K$ over $C_K$ whose $s_{A_K}$-pullback is $P_K$ (resp., a $G$-torsor $\cP$ over $C$ whose $s$-pullback is $P$) and whose restriction to $C_K\backslash Z_K$ (resp., $C\backslash Z$) reduces to a $\rad^u(B)$-torsor. Now, since $Z= Z' \times_{\bA_A^1}C\simeq Z'$, \cite{Ces22a}*{Lemma~7.2} (the Noetherian hypothesis is not needed) implies the pullback maps
\[
H^1_{\et}(\bA_A^1\backslash Z',\rad^u(G))\twoheadrightarrow H^1_{\et}(C\backslash Z,\rad^u(G))
\]
and
\[
H^1_{\et}(\bA_{A_K}^1\backslash Z_K',\rad^u(G))\twoheadrightarrow H^1_{\et}(C_K\backslash Z_K,\rad^u(G))
\]
are surjective. Combining these, we see that $\cP_K|_{C_K\backslash Z_K}$ (resp., $\cP|_{C\backslash Z}$) descends to a $G$-torsor $\cQ_K$ (resp., $\cQ$) over $\bA_{A_K}^1\backslash Z_K'$ (resp., $\bA_A^1\backslash Z'$) that reduces to a $\rad^u(B)$-torsor. By \cite{Ces22a}*{Lemma~7.1}, we may (non-canonically) glue $\cP_K$ with $\cQ_K$ (resp., $\cP$ with $\cQ$) to descend $\cP_K$ (resp., $\cP$) to a $G$-torsor $\wt{\cP_K}$ (resp., $\wt{\cP}$) over $\bA_{A_K}^1$ (resp., over $\bA_{A}^1$) that reduces to a $\rad^u(B)$-torsor over $\bA_{A_K}^1\backslash Z_K'$ (resp., over $\bA_A^1\backslash Z'$).
It remains to replace $C$ by $\bA_A^1$, $Z$ by $Z'$, $s\in C(A)$ by its image in $\bA_A^1(A)$, and $\cP_K$ by $\wt{\cP_K}$ (resp., $\cP$ by $\wt{\cP}$). Finally,  by shifting, we may assume that $s=0\in \bA_A^1(A)$.
\epf

\bpf[Proof of \Cref{qs-torsors}]
Let $P_K$ (resp., $P$) be a generically trivial $G_{A_K}$-torsor (resp., $G$-torsor). By the reduction \Cref{Desc. to A^1 via patching}, we get an $A$-finite closed subscheme $Z\subset \bA^1_A$, and a $G_{\mathbb{A}_{A_K}^1}$-torsor $\cP_K$ (resp., $G_{\bA_A^1}$-torsor $\cP$) whose pullback along the zero section is $P_K$ (resp., $P$) such that $(\cP_K)|_{\bA_{A_K}^1 \backslash Z_K}$ (resp., $\cP|_{\bA_A^1 \backslash Z}$) reduces to a $\rad^u(B)$-torsor. Since any $A$-finite closed subscheme of $\bA_A^1$ is contained in $\{f=0\}$ for some monic polynomial $f\in A[t]$, we may enlarge $Z$ to assume that $\bA_A^1\backslash Z$ is affine, to the effect that any $\rad^u(B)$-torsor over $\bA_{A_K}^1\backslash Z_K$ (resp., over $\bA_A^1\backslash Z$), such as $(\cP_K)|_{\bA_{A_K}^1 \backslash Z_K}$ (resp., $\cP|_{\bA_A^1 \backslash Z}$), is trivial. By section theorem \cite{GL23}*{Theorem 5.1}, the pullback of $\cP_K$ (resp., of $\cP$) along the zero section is trivial, that is, $P_K$ (resp., $P$) is trivial, as desired.
\epf
\newpage
\begin{appendix}

\section{Regular coherent rings} \label{apdxA}
In this appendix, we delve into the homological properties of coherent regular rings.
A coherent ring is \emph{regular} if its every finitely generated ideal has finite projective dimension.
Localizations of a coherent regular ring remain coherent regular. This is due to the stability of coherence under localization, combined with the fact that every finitely generated ideal in the localization is the localization of a finitely generated ideal. Moreover, one can verify coherent regularity Zariski locally, as can be deduced from the following more general result.

\blemt [Faithfully flat descent for coherent regularity] \label{ff descent coh reg}
For a faithfully flat ring map $A\to B$, if $B$ is coherent regular, then so is $A$.
\elemt
\bpf
The coherence of $A$ follows since the property of being finitely presented satisfies faithfully flat descent.
To see the regularity of $A$, we let $I\subset A$ be a finitely generated ideal and pick a resolution $P_{\bullet}\to I$ with each $P_i$ finite free over $A$ (using the coherence of $A$). Since $A\to B$ is flat, $P_{\bullet}\otimes_AB\to I\otimes_AB\simeq IB$ is a resolution of the $B$-module $IB$. As $B$ is regular, $\mathrm{Im}(P_{n+1}\otimes_AB\to P_{n}\otimes_AB)\simeq \mathrm{Im}(P_{n+1}\to P_n)\otimes_AB$ is finite projective over $B$ for some $n\ge 0$, and so $\mathrm{Im}(P_{n+1}\to P_n)$ is finite projective over $A$ by faithfully flat descent.
\epf

\blemt [\'Etale-local nature of coherent regularity] \label{et loc coh reg}
Let $A\to B$ be an \'etale ring map.
If $A$ is coherent regular, then so is $B$. The converse holds if $A\to B$ is faithfully flat and \'etale.
\elemt
\bpf
In light of \Cref{ff descent coh reg}, it is enough to show the coherent regularity of $B$. By Zariski descent, we may work Zariski locally to assume that $B$ is a (principal) localization of a finite free $A$-algebra $C$ (using the local structure of \'etale algebras). In this scenario, every finitely generated ideal of $C$ is finitely presented over $A$ and so is it over $C$.
Thus, $C$ is coherent, and so is its localization $B$.
Now, to show the regularity of $B$, our goal is to prove that $\fldim_B I<\infty$ for every finitely generated ideal $I\subset B$.

There exists a finitely generated ideal $J\subset C$ such that $J\otimes_CB\simeq JB=I$.
As $A$ is coherent regular, we have $\fldim_A J<\infty$. But $I$ is a filtered colimit of copies of $J$, so $\fldim_A I<\infty$. Since $B$ is $A$-flat, choosing a partial flat resolution of the $B$-module $I$, we are reduced to the following claim (\emph{cf.}~\SP{05B9}). \epf
\bcl
If $A\to B$ is an \'etale ring map and $M$ is a $B$-module which is $A$-flat, then $M$ is $B$-flat.
\ecl
\bpf[Proof of the claim]
We will argue by induction on the supremum of the cardinality of the geometric points in the fibres of $\Spec B\to \Spec A$. Since flatness can be checked Zariski locally, we may localize to assume that $\Spec B\to \Spec A$ is \'etale surjective. By the faithfully flat descent of flatness, we can replace $A$ with $B$, $B$ with $B\otimes_AB$, and $M$ with $M\otimes_AB$, thus assuming that $A\to B$ possesses a retraction $B\to A$. In this scenario, we have $B\simeq A \times B_1$, where $B_1$ is an \'etale $A$-algebra. Consequently, $M\simeq M_0 \times M_1$, with $M_0=M\otimes_BA$ and $M_1=M\otimes_BB_1$. By assumption, both $M_0$ and $M_1$ are $A$-flat.

Now, since the supremum of the cardinality of the geometric points in the fibres of $\Spec B_1\to \Spec A$ is strictly less than that of $\Spec B\to \Spec A$, our induction hypothesis gives the $B_1$-flatness of $M_1$, and consequently, the $B$-flatness of $M$.
\epf

\blemt \label{fibre criterion coh reg}
Let $f\colon A\to B$ be a flat local map of local rings and $\kappa_A$ the residue field of $A$.
Assume that $B$ is coherent. Let $M$ be a finitely presented $B$-module.
\benumr
\item If $M$ is $A$-flat, then we have $\mathrm{proj. dim}_B(M)\le \mathrm{fl.}\dim_{B\otimes_A\kappa_A}(M\otimes_A\kappa_A)$.
\item In general, we have $\mathrm{proj. dim}_B(M) \le \mathrm{fl.}\dim_A(M)+\wdim(B\otimes_A\kappa_A)$. Consequently, we have
\[
  \wdim(B) \le \wdim(A)+\wdim(B\otimes_A\kappa_A).
\]
\item  {If every finitely presented $B$-module, considered as an $A$-module, and} every finitely presented $B\otimes_A\kappa_A$-module have finite flat dimensions, then $B$ is coherent regular.
\eenum
\elemt
\bpf
For (i), set $l\ce \mathrm{fl.}\dim_{B\otimes_A\kappa_A}(M\otimes_A\kappa_A)$.
If $l=\infty$, there is nothing to show.
Otherwise, we choose a partial resolution over $B$
\[
  0\to M\pr \to P_{l-1}\to \cdots \to  P_0\to M \to 0
\]
with each $P_i$ finite free and $M\pr$ finitely presented (\Cref{lem-coh}~\ref{ker-coker-coh}).
By assumption, $M$ and $B$ are $A$-flat, so by \SP{03EY}, this is an $A$-flat resolution of $M$.
Tensoring it with $\kappa_A\otimes_A(-)$ gives a partial flat resolution of the $B\otimes_A\kappa_A$-module $M\otimes_A\kappa_A$.
The latter has flat dimension $l$, so that $M'\otimes_A\kappa_A$ is flat hence free over $B\otimes_A\kappa_A$.
Since $M\pr$ is $A$-flat, by \cite{EGAIV3}*{Proposition~11.3.7} and Nakayama's lemma, $M\pr$ is free over $B$ of the same rank.
This proves that $\text{proj.}\dim_B(M)\le l$, as desired.

In general, although $M$ may not be $A$-flat, we can first choose a partial resolution
$$
0\to M\prpr\to F_{s-1} \to \cdots \to F_0 \to M\to 0
$$
with each $F_i$ finite free over $B$ and $s\ce \mathrm{fl.}\dim_A(M)$ (if it is finite). Then, the finitely presented $B$-module $M\prpr$ is $A$-flat, so we can apply the part (i) to $M\prpr$ and obtain
\[
\mathrm{proj. dim}_B(M) \le s +\mathrm{proj. dim}_B(M\prpr) \le s+ \mathrm{fl.}\dim_{B\otimes_A\kappa_A}(M\prpr\otimes_A\kappa_A).
\]
This simultaneously proves the assertions (ii)--(iii).
\epf
Combining \Cref{fibre criterion coh reg} and \Cref{ring property coh reg}~\ref{qc coh reg has fin wdim}, we obtain the following.
\bcort\label{reg-reg-reg}
Let $f\colon A\to B$ be a local map of local rings and $\kappa_A$ the residue field of $A$.
Assume that
\benumr
\item $f$ is flat,
\item $A$ has finite weak dimension,
\item $B\otimes_A\kappa_A$ is coherent regular, and
\item $B$ is coherent,
\eenum
then $B$ is regular. The same conclusion holds if $\mathrm{(ii)}$ is replaced by
\begin{itemize}
  \item [$\mathrm{(ii)'}$] $A$ is coherent regular with a quasi-compact punctured spectrum.
\end{itemize}
\ecort

\begt
Let $A\to B$ be a flat, regular ring map of coherent rings.
If wdim$(A)<\infty$ (this implies that $A$ is regular; \emph{e.g.} when $A$ is a valuation ring), then $B$ is also coherent regular.
\eegt

Finally, we aim to prove \Cref{Auslander-Bachsbaum formula}, which serves as an extension of the classical Auslander--Buchsbaum formula to general non-Noetherian rings, drawing parallels with the conventional regular scenario established in \cite{AB57}*{Theorem~3.7}.

The following simple lemma from linear algebra will be used in the proof of \Cref{Auslander-Bachsbaum formula}.
\blemt \label{a lem on linear algebra}
For a local ring $(A,\mathfrak{m}_A)$, a nonzero $A$-module $M$ supported on $\{\mathfrak{m}_A\}$, and a matrix $H\in \text{M}_{m\times n}(A)$, if the $A$-linear map $H_M:M^{\oplus n}\to M^{\oplus m}$ induced by $H$ (via left multiplication) is injective, then $H$ admits a left inverse, or, equivalently, $H$ exhibits $A^{\oplus n}$ as a direct summand of $A^{\oplus m}$.
\elemt
\bpf
Recall \SP{0953} that the assumption on the support of $M$ means that, for any $w\in M$ and any finitely generated ideal $I\subset A$, we have $I^N\cdot w=0$ for large enough $N$.
Denote $H=(h_{ij})$.
We observe that at least one of $h_{ij}$ is invertible.
Otherwise, the entries $h_{ij}$ generate a proper ideal $I$ of $A$; pick a nonzero element $w\in M$ and let $N\ge 0$ be the smallest integer such that $I^N\cdot w\neq 0$, then $H_M((I^N\cdot w)^{\oplus n})=0$, contradicting our assumption that $H_M$ is injective.
Without loss of generality, we may assume that $h_{11} \in A^{\times}$. By subtracting suitable multiples of the first row of $H$ to other rows (resp. the first column of $H$ to other columns), we may also assume that $h_{1j}=0$ for $1<j\le n$ and $h_{i1}=0$ for $1<i\le m$ (the assumption and conclusion of the lemma are preserved if we replace $H$ by $H_1HH_2$, where $H_1\in M_{m\times m}(A)$ and $H_2\in M_{n\times n}(A)$). In other words, we have $H=(h_{11}) \oplus H'$, where $H'\in M_{(m-1)\times (n-1)}(A)$. Then the map  $H'_{M}: M^{\oplus (n-1)}\to M^{\oplus (m-1)}$ induced by $H'$ is also injective. So we may assume by induction that $H'$ admits a left inverse $H''\in M_{(n-1)\times(m-1)}(A)$. Then $(h_{11}^{-1}) \oplus H''$ is a left inverse of $H$.
\epf
{As L. Moret-Bailly pointed out, since $\Ker(H_M)=\Hom_A(\Coker(H^{t}), M)$, where $H^t$ is the transpose of $H$, the above lemma is also a direct consequence of \Cref{HomAss-local}.}

 Now, we acquire the the Auslander--Buchsbaum formula (\emph{cf.}~\cite{AB57}*{Theorem~3.7}) for general local rings.
{Since our depth is not well-behaved over non-quasicompact punctured spectra, we need to assume this condition in the sequel.}
    \bthmt[Auslander--Buchsbaum formula]\label{Auslander-Bachsbaum formula}
    \hfill
    \benumr
    \item Let $(A,\fm_A)$ be a local ring with a quasi-compact punctured spectrum. Let $M$ be an $A$-module having a finite resolution by finite free $A$-modules. Then we have
        \[
        \mathrm{proj.}\dim_A(M)+\depth_A(M)=\depth_A(A).
        \]
    \item Let $V$ be a valuation ring such that $s:=\fm_V$ is the radical of a finitely generated ideal.
    Let $X$ be a $V$-flat finite type scheme and $x\in X$ a point lying over $s\in \Spec\, V$. Assume that the local ring $\sO_{X_s,x}$ is regular. Denote $A\ce \sO_{X,x}$.
    Then, for every finitely presented $A$-module $M$,
     \[
     \mathrm{proj.}\dim_A(M)+\depth_A(M)=\depth_A(A)=\dim (\sO_{X_s,x})+1.
     \]
     (By convention, $\mathrm{proj.}\dim_A(0)=-\infty$)
     \eenum
    \ethmt
     \bpf
       By \Cref{coh of V-flat ft alg}~\ref{coherence of O}--\ref{G-R 17.4.1 bound}, the local ring $\sO_{X,x}$ in (ii) is coherent regular. The assumption on $V$ implies that $\Spec \sO_{X,x}\backslash \{x\}$ is quasi-compact, so, by \Cref{depth of O_{X,x}}, $\sO_{X,x}$ has depth $\dim (\sO_{X_s,x})+1$. Therefore, it is enough to prove part (i).

        For (i), consider first the case where $\depth_A(A)=0$. We claim that every $A$-module $M$ having a finite resolution by finite free $A$-modules is free; thus in this case the formula in (i) holds. Clearly, it suffices to prove that every short exact sequence of the form $0\to A^{\oplus m}\xrightarrow{s} A^{\oplus n}\to N\to 0$ splits. Indeed, since $R^0\GG_{x}A\neq 0$ and the map $R^0\GG_{x}(s)$ is injective, the last statement follows from \Cref{a lem on linear algebra}.
       Here and in what follows, we redefine $x\ce \fm_A$.

       Assume now that $\depth_A(A)\ge 1$.
       Set $d\ce \depth_A(A)-1$.
        We will induct on $\mathrm{proj.}\dim_A(M)$ to verify the formula in (i).
      If $\mathrm{proj.}\dim_A(M)=0$, that is, $M$ is $A$-free, then it is clear that the formula holds.

     Next, assume that $\mathrm{proj.}\dim_A(M)\ge 1$, so every partial resolution $0\to M\pr  \xrightarrow{\iota} A^{\oplus n}\to M\to 0$ is non-split and satisfies $\mathrm{proj.}\dim_A(M')=\mathrm{proj.}\dim_A(M)-1$.
     It is a standard fact that $M'$ also has a finite resolution by finite free $A$-modules. We exploit the associated long exact sequence
\[
\cdots \to  R^{i-1}\GG_{x}M\pr \to R^{i-1}\GG_{x}A^{\oplus n} \to R^{i-1}\GG_{x} M \to R^i\GG_{x}M' \to R^i\GG_{x} A^{\oplus n} \to \cdots.
\]
     If $\mathrm{proj.}\dim_A(M)=1$, then $M\pr \simeq A^{\oplus m}$ for some $m\ge 1$. 
     {We have seen that $\depth_A(M\pr)=d+1$}, and so $R^{i}\GG_{x} M\pr=0$ for all $i\le d$. It follows from the above long exact sequence that $R^{i}\GG_{x} M=0$ for all $i\le d-1$. If $R^{d}\GG_{x}M $ were zero, then
    $$
    R^{d+1}\GG_x(\iota):\left(R^{d+1}\GG_{x} A\right)^{\oplus m} \simeq R^{d+1}\GG_{x}M' \hookrightarrow \left(R^{d+1}\GG_{x} A\right)^{\oplus n}$$
   is injective. Since $R^{d+1}\GG_{x} A $ is nonzero and supported on $\{x\}$, we deduce from \Cref{a lem on linear algebra} that $\iota$ splits, and so $M$ is $A$-free. This contradicts our assumption that $\mathrm{proj.}\dim_A(M)=1$. Therefore, $\depth_A(M)=d$, leading to the desired formula $\mathrm{proj.}\dim_A(M)+\depth_A(M)=d+1$.

     If $\mathrm{proj.}\dim_A(M)>1$, then $\mathrm{proj.}\dim_A(M')=\mathrm{proj.}\dim_A(M)-1$. Applying the induction hypothesis to $M\pr$, we
     obtain that 
     \begin{equation*}
     \begin{split}
     \depth_A(M\pr) & = d+1-(\mathrm{proj.}\dim_A(M)-1) \\
     & =d+2-\mathrm{proj.}\dim_A(M), \q \text{ which is } \q \le d.
     \end{split} 
     \end{equation*}
     It follows from the above long exact sequence that $R^{i-1}\GG_{x} M \simeq R^{i}\GG_{x} M\pr$ for all $i\le d$. Combining this with the bound $\depth_A(M\pr)\le d$, we deduce that $\depth_A(M)=\depth_A(M\pr)-1$. Therefore, by induction hypothesis, we have
     \[
     \mathrm{proj.}\dim_A(M)+\depth_A(M) =\left(\mathrm{proj.}\dim_A(M')+1\right)+\left(\depth_A(M\pr)-1\right)=d+1.
     \]
This finishes the induction step.\qedhere
     \epf
{The following corollary fails without the quasi-compactness assumption, see Example~\ref{exa-val-depth>=2}.}
\bcort \label{wdim=depth}
If $A$ is a regular coherent local ring with a quasi-compact punctured spectrum, then
\[
   \wdim(A)=\depth_{A}(A)<\infty.
\]
\ecort
\bpf
The Auslander--Buchsbaum formula implies that every finitely presented $A$-module has projective dimension at most $\depth_{A}(A)$. Since for any ring $R$, we have
$$
\wdim(R)=\sup\{\fldim(R/J) \,|\, \x{finitely generated ideal $J\subset R$}\},$$
 we obtain that $\wdim(A)\le \depth_{A}(A)$.
On the other hand, the quasi-compactness assumption implies that there is a finitely generated ideal $I\subset A$ contained in the maximal ideal $\fm_A$ such that $\sqrt{I}=\fm_A$.
Now we let $M\ce A/I$; it is a coherent $A$-module and therefore perfect with depth zero.
Applying the Auslander--Buchsbaum formula, we have $\pd_A(A/I)=\depth_A(A)$, which implies:
  \[
  \wdim(A)=\depth_A(A)=\pd_A(A/I)<\infty. \qedhere
   \]
\epf


\end{appendix}


\bib{Bou98}{book}{
  author={Bourbaki, Nicolas},
  title={Commutative algebra. {C}hapters 1--7},
  series={Elements of Mathematics (Berlin)},
  note={Translated from the French, Reprint of the 1989 English translation},
  publisher={Springer-Verlag, Berlin},
  year={1998},
  pages={xxiv+625},
  isbn={3-540-64239-0},
  label={BouAC},
}

\bib{EGAI}{article}{
  author={Grothendieck, A.},
  author={Dieudonn\'{e}, J.},
  title={\'El\'ements de g\'eom\'etrie alg\'ebrique. I. Le langage des sch\'emas},
  journal={Inst. Hautes \'Etudes Sci. Publ. Math.},
  number={4},
  date={1960},
  pages={228},
  issn={0073-8301},
  review={\MR {0217083 (36 \#177a)}},
  label={EGA~I},
}

\bib{EGAIV2}{article}{
  author={Grothendieck, Alexander},
  author={Dieudonn\'{e}, Jean},
  title={\'El\'ements de g\'eom\'etrie alg\'ebrique. IV. \'Etude locale des sch\'emas et des morphismes de sch\'emas. II},
  language={French},
  journal={Inst. Hautes \'Etudes Sci. Publ. Math.},
  number={24},
  date={1965},
  pages={231},
  issn={0073-8301},
  review={\MR {0199181 (33 \#7330)}},
  label={EGA~IV$_{2}$},
}

\bib{EGAIV3}{article}{
  author={Grothendieck, A.},
  author={Dieudonn\'{e}, J.},
  title={\'El\'ements de g\'eom\'etrie alg\'ebrique. IV. \'Etude locale des sch\'emas et des morphismes de sch\'emas. III},
  journal={Inst. Hautes \'Etudes Sci. Publ. Math.},
  number={28},
  date={1966},
  pages={255},
  issn={0073-8301},
  review={\MR {0217086 (36 \#178)}},
  label={EGA~IV$_{3}$},
}

\bib{EGAIV4}{article}{
  author={Grothendieck, Alexander},
  author={Dieudonn\'{e}, Jean Alexandre Eug\`ene},
  title={\'El\'ements de g\'eom\'etrie alg\'ebrique. IV. \'Etude locale des sch\'emas et des morphismes de sch\'emas IV},
  language={French},
  journal={Inst. Hautes \'Etudes Sci. Publ. Math.},
  number={32},
  date={1967},
  pages={361},
  issn={0073-8301},
  review={\MR {0238860 (39 \#220)}},
  label={EGA~IV$_{4}$},
}

\bib{SGA2new}{book}{
  author={Grothendieck, Alexander},
  title={Cohomologie locale des faisceaux coh\'erents et th\'eor\`emes de Lefschetz locaux et globaux (SGA 2)},
  language={French},
  series={Documents Math\'ematiques (Paris) [Mathematical Documents (Paris)]},
  volume={4},
  note={S\'eminaire de G\'eom\'etrie Alg\'ebrique du Bois Marie, 1962; Augment\'e d'un expos\'e de Mich\`ele Raynaud. [With an expos\'e by Mich\`ele Raynaud]; With a preface and edited by Yves Laszlo; Revised reprint of the 1968 French original},
  publisher={Soci\'et\'e Math\'ematique de France, Paris},
  date={2005},
  pages={x+208},
  isbn={2-85629-169-4},
  review={\MR {2171939}},
  label={SGA~2$_{\text {new}}$},
}

\bib{SGA3IIInew}{collection}{
  title={Sch\'emas en groupes (SGA 3). Tome III. Structure des sch\'emas en groupes r\'eductifs},
  language={French},
  series={Documents Math\'ematiques (Paris) [Mathematical Documents (Paris)], 8},
  editor={Gille, Philippe},
  editor={Polo, Patrick},
  note={S\'eminaire de G\'eom\'etrie Alg\'ebrique du Bois Marie 1962--64. [Algebraic Geometry Seminar of Bois Marie 1962--64]; A seminar directed by M. Demazure and A. Grothendieck with the collaboration of M. Artin, J.-E. Bertin, P. Gabriel, M. Raynaud and J-P. Serre; Revised and annotated edition of the 1970 French original},
  publisher={Soci\'et\'e Math\'ematique de France, Paris},
  date={2011},
  pages={lvi+337},
  isbn={978-2-85629-324-9},
  review={\MR {2867622}},
  label={SGA~3$_{\text {III new}}$},
}

\bib{SGA4II}{book}{
  title={Th\'eorie des topos et cohomologie \'etale des sch\'emas. Tome 2},
  language={French},
  series={Lecture Notes in Mathematics, Vol. 270},
  note={S\'eminaire de G\'eom\'etrie Alg\'ebrique du Bois-Marie 1963--1964 (SGA 4); Dirig\'e par M. Artin, A. Grothendieck et J. L. Verdier. Avec la collaboration de N. Bourbaki, P. Deligne et B. Saint-Donat},
  publisher={Springer-Verlag},
  place={Berlin},
  date={1972},
  pages={iv+418},
  review={\MR {0354653 (50 \#7131)}},
  label={SGA~4$_{\text {II}}$},
}

\bib{Alp14}{article}{
  author={Alper, Jarod},
  title={Adequate moduli spaces and geometrically reductive group schemes},
  journal={Algebr. Geom.},
  volume={1},
  year={2014},
  number={4},
  pages={489--531},
  issn={2313-1691},
  doi={10.14231/AG-2014-022},
  url={https://doi-org.revues.math.u-psud.fr/10.14231/AG-2014-022},
}

\bib{AB57}{article}{
  author={Auslander, Maurice},
  author={Buchsbaum, David A.},
  title={Homological dimension in local rings},
  journal={Trans. Am. Math. Soc.},
  issn={0002-9947},
  volume={85},
  pages={390--405},
  year={1957},
  language={English},
  doi={10.2307/1992937},
}

\bib{Aus62}{article}{
  author={Auslander, M.},
  title={On the purity of the branch locus},
  journal={Am. J. Math.},
  issn={0002-9327},
  volume={84},
  pages={116--125},
  year={1962},
  language={English},
  doi={10.2307/2372807},
}

\bib{Ber71}{article}{
  author={Bertin, Jos{\'e}},
  title={Anneaux coh{\'e}rents r{\'e}guliers. ({Regular} coherent rings)},
  journal={C. R. Acad. Sci., Paris, S{\'e}r. A},
  issn={0366-6034},
  volume={273},
  pages={590--591},
  year={1971},
  language={French},
  keywords={13F15,13C15,13E15},
}

\bib{Ber72}{misc}{
  author={Bertin, Jose},
  title={Anneaux coh{\'e}rents de dimension homologique finie},
  year={1972},
  language={French},
  howpublished={Colloque {Algebre} commut., {Rennes} 1972, {Publ}. {Sem}. {Math}. {Univ}. {Rennes}, {No}. 16, 12 p. (1972).},
  keywords={13H05,13C10,13D05},
}

\bib{BS15}{incollection}{
  author={Bhatt, Bhargav},
  author={Scholze, Peter},
  title={The pro-{\'e}tale topology for schemes},
  booktitle={Astérisque (volume~369), De la g\'eom\'etrie alg\'ebrique aux formes automorphes (I). Une collection d'articles en l'honneur du soixanti\`eme anniversaire de G\'erard Laumon},
  isbn={978-2-85629-805-3},
  pages={99--201},
  year={2015},
  publisher={Paris: Soci{\'e}t{\'e} Math{\'e}matique de France (SMF)},
  language={English},
  keywords={19F27,18F10,14H30},
}

\bib{BR83}{article}{
  author={Bhatwadekar, S. M. },
  author={Rao, R. A.},
  title={On a question of {Quillen}},
  journal={Trans. Am. Math. Soc.},
  issn={0002-9947},
  volume={279},
  pages={801--810},
  year={1983},
  language={English},
  doi={10.2307/1999568},
  keywords={13C10,18F25,13D15,13H05},
}

\bib{BC22}{article}{
  title={Torsors on loop groups and the Hitchin fibration},
  author={Bouthier, Alexis},
  author={{\v {C}}esnavi{\v {c}}ius, K{\k {e}}stutis},
  journal={Annales scientifiques de l'École normale supérieure},
  number={3},
  pages={791–864},
  volume={55},
  year={2022},
  doi={10.24033/asens.2506},
}

\bib{Ces22a}{article}{
  title={Grothendieck--Serre in the quasi-split unramified case},
  author={{\v {C}}esnavi{\v {c}}ius, K{\k {e}}stutis},
  journal={Forum Math. Pi},
  year={2022},
  number={e9, 30},
  volume={10},
  doi={10.1017/fmp.2022.5},
}

\bib{Ces22}{article}{
  author={{\v {C}}esnavi{\v {c}}ius, K{\k {e}}stutis},
  title={Problems about torsors over regular rings},
  journal={Acta Math. Vietnam.},
  issn={0251-4184},
  volume={47},
  number={1},
  pages={39--107},
  year={2022},
  language={English},
  doi={10.1007/s40306-022-00477-y},
  keywords={14M17,14L15,20G10},
}

\bib{Ces22b}{article}{
  title={Torsors on the complement of a smooth divisor},
  author={{\v {C}}esnavi{\v {c}}ius, K{\k {e}}stutis},
  journal={Available at \url {https://www.imo.universite-paris-saclay.fr/~cesnavicius/torsors-complement.pdf}},
  year={2022},
}

\bib{CS21}{article}{
  author={{\v {C}}esnavi{\v {c}}ius, K{\k {e}}stutis},
  author={Scholze, Peter},
  title={Purity for flat cohomology},
  journal={Ann. Math. (2)},
  issn={0003-486X},
  volume={199},
  number={1},
  pages={51--180},
  year={2024},
  language={English},
  doi={10.4007/annals.2024.199.1.2},
  keywords={14F20,14F22,14F30,14H20},
}

\bib{CTS79}{article}{
  author={Colliot-Th\'el\`ene, Jean-Louis},
  author={Sansuc, Jean-Jacques},
  title={Fibr\'{e}s quadratiques et composantes connexes r\'{e}elles},
  journal={Math. Ann.},
  fjournal={Mathematische Annalen},
  volume={244},
  year={1979},
  number={2},
  pages={105--134},
  issn={0025-5831},
  doi={10.1007/BF01420486},
  url={https://doi.org/10.1007/BF01420486},
}

\bib{Fed21}{article}{
  title={On the purity conjecture of Nisnevich for torsors under reductive group schemes},
  author={Fedorov, Roman},
  journal={Available at \url {https://arxiv.org/pdf/2109.10332.pdf}},
  year={2021},
}

\bib{Fed22b}{article}{
  author={Fedorov, Roman},
  title={On the {Grothendieck}-{Serre} conjecture on principal bundles in mixed characteristic},
  journal={Trans. Am. Math. Soc.},
  issn={0002-9947},
  volume={375},
  number={1},
  pages={559--586},
  year={2022},
  language={English},
  doi={10.1090/tran/8490},
  keywords={14L15,14D10,14M17,20G35,20G41,11E08},
}

\bib{Gab81}{misc}{
  author={Gabber, Ofer},
  title={Some theorems on {Azumaya} algebras},
  year={1981},
  language={English},
  howpublished={Groupe de {Brauer}, {S{\'e}min}., {Les} {Plans}-sur-{Bex} 1980, {Lect}. {Notes} {Math}. 844, 129-209 (1981).},
  doi={10.1007/BFb0090480},
  keywords={14F05,14F22,16H05,14F20},
}

\bib{GR18}{article}{
  title={Foundations for almost ring theory},
  author={Gabber, Ofer},
  author={Ramero, Lorenzo},
  journal={Available at \url {https://arxiv.org/abs/math/0409584}},
  year={2018},
}

\bib{Gla89}{book}{
  author={Glaz, Sarah},
  title={Commutative coherent rings},
  series={Lect. Notes Math.},
  issn={0075-8434},
  volume={1371},
  isbn={3-540-51115-6},
  year={1989},
  publisher={Berlin etc.: Springer-Verlag},
  language={English},
  doi={10.1007/BFb0084570},
  keywords={13E15,13-01,13-02,13D05,13B02},
}

\bib{Gil21}{article}{
  author={Gille, Philippe},
  title={When is a reductive group scheme linear?},
  journal={Mich. Math. J.},
  issn={0026-2285},
  volume={72},
  pages={439--448},
  year={2022},
  language={English},
  doi={10.1307/mmj/20217208},
  keywords={14L15,20G35},
  url={projecteuclid.org/journals/michigan-mathematical-journal/volume-72/issue-none/When-Is-a-Reductive-Group-Scheme-Linear/10.1307/mmj/20217208.full},
}

\bib{Guo20b}{article}{
  title={The Grothendieck--Serre conjecture over valuation rings},
  author={Guo, Ning},
  journal={Compos. Math.},
  year={2024},
  volume={160},
  number={2},
  pages={317--355},
  doi={10.1112/S0010437X23007583},
  language={English},
  publisher={London Mathematical Society},
}

\bib{GL23}{article}{
  title={Grothendieck--Serre for constant reductive group schemes},
  author={Guo, Ning},
  author={Liu, Fei},
  journal={Available at \url {https://arxiv.org/abs/2301.12460}},
  year={2023},
}

\bib{Kna03}{incollection}{
  author={Knaf, Hagen},
  title={Regular curves over {Pr{\"u}fer} domains},
  booktitle={Valuation theory and its applications. Volume II. Proceedings of the international conference and workshop, University of Saskatchewan, Saskatoon, Canada, July 28--August 11, 1999},
  isbn={0-8218-3206-9},
  pages={89--106},
  year={2003},
  publisher={Providence, RI: American Mathematical Society (AMS)},
  language={English},
  keywords={13F05,13D05},
}

\bib{Kna08}{article}{
  author={Knaf, Hagen},
  title={Regular local algebras over valuation domains: weak dimension and regular sequences},
  journal={J. Algebra Appl.},
  issn={0219-4988},
  volume={7},
  number={5},
  pages={575--591},
  year={2008},
  language={English},
  doi={10.1142/S0219498808003004},
  keywords={13D05,13H05},
}

\bib{KM76}{article}{
  author={Knudsen, Finn},
  author={Mumford, David},
  title={The projectivity of the moduli space of stable curves. {I}: {Preliminaries} on ``det'' and ``{Div}''},
  journal={Math. Scand.},
  issn={0025-5521},
  volume={39},
  pages={19--55},
  year={1976},
  language={English},
  doi={10.7146/math.scand.a-11642},
  keywords={14H10,14D20},
}

\bib{Mar16}{article}{
  title={A purity theorem for torsors},
  author={Marrama, Andrea},
  journal={ALGANT master thesis},
  year={2016},
}

\bib{MB22}{article}{
  author={Moret-Bailly, Laurent},
  title={A construction of weakly unramified extensions of a valuation ring},
  journal={Rend. Semin. Mat. Univ. Padova},
  issn={0041-8994},
  volume={147},
  pages={139--151},
  year={2022},
  language={French},
  doi={10.4171/RSMUP/94},
  keywords={13F30,12J10,12J25},
}

\bib{Nag59}{article}{
  author={Nagata, Masayoshi},
  title={On the purity of branch loci in regular local rings},
  journal={Ill. J. Math.},
  issn={0019-2082},
  volume={3},
  pages={328--333},
  year={1959},
  language={English},
}

\bib{Nag66}{article}{
  author={Nagata, Masayoshi},
  title={Finitely generated rings over a valuation ring},
  journal={J. Math. Kyoto Univ.},
  fjournal={Journal of Mathematics of Kyoto University},
  volume={5},
  year={1966},
  pages={163--169},
  issn={0023-608X},
  doi={10.1215/kjm/1250524533},
  url={https://doi-org.revues.math.u-psud.fr/10.1215/kjm/1250524533},
}

\bib{Nis89}{article}{
  author={Nisnevich, Yevsey A.},
  title={Rationally trivial principal homogeneous spaces, purity and arithmetic of reductive group schemes over extensions of two-dimensional regular local rings},
  language={English, with French summary},
  journal={C. R. Acad. Sci. Paris S\'er. I Math.},
  volume={309},
  date={1989},
  number={10},
  pages={651--655},
  issn={0764-4442},
  review={\MR {1054270}},
}

\bib{Pan18}{article}{
  author={Panin, Ivan},
  title={On Grothendieck-Serre conjecture concerning principal bundles},
  conference={ title={Proceedings of the International Congress of Mathematicians---Rio de Janeiro 2018. Vol. II. Invited lectures}, },
  book={ publisher={World Sci. Publ., Hackensack, NJ}, },
  date={2018},
  pages={201--221},
  review={\MR {3966763}},
}

\bib{Pop02}{article}{
  author={Popescu, Dorin},
  title={On a question of {Quillen}},
  journal={Bull. Math. Soc. Sci. Math. Roum., Nouv. S{\'e}r.},
  issn={1220-3874},
  volume={45},
  number={3-4},
  pages={209--212},
  year={2002},
  language={English},
  keywords={13C10,13D15,13H05},
}

\bib{Que71}{article}{
  author={Quentel, Yann},
  title={Sur le th{\'e}or{\`e}me d'{Auslander}-{Buchsbaum}. ({On} the {Auslander} {Buchsbaum} theorem.)},
  journal={C. R. Acad. Sci., Paris, S{\'e}r. A},
  issn={0366-6034},
  volume={273},
  pages={880--881},
  year={1971},
  language={French},
  keywords={13H05,13F15},
}

\bib{RG71}{article}{
  author={Raynaud, Michel},
  author={Gruson, Laurent},
  title={Crit{\`e}res de platitude et de projectivit{\'e}. {Techniques} de ''platification'' d'un module. ({Criterial} of flatness and projectivity. {Technics} of ''flatification of a module.)},
  journal={Invent. Math.},
  issn={0020-9910},
  volume={13},
  pages={1--89},
  year={1971},
  language={French},
  doi={10.1007/BF01390094},
  keywords={14B25,13C11,14F20,13J15,14F10,14F05},
}

\bib{Sam64}{article}{
  author={Samuel, Pierre},
  title={Anneaux gradu{\'e}s factoriels et modules r{\'e}flexifs},
  journal={Bull. Soc. Math. Fr.},
  issn={0037-9484},
  volume={92},
  pages={237--249},
  year={1964},
  language={French},
  doi={10.24033/bsmf.1608},
}

\bib{Ser56}{misc}{
  author={Serre, Jean-Pierre},
  title={Sur la dimension homologique des anneaux et des modules noeth{\'e}riens},
  year={1956},
  language={French},
  howpublished={Proc. {Int}. {Symp}. {Algebraic} {Number} {Theory} 1955, 175-189 (1956).},
  keywords={13D05,13E05,13H05},
}

\bib{SP}{misc}{
  shorthand={Stacks},
  author={The {Stacks Project Authors}},
  title={\textit {Stacks Project}},
  howpublished={\url {https://stacks.math.columbia.edu}},
  year={2018},
  label={SP},
}

\bib{Zaf78}{article}{
  author={Zafrullah, Muhammad},
  title={On finite conductor domains},
  journal={Manuscr. Math.},
  issn={0025-2611},
  volume={24},
  pages={191--204},
  year={1978},
  language={English},
  doi={10.1007/BF01310053},
}

\bib{Zar58}{article}{
  author={Zariski, Oscar},
  title={On the purity of the branch locus of algebraic functions},
  journal={Proc. Natl. Acad. Sci. USA},
  issn={0027-8424},
  volume={44},
  pages={791--796},
  year={1958},
  language={English},
  doi={10.1073/pnas.44.8.791},
}



\begin{bibdiv}
\begin{biblist}
\bibselect{bibliography}
\end{biblist}
\end{bibdiv}

\end{document}